\documentclass{amsart}

\usepackage[arrow,curve,matrix,tips,2cell]{xy}

\usepackage[alphabetic]{amsrefs}
\usepackage{amsmath}
\usepackage{amscd}
\usepackage{amssymb}
\usepackage{latexsym}
\usepackage{stmaryrd}

\newtheorem{theorem}{Theorem}[section]
\newtheorem*{theorem*}{Theorem}
\newtheorem{lemma}[theorem]{Lemma}
\newtheorem*{lemma*}{Lemma}
\newtheorem{corollary}[theorem]{Corollary}
\newtheorem*{corollary*}{Corollary}
\newtheorem{proposition}[theorem]{Proposition}

\newtheorem{remark}[theorem]{Remark}
\newtheorem{question}[theorem]{Question}
\newtheorem{definition}[theorem]{Definition}

\newtheorem{example}[theorem]{Example}
\newtheorem{examples}[theorem]{Examples}
\newcommand{\bgl}{\begin{equation}} 
\newcommand{\egl}{\end{equation}}
\newcommand{\bgloz}{\begin{equation*}} 
\newcommand{\egloz}{\end{equation*}}
\newcommand{\bgln}{\begin{eqnarray}} 
\newcommand{\egln}{\end{eqnarray}}
\newcommand{\bglnoz}{\begin{eqnarray*}} 
\newcommand{\eglnoz}{\end{eqnarray*}}
\newcommand{\btheo}{\begin{theorem}}
\newcommand{\etheo}{\end{theorem}}
\newcommand{\btheooz}{\begin{theorem*}}
\newcommand{\etheooz}{\end{theorem*}}
\newcommand{\blemma}{\begin{lemma}}
\newcommand{\elemma}{\end{lemma}}
\newcommand{\blemmaoz}{\begin{lemma*}}
\newcommand{\elemmaoz}{\end{lemma*}}
\newcommand{\bproof}{\begin{proof}}
\newcommand{\eproof}{\end{proof}}
\newcommand{\bbew}{\begin{beweis}}
\newcommand{\ebew}{\end{beweis}}
\newcommand{\bremark}{\begin{remark}\em}
\newcommand{\eremark}{\end{remark}}
\newcommand{\bquestion}{\begin{question}\em}
\newcommand{\equestion}{\end{question}}
\newcommand{\bdefin}{\begin{definition}}
\newcommand{\edefin}{\end{definition}}
\newcommand{\bprop}{\begin{proposition}}
\newcommand{\eprop}{\end{proposition}}
\newcommand{\bcor}{\begin{corollary}}
\newcommand{\ecor}{\end{corollary}}
\newcommand{\bcoroz}{\begin{corollary*}}
\newcommand{\ecoroz}{\end{corollary*}}
\newcommand{\bfa}{\begin{cases}} 
\newcommand{\efa}{\end{cases}}
\newcommand{\bexample}{\begin{example}\em}
\newcommand{\eexample}{\end{example}}
\newcommand{\bexamples}{\begin{examples}\em}
\newcommand{\eexamples}{\end{examples}}
\newcommand{\cB}{\mathcal B}
\newcommand{\cC}{\mathcal C}
\newcommand{\cD}{\mathcal D}

\newcommand{\cF}{\mathcal F}
\newcommand{\cG}{\mathcal G}

\newcommand{\cI}{\mathcal I}
\newcommand{\cJ}{\mathcal J}
\newcommand{\cK}{\mathcal K}
\newcommand{\cL}{\mathcal L}

\newcommand{\cO}{\mathcal O}
\newcommand{\cP}{\mathcal P}

\newcommand{\cS}{\mathcal S}

\def\Cz{\mathbb{C}}
\def\Fz{\mathbb{F}}

\def\Nz{\mathbb{N}}

\def\Qz{\mathbb{Q}}
\def\Rz{\mathbb{R}}

\def\Zz{\mathbb{Z}}

\def\1z{\mathbb{1}}

\newcommand{\fA}{\mathfrak A}

\newcommand{\fX}{\mathfrak X}
\newcommand{\mfa}{\mathfrak a}
\newcommand{\mfb}{\mathfrak b}

\newcommand{\mfp}{\mathfrak p}
\newcommand{\bI}{\mathbf{I}}
\newcommand{\an}[1]{``#1''} 
\newcommand{\ti}{\tilde}

\newcommand{\lori}{\longrightarrow}

\newcommand{\ma}{\mapsto} 
\newcommand\onto{\twoheadrightarrow} 
\newcommand\into{\hookrightarrow} 
\newcommand{\Rarr}{\Rightarrow} 
\newcommand{\Larr}{\Leftarrow} 
\newcommand{\LRarr}{\Leftrightarrow} 

\def\SEMI{\mbox{$\times\kern-2pt\vrule height5pt width.6pt \kern3pt $}}

\newcommand{\halb}{\tfrac{1}{2}}

\newcommand{\img}{{\rm im\,}}

\newcommand{\Spec}{{\rm Spec\,}} 

\newcommand{\id}{{\rm id}}

\newcommand{\lcm}{{\rm lcm}} 
\newcommand{\reg}{^\times} 
\newcommand{\lspan}{{\rm span}} 
\newcommand{\clspan}{\overline{\lspan}} 
\newcommand{\abs}[1]{\left|#1\right|} 
\newcommand{\norm}[1]{\left\|#1\right\|} 
\newcommand{\defeq}{\mathrel{:=}} 
\newcommand{\eqdef}{\mathrel{=:}} 

\newcommand{\dop}{\text{: }} 
\newcommand{\falls}{\text{ if }} 
\newcommand{\fa}{\text{ for all }} 

\newcommand{\dom}{{\rm dom}}

\newcommand{\lge}{\left\{} 
\newcommand{\rge}{\right\}} 
\newcommand{\lru}{\left(} 
\newcommand{\rru}{\right)} 
\newcommand{\leck}{\left[} 
\newcommand{\reck}{\right]} 
\newcommand{\lsp}{\left\langle} 
\newcommand{\rsp}{\right\rangle} 
\newcommand{\rukl}[1]{\lru #1 \rru} 
\newcommand{\eckl}[1]{\leck #1 \reck} 
\newcommand{\gekl}[1]{\lge #1 \rge} 
\newcommand{\spkl}[1]{\lsp #1 \rsp} 

\newcommand{\height}{\operatorname{ht}} 

\newcommand{\menge}[2]{\gekl{ #1 \, \dop #2 }} 
\begin{document}

\title{Semigroup C*-algebras}

\author{Xin Li}

\begin{abstract}
We give an overview of some recent developments in semigroup C*-algebras.
\end{abstract}

\maketitle 

\tableofcontents

\setlength{\parindent}{0cm} \setlength{\parskip}{0.5cm}

\section{Introduction}

A semigroup C*-algebra is the C*-algebra generated by the left regular representation of a left cancellative semigroup. In the case of groups, this is the classical construction of reduced group C*-algebras, which received great interest and serves as a motivating class of examples in operator algebras.

For semigroups which are far from being groups, we encounter completely new phenomena which are not visible in the group case. It is therefore a natural and interesting task to try to understand and explain these new phenomena. This challenge has been taken up by several authors in many pieces of work, and our present goal is to give a unified treatment of this endeavour.

We point out that particular classes of semigroups played a predominant role in the development, as they serve as our motivation and guide us towards important properties of semigroups which allow for a systematic study of their C*-algebras. The examples include positive cones in totally ordered groups, semigroups given by particular presentations and semigroups coming from rings of number-theoretic origin. Important properties that isolate from the general and wild class of all left cancellative semigroups a manageable subclass were first given by Nica's quasi-lattice order \cite{Nica} and later on by the independence condition \cite{LiSG} and the Toeplitz condition \cite{LiNuc}.

Aspects of semigroup C*-algebras which we would like to discuss in the following include descriptions as crossed products and groupoid C*-algebras, the connection between amenability and nuclearity, boundary quotients, and the classification problem for semigroup C*-algebras. The first three topics are discussed in detail, and we give a more or less self-contained presentation. The last topic puts together many results. In particular, it builds on the K-theory computations that are explained in detail by S. Echterhoff in \cite{EchKK}. Since a detailed account of classification results would take too much space, we just briefly summarize the main results, and refer the interested reader to the relevant papers for more details and complete proofs.

Our discussion of semigroup C*-algebras builds on previous work of J. Renault on groupoids and their C*-algebras \cite{Ren1}, and the work of R. Exel on C*-algebras of inverse semigroups, their quotients corresponding to tight representations of inverse semigroups, and on partial actions \cites{Ex2,Ex3,Ex4}. 

Inevitably, certain interesting aspects of semigroup C*-algebras are not covered in this book. This includes a discussion of C*-algebras of semigroups which do not embed into groups such as general right LCM semigroups (see \cite{Star}) or Zappa-Sz\'{e}p products (see \cite{BRRW}), or C*-algebras of certain topological semigroups (see \cites{ReSu,Sun}). Moreover, we do not discuss KMS-states in detail, but we refer the reader to \cites{LR2,BHLR,CDL,CHR} for more information. We also mention that in \cite{CunAlgAct}, J. Cuntz describes KMS-states for particular examples. We apologize for these omissions and try to make up for them by pointing the interested reader to the relevant literature. To this end, we have included a long (but not complete) list of references.

\section{C*-algebras generated by left regular representations}
\label{leftregrep}

Let $P$ be a semigroup. We assume that $P$ is left cancellative, i.e., for all $p, x, y \in P$, $px = py$ implies $x = y$. In other words, the map
$$
  P \to P, \, x \ma px
$$
given by left multiplication with $p \in P$ is injective for all $p \in P$.

The left regular representation of $P$ is given as follows: The Hilbert space $\ell^2 P$ comes with a canonical orthonormal basis $\menge{\delta_x}{x \in P}$. Here $\delta_x$ is the delta-function in $x \in P$, defined by
$$
  \delta_x(y) = 1 \ {\rm if} \ y = x \ {\rm and} \ 
  \delta_x(y) = 0 \ {\rm if} \ y \neq x.
$$
For every $p \in P$, the map
$$
  P \to P, \, x \ma px
$$
is injective by left cancellation, so that the mapping
$$
  \delta_x \ma \delta_{px} \ (x \in P)
$$
extends (uniquely) to an isometry
$$
  V_p : \: \ell^2 P \to \ell^2 P.
$$
The assignment
$$
  p \ma V_p \ (p \in P)
$$
represents our semigroup $P$ as isometries on $\ell^2 P$. This is called the left regular representation of $P$. It generates the following C*-algebra:

\bdefin
\label{sgpCSTAR}
$$
  C^*_{\lambda}(P) \defeq C^*(\menge{V_p}{p \in P}) \subseteq \cL(\ell^2 P).
$$
\edefin

By definition, $C^*_{\lambda}(P)$ is the smallest subalgebra of $\cL(\ell^2 P)$ containing $\menge{V_p}{p \in P}$ which is invariant under forming adjoints and closed in the operator norm topology. We call $C^*_{\lambda}(P)$ the semigroup C*-algebra of $P$, or more precisely, the left reduced semigroup C*-algebra of $P$.

Note that left cancellation is a crucial assumption for our construction. In general, without left cancellation, the mapping $\delta_x \ma \delta_{px}$ does not even extend to a bounded linear operator on $\ell^2 P$. Moreover, we point out that we view our semigroups as discrete objects. Our construction, and some of the analysis, carries over to certain topological semigroups (see \cites{ReSu,Sun}). Finally, $C^*_{\lambda}(P)$ will be separable if $P$ is countable. This helps to exclude pathological cases. Therefore, for convenience, we assume from now on that all our semigroups are countable, although this is not always necessary in our discussion.

\section{Examples}
\label{sec:Ex}

We have already pointed out the importance of examples. Therefore, it is appropriate to start with a list of examples of semigroups where we can apply our construction. All our examples are actually semigroups with an identity, so that they are all monoids.

\subsection{The natural numbers}

Our first example is given by $P = \Nz = \gekl{0,1,2,\dotsc}$, the set of natural numbers including zero, viewed as an additive monoid. By construction, $V_1$ is the unilateral shift. Since $\Nz$ is generated by $1$ as a monoid, it is clear that $C^*_{\lambda}(\Nz)$ is generated as a C*-algebra by the unilateral shift. This C*-algebra has been studied by Coburn (see \cites{Cob1,Cob2}). It turns out that it is the universal C*-algebra generated by one isometry, i.e.,
$$
  C^*_{\lambda}(\Nz) \cong C^*(v \ \vert \ v^* v = 1), \, V_1 \ma v.
$$

$C^*_{\lambda}(\Nz)$ is also called the Toeplitz algebra. This name comes from the observation that $C^*_{\lambda}(\Nz)$ can also be described as the C*-algebra of Toeplitz operators on the Hardy space, defined on the circle. This interpretation connects our semigroup C*-algebra $C^*_{\lambda}(\Nz)$ with index theory and K-theory.

\subsection{Positive cones in totally ordered groups}
\label{ex:totalorder}

Motivated by connections to index theory and K-theory, several authors including Coburn and Douglas studied the following examples in \cites{CD,CDSS,Dou,DH}:

Let $G$ be a subgroup of $(\Rz,+)$, and consider the additive monoid $P = [0,\infty) \cap G$. The case $G = \Zz$ gives our previous example $P = \Nz$. The case where $G = \Zz[\lambda,\lambda^{-1}]$ for some positive real number $\lambda$ is discussed in \cites{CPPR,Li4}.

These examples belong to the bigger class of positive cones in totally ordered groups. A left invariant total order on a group $G$ is a relation $\leq$ on $G$ such that
\begin{itemize}
\item For all $x, y \in G$, we have $x=y$ if and only if $x \leq y$ and $y \leq x$.
\item For all $x, y \in G$, we always have $x \leq y$ or $y \leq x$.
\item For all $x, y, z \in G$, $x \leq y$ and $y \leq z$ imply $x \leq z$.
\item For all $x, y, z \in G$, $x \leq y$ implies $zx \leq zy$.
\end{itemize}

Given a left invariant total order $\leq$ on $G$, define $P \defeq \menge{x \in G}{e \leq x}$. Here $e$ is the identity in $G$. $P$ is called the positive cone in $G$. It is a monoid satisfying
\bgl
\label{G=PuP}
G = P \cup P^{-1} \ {\rm and} \ P \cap P^{-1} = \gekl{e}.
\egl
Conversely, every submonoid $P \subseteq G$ of a group $G$ satisfying \eqref{G=PuP} gives rise to a left invariant total order $\leq$ by setting, for $x, y \in G$, $x \leq y$ if $y \in xP$. Here $x P = \menge{x p}{p \in P} \subseteq G$.

In the examples mentioned above of subgroups of $(\Rz,+)$, we have canonical left invariant total orders given by restricting the canonical order on $(\Rz,+)$.

The study of left invariant total orders on group is of great interest in group theory. For instance, the existence of a left invariant total order on a group $G$ implies the Kaplansky conjecture for $G$. This conjecture says that for a torsion-free group $G$ and a ring $R$, the group ring $RG$ does not have zero-divisors if $R$ does not have zero-divisors. We refer to \cites{MuRh,DNR} for more details.

While it is known that every torsion-free nilpotent group admits a left invariant total order, it is an open conjecture that lattices in simple Lie groups of rank at least two have no left invariant total order. It is also an open question whether an infinite property (T) group can admit a left invariant total order (see \cite{DNR} for more details).

\subsection{Monoids given by presentations}
\label{ex:presentations}

Another source for examples of monoids comes from group presentations. One way to define a group is to give a presentation, i.e., generators and relations. For instance, the additive group of integers is the group generated by one element with no relation, $\Zz = \spkl{a}$. The non-abelian free group on two generators is the group generated by two elements with no relations, $\Fz_2 = \Zz * \Zz = \spkl{a,b}$. And $\Zz \times \Zz$ is the group generated by two elements which commute, $\Zz^2 = \Zz \times \Zz = \spkl{a,b \ \vert \ ab = ba}$. If we look at the semigroups (or rather monoids) defined by the same presentations, we get $\Nz = \spkl{a}^+$, $\Nz^{*2} = \Nz * \Nz = \spkl{a,b}^+$, $\Nz^2 = \Nz \times \Nz = \spkl{a,b \ \vert \ ab = ba}^+$. Here, we write $\spkl{\cdot \ \vert \ \cdot}^+$ for the universal monid given by a particular presentation, while we write $\spkl{\cdot \ \vert \ \cdot}$ for the universal group given by a particular representation. This is to distinguish between group presentations and monoid presentations.

Of course, in general, it is not clear whether this procedure of taking generators and relations from group presentations to define monoids leads to interesting semigroups, or whether we can apply our C*-algebraic construction to the resulting semigroups. For instance, it could be that the monoid given by a presentation actually coincides with the group given by the same presentation. Another problem that might arise is that the canonical homomorphism from the monoid to the group given by the same presentation, sending generator to generator, is not injective. In that case, our monoid might not even be left cancellative. However, there are conditions on our presentations which ensure that these problems do not appear. There is for instance the notion of completeness (see \cite{Deh}), explained in \S~\ref{independence}. Now let us just give a list of examples.

The presentations for $\Zz$, $\Fz_2$ and $\Zz^2$ all have in common that two generators either commute or satisfy no relation (i.e., they are free), and these are the only relations we impose. This can be generalized. Let $\Gamma = (V,E)$ be an undirected graph, where we connect two vertices by at most one edge and no vertex to itself. This means that we can think of $E$ as a subset of $V \times V$.

We then define
$$A_{\Gamma} \defeq \spkl{\menge{\sigma_v}{v \in V} \ \vert \ \sigma_v \sigma_w = \sigma_w \sigma_v \ {\rm for \ all} \ (v,w) \in E},$$
$$A_{\Gamma}^+ \defeq \spkl{\menge{\sigma_v}{v \in V} \ \vert \ \sigma_v \sigma_w = \sigma_w \sigma_v \ {\rm for \ all} \ (v,w) \in E}^+.$$

For instance, the graph for $\Zz$ only consists of one vertex and no edge, the graph for $\Fz_2$ consists of two vertices and no edges, and the graph for $\Zz^2$ consists of two vertices and one edge joining them.

The groups $A_{\Gamma}$ are called right-angled Artin groups and the monoids $A_{\Gamma}^+$ are called right-angled Artin monoids. Their C*-algebras are discussed in \cites{CrLa1,CrLa2,Iva,ELR}.

Right-angled Artin monoids and the corresponding groups are special cases of graph products. Let $\Gamma = (V,E)$ be a graph as above, with $E \subseteq V \times V$. Assume that for every $v \in V$, $G_v$ is a group containing a submonoid $P_v$. Then let $\Gamma_{v \in V} G_v$ be the group obtained from the free product $*_{v \in V} G_v$ by introducing the relations $xy = yx$ for all $x \in G_v$ and $y \in G_w$ with $(v,w) \in E$. Similarly, define $\Gamma_{v \in V} P_v$ as the monoid obtained from the free product $*_{v \in V} P_v$ by introducing the relations $xy = yx$ for all $x \in P_v$ and $y \in P_w$ with $(v,w) \in E$. It is explained in \cite{CrLa1} (see also \cites{GreenPhD,HM}) that the embeddings $P_v \into G_v$ induce an embedding
$$
  \Gamma_{v \in V} P_v \into \Gamma_{v \in V} G_v.
$$
In the case that $P_v \subseteq G_v$ is given by $\Nz \subseteq \Zz$ for all $v \in V$, we obtain right-angled Artin monoids and the corresponding groups.

We will have more to say about general graph products in \S~\ref{ss:GraphProducts} and \S~\ref{sec:GraphProducts}.

As the name suggests, there is a more general class of Artin groups which contains right-angled Artin groups. Let $I$ be a countable index set,
$$\menge{m_{ij} \in \gekl{2, 3, 4, \dotsc} \cup \gekl{\infty}}{i,j \in I, \, i \neq j}$$ be such that $m_{ij} = m_{ji}$ for all $i$ and $j$. Then define
$$G \defeq \spkl{\menge{\sigma_i}{i \in I} \ \vert \ \underbrace{\sigma_i \sigma_j \sigma_i \sigma_j \dotsm}_{m_{ij}} = \underbrace{\sigma_j \sigma_i \sigma_j \sigma_i \dotsm}_{m_{ji}} \ {\rm for \ all} \ i,j \in I, \, i \neq j}.$$
For $m_{ij} = \infty$, there is no relation involving $\sigma_i$ and $\sigma_j$, i.e., $\sigma_i$ and $\sigma_j$ are free. And define
$$P \defeq \spkl{\menge{\sigma_i}{i \in I} \ \vert \ \underbrace{\sigma_i \sigma_j \sigma_i \sigma_j \dotsm}_{m_{ij}} = \underbrace{\sigma_j \sigma_i \sigma_j \sigma_i \dotsm}_{m_{ji}} \ {\rm for \ all} \ i,j \in I, \, i \neq j}^+.$$
If $m_{ij} \in \gekl{2,\infty}$ for all $i$ and $j$, then we get right-angled Artin groups and monoids. To see some other groups, take for instance $I = \gekl{1,2}$ and $m_{1,2} = m_{2,1} = 3$. We get the (third) Braid group and the corresponding Braid monoid
$$B_3 \defeq \spkl{\sigma_1, \sigma_2 \ \vert \ \sigma_1 \sigma_2 \sigma_1 = \sigma_2 \sigma_1 \sigma_2},$$
$$B_3^+ \defeq \spkl{\sigma_1, \sigma_2 \ \vert \ \sigma_1 \sigma_2 \sigma_1 = \sigma_2 \sigma_1 \sigma_2}^+.$$
In general, for $n \geq 1$, the braid group $B_n$ and the corresponding braid monoid $B_n^+$ are given by
$$B_n \defeq \spkl{\sigma_1, \dotsc, \sigma_{n-1} \ \vline 
  \begin{array}{c}
  \sigma_i \sigma_{i+1} \sigma_i 
  = \sigma_{i+1} \sigma_i \sigma_{i+1} 
  \ {\rm for} \ 1 \leq i \leq n-2,
  \\
  \sigma_i \sigma_j 
  = \sigma_j \sigma_i \ {\rm for} \ \abs{i-j} \geq 2
  \end{array}},$$
$$B_n^+ \defeq \spkl{\sigma_1, \dotsc, \sigma_{n-1} \ \vline 
  \begin{array}{c}
  \sigma_i \sigma_{i+1} \sigma_i 
  = \sigma_{i+1} \sigma_i \sigma_{i+1} 
  \ {\rm for} \ 1 \leq i \leq n-2,
  \\
  \sigma_i \sigma_j 
  = \sigma_j \sigma_i \ {\rm for} \ \abs{i-j} \geq 2
  \end{array}}^+.$$
This corresponds to the case where $I = \gekl{1, \dotsc, n-1}$ and $m_{i,i+1} = m_{i+1,i} = 3$ for all $1 \leq i \leq n-2$ and $m_{i,j} = m_{j,i} = 2$ for all $\abs{i-j} \geq 2$.

These Artin groups form an interesting class of examples which is of interest for group theorists.

Another family of examples is given by Baumslag-Solitar groups and their presentations: For $k, l \geq 1$, define the group
$$
  B_{k,l} \defeq \spkl{a,b \ \vert \ a b^k = b^l a}
$$
and the monoid
$$
  B_{k,l}^+ \defeq \spkl{a,b \ \vert \ a b^k = b^l a}^+.
$$
Also, again for $k, l \geq 1$, define the group
$$
  B_{-k,l} \defeq \spkl{a,b \ \vert \ a = b^l a b^k}
$$
and the monoid
$$
  B_{-k,l}^+ \defeq \spkl{a,b \ \vert \ a = b^l a b^k}^+.
$$
These are the Baumslag-Solitar groups and the Baumslag-Solitar monoids. The reader may find more about the semigroup C*-algebras attached to Baumslag-Solitar monoids in \cites{Sp1,Sp2}.

Finally, let us mention the Thompson group and the Thompson monoid. The Thompson group is given by
$$
  F \defeq \spkl{x_0, x_1, \dotsc \ \vert \ x_n x_k = x_k x_{n+1} \ {\rm for} \ k < n}.
$$
This is just one possible presentation defining the Thompson group. There are other, for instance
$$
  F = \spkl{A,B \ \vert \ [AB^{-1},A^{-1}BA] = [AB^{-1},A^{-2}BA^2] = e}.
$$
The first presentation however has the advantage that it leads naturally to the definition of the Thompson monoid as
$$
  F^+ \defeq \spkl{x_0, x_1, \dotsc \ \vert \ x_n x_k = x_k x_{n+1} \ {\rm for} \ k < n}^+.
$$
The Thompson group is of great interest in group theory, in particular the question whether it is amenable or not is currently attracting a lot of attention. Therefore, it would be very interesting to study the Thompson monoid and its semigroup C*-algebra.

\subsection{Examples from rings in general, and number theory in particular}
\label{ex:rings}

Let us present another source for examples. This time, our semigroups come from rings. Let $R$ be a ring without zero-divisors ($x \neq 0$ is a zero-divisor if there exists $0 \neq y \in R$ with $xy = 0$). Then $R\reg = R \setminus \gekl{0}$ is a cancellative semigroup with respect to multiplication.

We can also construct the $ax+b$-semigroup $R \rtimes R\reg$. The underlying set is $R \times R\reg$, and multiplication is given by $(d,c)(b,a) = (d+cb,ca)$. It is a semidirect product for the canonical multiplicative action of $R\reg$ on $R$.

Another possibility would be to take an integral domain $R$, i.e., a commutative ring with unit not containing zero-divisors, and form the semigroup $M_n(R)\reg$ of $n \times n$-matrices over $R$ with non-vanishing determinant. We could also form the semidirect product $M_n(R) \rtimes M_n(R)\reg$ for the canonical multiplicative action of $M_n(R)\reg$ on $M_n(R)$. 

In particular, rings from number theory are interesting. Let $K$ be a number field, i.e., a finite extension of $\Qz$. Then the ring of algebraic integers $R$ in $K$ is given by
$$
  \menge{x \in K}{{\rm There} \ {\rm are} \ n \geq 1, \, a_{n-1}, \dotsc, a_0 \in \Zz \ {\rm with} \ x^n + a_{n-1} x^{n-1} + \dotso + a_0 = 0}.
$$
For instance, for the classical case $K = \Qz$, the ring of algebraic integers is given by the usual integers, $R = \Zz$. For the number field of Gaussian numbers, $K = \Qz[i]$, the ring of algebraic integers are given by the Gaussian integers, $R = \Zz[i]$. More generally, for the number field $K = \Qz[\zeta]$ generated by a root of unity $\zeta$, the ring of algebraic integers is given by $R = \Zz[\zeta]$. For the real quadratic number field $K = \Qz[\sqrt{2}]$, the ring of algebraic integers is given by $\Zz[\sqrt{2}]$, while for the real quadratic number field $K = \Qz[\sqrt{5}]$, the ring of algebraic integers is given by $R = \Zz[\frac{1 + \sqrt{5}}{2}]$.

Let us briefly mention an interesting invariant of number fields. Let $K$ be a number field with ring of algebraic integers $R$. We introduce an equivalence relation for non-zero ideals $\cI$ of $R$ by saying that $\mfa \sim \mfb$ if there exist $a,b \in R\reg$ with $b \mfa = a \mfb$. It turns out that with respect to multiplication of ideals, $\cI / {}_{\sim}$ becomes a finite abelian group. This is the class group $Cl_K$ of $K$. An outstanding open question in number theory is how to compute $Cl_K$, or even just the class number $h_K = \# Cl_K$, in a systematic and efficient way. It is not even known whether there are infinitely many (non-isomorphic) number fields with trivial class group (i.e., class number one). We refer the interested reader to \cite{Neu} for more details.

It is possible to consider more general semidirect products, in the more flexible setting of semigroups acting by endomorphisms on a group. Particular cases are discussed in \cite{CunAlgAct}. We also refer to \cites{CuVe,BLS,BS,Stamm} and the references therein for more examples and for results on the corresponding semigroup C*-algebras.

\subsection{Finitely generated abelian cancellative semigroups}

Finally, one more class of examples which illustrates quite well that the world of semigroups can be much more complicated than the world of groups: Consider finitely generated abelian cancellative semigroups, or monoids. For groups, we have a well understood structure theorem for finitely generated abelian groups. But for semigroups, this class of examples is interesting and challenging to understand. For instance, particular examples are given by numerical semigroups, i.e., semigroups of the type $P = \Nz \setminus F$, where $F$ is a finite subset of $\Nz$ such that $\Nz \setminus F$ is additively closed. For instance, we could take $F = \gekl{1}$ or $F = \gekl{1,3}$. We refer the interested reader to \cite{RG} and the references therein for more about numerical semigroups, and also to \cite{CunToricVar}.

\section{Preliminaries}

\subsection{Embedding semigroups into groups}
\label{PinG}

As we mentioned earlier, we need left cancellation for semigroups in our construction of semigroup C*-algebras. One way to ensure cancellation is to embed our semigroups into groups, i.e., to find an injective semigroup homomorphism from our semigroup into a group. In general, the question which semigroups embed into groups is quite complicated. Cancellation is necessary but not sufficient. Malcev gave the complete answer. He found an infinite list of conditions which are necessary and sufficient for group embeddability, and showed that any finite subset of his list is no longer sufficient. His list includes cancellation, which means both left cancellation and right cancellation. The latter means that for every $p,x,y \in P$, $xp = yp$ implies $x=y$. But Malcev's list also consists of conditions like the following:

\begin{center}
For every $a,b,c,d,u,v,x,y \in P$,
\setlength{\parindent}{0cm} \setlength{\parskip}{0cm}
 
$xa = yb$, $xc = yd$, $ua = vb$ implies $uc = vd$.
\setlength{\parindent}{0cm} \setlength{\parskip}{0.5cm}
\end{center}
We refer to \cite[\S~12]{CPII} for more details.

As explained in \cite[\S~12]{CPII}, if a semigroup $P$ embeds into a group, then there is a universal group embedding $P \into G_{\rm univ}$, meaning that for every homomorphism $P \to G$ of the semigroup $P$ to a group $G$, there is a unique homomorphism $G_{\rm univ} \to G$ which makes the diagram
\bgloz
  \xymatrix@C=20mm{
  P \ar@{^{(}->}[r] \ar[dr] & G_{\rm univ} \ar[d] \\
   & G
  }
\egloz
commutative.

Group embeddability is in general a complicated issue. Therefore, whenever it is convenient, we will simply assume that our semigroups can be embedded into groups. Verifying this assumption might be a challenge, for instance in the case of Artin monoids (compare \cite{Par}).

However, we would like to mention one sufficient condition for group embeddability. Let $P$ be a cancellative semigroup, i.e., $P$ is left and right cancellative. Furthermore, assume that $P$ is right reversible, i.e., for every $p,q \in P$, we have $Pp \cap Pq \neq \emptyset$. Here $Pp = \menge{xp}{x \in P}$. Then $P$ embeds into a group. Actually, the universal group in the universal group embedding of $P$ is given by an explicit construction as the group $G$ of left quotients. This means that $G$ consists of formal quotients of the form $q^{-1}p$, for all $q \in P$ and $p \in P$. We say that two such formal expressions $\ti{q}^{-1} \ti{p}$ and $q^{-1}p$ represent the same element in $G$ if there is $r \in P$ with $\ti{q} = rq$ and $\ti{p} = rp$. To multiply elements in $G$, we make use of right reversibility: Given $p,q,r,s \in P$, suppose we want to multiply $s^{-1}r$ with $q^{-1}p$. As $Pq \cap Pr \neq \emptyset$, there exist $x$ and $y$ in $P$ with $q = yr$. Thus $q^{-1}p = (xq)^{-1} (xp) = (yr)^{-1} (xp)$. Let us now make the following formal computation:
$$
  (s^{-1}r) (q^{-1}p) = (s^{-1}r)(yr)^{-1}(xp) = s^{-1} r r^{-1} y^{-1} (xp) = s^{-1} y^{-1} (xp) = (ys)^{-1} (xp).
$$
Motivated by this computation, we set
$$
(s^{-1}r)(q^{-1}p) \defeq (ys)^{-1}(xp).
$$
It is now straightforward to check that this indeed defines a group $G = P^{-1}P$, and that $P \to G, \, p \ma e^{-1} p$ is an embedding of our semigroup $P$ into our group $G$. Here $e$ is the identity of $P$. We can always arrange that $P$ has an identity by simply adjoining one if necessary. It is easy to see that this group embedding which we just constructed is actually the universal group embedding for $P$. We refer the reader to \cite[\S~1.10]{CPI} for more details.

Obviously, by symmetry, we also obtain that a cancellative semigroup $P$ embeds into a group, if $P$ is left reversible, i.e., if for every $p,q \in P$, we have $pP \cap qP \neq \emptyset$. In that case, $P$ embeds into its group $G$ of right quotients, $G = P P^{-1}$, and this is the universal group embedding for $P$.

For instance, both of these necessary conditions for group embeddability are satisfied for cancellative abelian semigroups. They are also satisfied for the Braid monoids $B_n^+$ introduced above.

The $ax+b$-semigroup $R \rtimes R\reg$ over an integral domain $R$ is right reversible, but if $R$ is not a field, then $R \rtimes R\reg$ is not left reversible.

The Thompson monoid is left reversible but not right reversible.

Finally, the non-abelian free monoid $\Nz^{*n}$ is neither left nor right reversible.

\subsection{Graph products}
\label{ss:GraphProducts}

We collect some basic facts about graph products which we will use later on in \S~\ref{sec:GraphProducts}. Basically, we follow \cite[\S~2]{CrLa1}. Let $\Gamma = (V,E)$ be a graph with vertices $V$ and edges $E$. Two vertices in $V$ are connected by at most one edge, and no vertex is connected to itself. Hence we view $E$ as a subset of $V \times V$. For every $v \in V$, assume that we are given a submonoid $P_v$ of a group $G_v$. We can then form the graph products $P \defeq \Gamma_{v \in V} P_v$ and $G \defeq \Gamma_{v \in V} G_v$. As we explained, the group $G$ is obtained from the free product $*_{v \in V} G_v$ by introducing the relations $xy = yx$ for all $x \in G_v$ and $y \in G_w$ with $(v,w) \in E$. Similarly, $P$ is defined as the monoid obtained from the free product $*_{v \in V} P_v$ by introducing the relations $xy = yx$ for all $x \in P_v$ and $y \in P_w$ with $(v,w) \in E$. As explained in \cite{CrLa1}, it turns out that for every $v$, the monoid $P_v$ sits in a canonical way as a submonoid inside the monoid $\Gamma_{v \in V} P_v$. Similarly, for each $v \in V$, the group $G_v$ sits in a canonical way as a subgroup inside the group $\Gamma_{v \in V} G_v$. Moreover, the monoid $P = \Gamma_{v \in V} P_v$ can be canonically embedded as a submonoid of the group $G = \Gamma_{v \in V} G_v$.

A typical element $g$ of $G = \Gamma_{v \in V} G_v$ is a product of the form $x_1 x_2 \dotsm x_l$, where $x_i \in G_{v_i}$ are all non-trivial. (To obtain the identity, we would have to allow the empty word, i.e., the case $l = 0$.) We distinguish between words like $x_1 x_2 \dotsm x_l$ and the element $g$ they represent in the graph product $G$ by saying that $x_1 x_2 \dotsm x_l$ is an expression for $g$. Let us now explain when two words are expressions for the same group element.

First of all, for a word like $x_1 x_2 \dotsm x_l$, we call the $x_i$s the syllables and $l$ the length of the word. We write $v(x_i)$ for the vertex $v_i \in V$ with the property that $x_i$ lies in $G_{v_i}$. Given a word
$$
  x_1 \dotsm x_i x_{i+1} \dotsm x_l
$$
with the property that $(v(x_i),v(x_{i+1})) \in E$, we can replace the subword $x_i x_{i+1}$ by $x_{i+1} x_i$. In this way, we transform the original word
$$
  x_1 \dotsm x_i x_{i+1} \dotsm x_l
$$
to the new word
$$
  x_1 \dotsm x_{i+1} x_i \dotsm x_l.
$$
This procedure is called a shuffle. Two words are called shuffle equivalent if one can be obtained from the other by performing finitely many shuffles.

Moreover, given a word
$$
  x_1 \dotsm x_i x_{i+1} \dotsm x_l
$$
with the property that $v(x_i) = v(x_{i+1})$, then we say that our word admits an amalgamation. In that case, we can replace the subword $x_i x_{i+1}$ by the product $x_i \cdot x_{i+1} \in G_{v_i}$, where $v_i = v(x_i) = v(x_{i+1})$. Furthermore, if $x_i \cdot x_{i+1} = e$ in $G_{v_i}$, then we delete this part of our word. In this way, we transform the original word
$$
  x_1 \dotsm x_i x_{i+1} \dotsm x_l
$$
to the new word
$$
  x_1 \dotsm (x_i \cdot x_{i+1}) \dotsm x_l
$$
if $x_i \cdot x_{i+1} \neq e$ in $G_{v_i}$ and
$$
  x_1 \dotsm x_{i-1} x_{i+2} \dotsm x_l
$$
if $x_i \cdot x_{i+1} = e$ in $G_{v_i}$. This procedure is called an amalgamation.

Finally, we say that a word is reduced if it is not shuffle equivalent to a word which admits an amalgamation.

We have the following
\blemma[Lemma~1 in \cite{CrLa1}]
\label{reduced-criterion}
A word
$$
  x_1 \dotsm x_l
$$
is reduced if and only if for all $1 \leq i < j \leq l$ with $v(x_i) = v(x_j)$, there exists $1 \leq k \leq l$ with $i < k < j$ such that $(v(x_i),v(x_k)) \notin E$.
\elemma

Suppose that we are given two words, and we can transform one word into the other by finitely many shuffles and amalgamations. Then it is clear that these two words are expressions for the same element in our group $G$. The converse is also true, this is the following result due to Green (see \cite{GreenPhD}):
\btheo[Theorem~2 in \cite{CrLa1}]
\label{graphprod-wordproblem}
Any two reduced words which are expressions for the same group element in $G$ are shuffle equivalent.
\etheo
In other words, two words which are expressions for the same group element in $G$ can be transformed into one another by finitely many shuffles and amalgamations. This is because, with the help of Lemma~\ref{reduced-criterion}, it is easy to see that every word can be transformed into a reduced one by finitely many shuffles and amalgamations.

Because of Theorem~\ref{graphprod-wordproblem}, we may introduce the notion of length:
\bdefin
The length of an element $g$ in our graph product $G$ is the length of a reduced word which is an expression for $g$.
\edefin

We also introduce the following
\bdefin
Suppose we are given a reduced word
$$
  x = x_1 \dotsm x_l.
$$
Then we call $x_i$ an initial syllable and $v(x_i)$ an initial vertex of our word, if for every $1 \leq h < i$, $(v(x_h),v(x_i)) \in E$. The set of all initial vertices of $x$ is denoted by $V^i(x)$ (in \cite{CrLa1}, the notation $\Delta(x)$ is used).

Similarly, we call $x_j$ a final syllable and $v(x_j)$ a final vertex of our word, if for every $j < k \leq l$, $(v(x_j),v(x_k)) \in E$. The set of all final vertices of $x$ is denoted by $V^f(x)$ (it is denoted by $\Delta^r(x)$ in \cite{CrLa1}).
\edefin

The following is an easy observation:
\blemma[Lemma~3 in \cite{CrLa1}]
Let
$$
  x = x_1 \dotsm x_l
$$
be a reduced word.

If $x_i$ is an initial syllable of $x$, then $x$ is shuffle equivalent to $x_i x_1 \dotsm x_{i-1} x_{i+1} \dotsm x_l$.

For all $v,w \in V^i(x)$, we have $(v,w) \in E$.

For every $v \in V^i(x)$, there is a unique initial syllable $x_i$ of $x$ with $v(x_i) = v$. Let us denote this syllable by $S_v^i(x)$.

If $x'$ is shuffle equivalent to $x$, then $V^i(x) = V^i(x')$ and for every $v \in V^i(x) = V^i(x')$,
$$S_v^i(x) = S_v^i(x').$$
\elemma

The last three statements are also true for final vertices and final syllables. So we denote for a reduced word $x$ with final syllable $v$ the unique final syllable $x_j$ of $x$ with $v(x_j) = v$ by $S_v^f(x)$.

\bdefin
Let $g$ be an element in our graph product $G$, and let $x$ be a reduced word which is an expression for $g$. Then we set
$$
  V^i(g) \defeq V^i(x),
$$
and for $v \in V$,
$$
  S_v^i(g) \defeq
  \bfa
    S_v^i(x) & {\rm if} \ v \in V^i(g)\\
    e & {\rm if} \ v \notin V^i(g).
  \efa
$$
Similarly, we define
$$
  V^f(g) \defeq V^f(x),
$$
and for $v \in V$,
$$
  S_v^f(g) \defeq
  \bfa
    S_v^f(x) & {\rm if} \ v \in V^f(g)\\
    e & {\rm if} \ v \notin V^f(g).
  \efa
$$
\edefin

We need the following
\blemma[Lemma~5 in \cite{CrLa1}]
\label{XYZ}
Given $g$ and $h$ in our graph product $G$, let
$$W \defeq V^f(g) \cap V^i(h),$$
and suppose that
$$z_w \defeq S_w^f(g) S_w^i(h) \neq e$$
for all $w \in W$. Define
$$z \defeq \prod_{w \in W} z_w,$$
in any order.

Then, if $x \cdot \prod_{w \in W} S_w^f(g)$ is a reduced expression for $g$ and $ \prod_{w \in W} S_w^i(h) \cdot y$ is a reduced expression for $h$, then $x \cdot z \cdot y$ is a reduced expression for $g \cdot h$.
\elemma

\subsection{Krull rings}
\label{ss:Krull}

Since we want to study $ax+b$-semigroups over integral domains and their semigroup C*-algebras later on, we collect a few basic facts in this context.

Let $R$ be an integral domain.
\bdefin
The constructible (ring-theoretic) ideals of $R$ are given by
$$
  \cI(R) \defeq \menge{c^{-1} \rukl{\bigcap_{i=1}^n a_i R}}{a_1, \dotsc, a_n, c \in R\reg}.
$$
\edefin
Here, for $c \in R\reg$ and an ideal $I$ of $R$, we set 
$$c^{-1}I \defeq \menge{r \in R}{cr \in I}.$$

Now let $Q$ be the quotient field of $R$.
\bgl
\label{const-ring-ideals}
  \cI(R \subseteq Q) \defeq \menge{(x_1 \cdot R) \cap \dotso \cap (x_n \cdot R)}{x_i \in Q\reg}.
\egl
Note that for $c \in R\reg$ and $X \subseteq R$, we set 
$$c^{-1} X = \menge{r \in R}{cr \in X}, \ {\rm but} \ {c^{-1} \cdot X} = \menge{c^{-1} x}{x \in X}.$$
Moreover, note that $\cI(R) = \menge{J \cap R}{J \in \cI(R \subseteq Q)}$.

By construction, the family $\cI(R)$ consists of integral divisorial ideals of $R$, and $\cI(R \subseteq Q)$ consists of divisorial ideals of $R$. By definition, a divisorial ideal of an integral domain $R$ is a fractional ideal $I$ that satisfies $I = (R:(R:I))$, where $(R:J) = \menge{q \in Q}{qJ \subseteq R}$. Equivalently, divisorial ideals are non-zero intersections of some non-empty family of principal fractional ideals (ideals of the form $qR$, $q \in Q$). Let $D(R)$ be the set of divisorial ideals of $R$. In our situation, we only consider finite intersections of principal fractional ideals (see \eqref{const-ring-ideals}). So in general, our family $\cI(R \subseteq Q)$ will only be a proper subset of $D(R)$.

However, for certain rings, the set $\cI(R \subseteq Q)$ coincides with $D(R)$. For instance, this happens for noetherian rings. It also happens for Krull rings. The latter have a number of additional favourable properties which are very helpful for our purposes. Let us start with the following
\bdefin
\label{Krull}
An integral domain $R$ is called a Krull ring if there exists a family of discrete valuations $(v_i)_{i \in I}$ of the quotient field $Q$ of $R$ such that
\begin{enumerate}
\item[(K1)] $R = \menge{x \in Q}{v_i(x) \geq 0 \fa i \in I}$,
\item[(K2)] for every $0 \neq x \in Q$, there are only finitely many valuations in $(v_i)_i$ such that $v_i(x) \neq 0$.
\end{enumerate}
\edefin

The following result gives us many examples of Krull rings.
\btheo{\cite[Chapitre~VII, \S~1.3, Corollaire]{Bour2}}
A noetherian integral domain is a Krull ring if and only if it is integrally closed.
\etheo

Let us collect some basic properties of Krull rings:

\cite[Chapitre~VII, \S~1.5, Corollaire~2]{Bour2} yields
\blemma
For a Krull ring $R$, $\cI(R \subseteq Q) = D(R)$ and $\cI(R)$ is the set of integral divisorial ideals.
\elemma
Moreover, the prime ideals of height $1$ play a distinguished role in a Krull ring.
\btheo{\cite[Chapitre~VII, \S~1.6, Th\'{e}or\`{e}me~3 and Chapitre~VII, \S~1.7, Th\'{e}or\`{e}me~4]{Bour2}}
Let $R$ be a Krull ring. Every prime ideal of height $1$ of $R$ is a divisorial ideal. Let
$$\cP(R) = \menge{\mfp \triangleleft R \text{ prime}}{\height(\mfp)=1}.$$
For every $\mfp \in \cP(R)$, the localization $R_{\mfp} = (R \setminus \mfp)^{-1} R$ is a principal valuation ring. Let $v_{\mfp}$ be the corresponding (discrete) valuations of the quotient field $Q$ of $R$. Then the family $(v_{\mfp})_{\mfp \in \cP(R)}$ satisfies the conditions (K1) and (K2) from Definition~\ref{Krull}.
\etheo
\bprop{\cite[Chapitre~VII, \S~1.5, Proposition~9]{Bour2}}
\label{approx}
Let $R$ be a Krull ring and $(v_{\mfp})_{\mfp \in \cP(R)}$ be the valuations from the previous theorem. Given finitely many integers $n_1$, ..., $n_r$ and finitely many prime ideals $\mfp_1$, ..., $\mfp_r$ in $\cP(R)$, there exists $x$ in the quotient field $Q$ of $R$ with
$$v_{\mfp_i}(x) = n_i \ {\rm for} \ {\rm all} \ 1 \leq i \leq r \ {\rm and} \ v_{\mfp}(x) \geq 0 \ {\rm for} \ {\rm all} \ \mfp \in \cP(R) \setminus \gekl{\mfp_1, \dotsc, \mfp_r}.
$$
\eprop
Moreover, given a fractional ideal $I$ of $R$, we let $I^\sim \defeq (R:(R:I))$ be the divisorial closure of $I$. $I^\sim$ is the smallest divisorial ideal of $R$ which contains $I$. We can now define the product of two divisorial ideals $I_1$ and $I_2$ to be the divisorial closure of the (usual ideal-theoretic) product of $I_1$ and $I_2$, i.e., $I_1 \bullet I_2 \defeq (I_1 \cdot I_2)^\sim$. $D(R)$ becomes a commutative monoid with this multiplication.
\btheo{\cite[Chapitre~VII, \S~1.2, Th\'{e}or\`{e}me~1; Chapitre~VII, \S~1.3, Th\'{e}or\`{e}me~2 and Chapitre~VII, \S~1.6, Th\'{e}or\`{e}me~3]{Bour2}}
For a Krull ring $R$, $(D(R),\bullet)$ is a group. It is the free abelian group with free generators given by $\cP(R)$, the set of prime ideals of $R$ which have height $1$.
\etheo
This means that every $I \in \cI(R \subseteq Q)$ ($Q$ is the quotient field of the Krull ring $R$) is of the form $I = \mfp_1^{(n_1)} \bullet \dotsm \bullet \mfp_r^{(n_r)}$, with $n_i \in \Zz$. Here for $\mfp \in \cP(R)$ and $n \in \Nz$, we write
$$\mfp^{(n)} \ {\rm for} \ \underbrace{\mfp \bullet \dotsm \bullet \mfp}_{n \text{ times}}, \ {\rm and} \ \mfp^{(-n)} \ {\rm for} \ \underbrace{\mfp^{-1} \bullet \dotsm \bullet \mfp^{-1}}_{n \text{ times}},$$
where $\mfp^{-1} = (R:\mfp)$. We set for $\mfp \in \cP(R)$:
$$
  v_{\mfp}(I) \defeq
  \bfa
    n_i \falls \mfp = \mfp_i, \\
    0 \falls \mfp \notin \gekl{\mfp_1, \dotsc, \mfp_r}.
  \efa
$$
With this notation, we have $I = \prod_{\mfp \in \cP(R)} \mfp^{(v_{\mfp}(I))}$, where the product is taken in $D(R)$. In addition, we have for $I \in \cI(R \subseteq Q)$ that $I \in \cI(R)$ if and only if $v_{\mfp}(I) \geq 0$ for all $\mfp \in \cP(R)$. And combining the last statement in \cite[Chapitre~VII, \S~1.3, Th\'{e}or\`{e}me~2]{Bour2} with \cite[Chapitre~VII, \S~1.4, Proposition~5]{Bour2}, we obtain for every $I \in \cI(R \subseteq Q)$:
\bgl
\label{I=}
  I = \menge{x \in Q}{v_{\mfp}(x) \geq v_{\mfp}(I) \text{ for all } \mfp \in \cP(R)}.
\egl
Finally, the principal fractional ideals $F(R)$ form a subgroup of $(D(R),\bullet)$ which is isomorphic to $Q\reg$. Suppose that $R$ is a Krull ring. Then the quotient group $C(R) \defeq D(R) / F(R)$ is called the divisor class group of $R$.

These were basic properties of Krull rings. We refer the interested reader to \cite[Chapitre~VII]{Bour2} or \cite{Fos} for more information.

\section{C*-algebras attached to inverse semigroups, partial dynamical systems, and groupoids}

We refer the interested reader to \cites{Ren1,Ex2,Ex4,Pat99} for more references for this section.

\subsection{Inverse semigroups}
\label{ss:ISGP}

Inverse semigroups play an important role in the study of semigroup C*-algebras.

\bdefin
An inverse semigroup is a semigroup $S$ with the property that for every $x \in S$, there is a unique $y \in S$ with $x = xyx$ and $y = yxy$. 
\edefin
We write $y = x^{-1}$ and call $y$ the inverse of $x$. 

\bdefin
An inverse semigroup $S$ is called an inverse semigroup with zero if there is a distinguished element $0 \in S$ satisfying $0 \cdot s = 0 = s \cdot 0$ for all $s \in S$. 
\edefin

Usually, if we write \an{inverse semigroup}, we mean an inverse semigroup with or without zero. Sometimes we write \an{inverse semigroups without zero} for ordinary inverse semigroups which do not have a distinguished zero element.

Every inverse semigroup can be realized as partial bijections on a fixed set. Multiplication is given by composition. However, a partial bijection is only defined on its domain. Therefore, if we want to compose the partial bijection $s: \: \dom(s) \to \img(s)$ with another partial bijection $t: \: \dom(t) \to \img(t)$, we have to restrict $t$ to
$$\dom(t) \cap t^{-1}(\dom(s))$$
to make sure that the image of the restriction of $t$ lies in the domain of $s$. Only then we can form $s \circ t$. The inverse of a partial bijection is the usual inverse, in the category of sets.

Inverse semigroups can also be realized as partial isometries on a Hilbert space. To make sure that the product of two partial isometries is again a partial isometry, we have to require that the source and range projections of our partial isometries commute. Then multiplication in the inverse semigroup is just the usual multiplication of operators on a fixed Hilbert space, i.e., composition of operators. The inverse in our inverse semigroup is given by the adjoint operation for operators in general or partial isometries in our particular situation.

Let us explain how to attach an inverse semigroup to a left cancellative semigroup. Assume that $P$ is a left cancellative semigroup. Its left inverse hull $I_l(P)$ is the inverse semigroup generated by the partial bijections
$$P \to pP, \, x \ma px,$$
whose domain is $P$ and whose image is $pP = \menge{px}{x \in P}$. Its inverse is given by
$$pP \to P, \, px \ma x.$$
So $I_l(P)$ is the smallest semigroup of partial bijections on $P$ which is closed under inverses and contains
$$\menge{P \to P, x \ma px}{p \in P}.$$
Given $p \in P$, we denote the partial bijection
$$P \to pP, \, x \ma px$$
by $p$. In this way, we obtain an embedding of $P$ into $I_l(P)$ by sending $p \in P$ to the partial bijection $p \in I_l(P)$. This allows us to view $P$ as a subsemigroup of $I_l(P)$. We say that $I_l(P)$ is an inverse semigroup with zero if the partial bijection which is nowhere defined, $\emptyset \to \emptyset$, is in $I_l(P)$. In that case, $\emptyset \to \emptyset$ is the distinguished zero element $0$.

Alternatively, we can also describe $I_l(P)$ as the smallest inverse semigroup of partial isometries on $\ell^2 P$ generated by the isometries $\menge{V_p}{p \in P}$. This means that $I_l(P)$ can be identified with the smallest semigroup of partial isometries on $\ell^2 P$ containing the isometries $\menge{V_p}{p \in P}$ and their adjoints $\menge{V_p^*}{p \in P}$ and which is closed under multiplication. In this picture, $I_l(P)$ is an inverse semigroup with zero if and only if the zero operator is in $I_l(P)$.

An important subsemigroup of an inverse semigroup $S$ is its semilattice of idempotents.
\bdefin
The semilattice $E$ of idempotents in an inverse semigroup $S$ is given by
$$E \defeq \menge{x^{-1} x}{x \in S} = \menge{x x^{-1}}{x \in S} = \menge{e \in S}{e = e^2}.$$
Define an order on $E$ by setting, for $e,f \in E$, $e \leq f$ if $e = ef$.

If $S$ is an inverse semigroup with zero, $E$ becomes a semilattice with zero, and the distinguished zero element of $S$ becomes the distinguished zero element of $E$.
\edefin

In the case of partial bijections, the semilattice of idempotents is given by all domains and images. Multiplication in this semilattice is intersection of sets, and $\leq$ is $\subseteq$ for sets, i.e., containment.

\bdefin
For the left inverse hull $I_l(P)$ attached to a left cancellative semigroup $P$, the semilattice of idempotents is denoted by $\cJ_P$.
\edefin

It is easy to see that $\cJ_P$ is given by
$$\cJ_P = \menge{p_n \dotsm q_1^{-1} p_1 (P)}{q_i, p_i \in P} \cup \menge{q_n^{-1} p_n \dotsm q_1^{-1} p_1 (P)}{q_i, p_i \in P}.$$
Here, for $X \subseteq P$ and $p,q \in P$, we write
$$p(X) = \menge{px}{x \in X}$$
and
$$q^{-1} (X) = \menge{y \in P}{qy \in X}.$$

Subsets of the form $p_n \dotsm q_1^{-1} p_1 (P)$ or $q_n^{-1} p_n \dotsm q_1^{-1} p_1 (P)$ are right ideals of $P$. Here, we call $X \subseteq P$ a right ideal if for every $x \in X$ and $r \in P$, we always have $xr \in X$.

\bdefin
The elements in $\cJ_P$ are called constructible right ideals of $P$.
\edefin
We will work out the set of constructible right ideals explicitly for classes of examples in \S~\ref{independence}.

There is a duality between semilattices, i.e., abelian semigroups of idempotents, and totally disconnected locally compact Hausdorff spaces. Given a semilattice $E$, we construct its space of characters $\widehat{E}$ as follows:
$$
  \widehat{E} = \gekl{\chi: \: E \to \gekl{0,1} \ \text{non-zero} \ {\rm semigroup} \ {\rm homomorphism}}.
$$
In other words, elements in $\widehat{E}$ are multiplicative maps from $E$ to $\gekl{0,1}$, where the latter set is equipped with the usual multiplication when we view it as a subspace of $\Rz$ (or $\Cz$). In addition, we require that these multiplicative maps must take the value $1$ for some element $e \in E$. If our semilattice $E$ is a semilattice with zero, and $0$ is its distinguished zero element, then we require that $\chi(0) = 0$ for all $\chi \in \widehat{E}$.

The topology on $\widehat{E}$ is given by pointwise convergence. Every $\chi \in \widehat{E}$ is uniquely determined by
$$\chi^{-1}(1) = \menge{e \in E}{\chi(e) = 1}.$$
$\chi^{-1}(1)$ is an $E$-valued filter (which we simply call filter from now on), i.e., a subset of $E$ satisfying:
\begin{itemize}
\item $\chi^{-1}(1) \neq \emptyset$.
\item For all $e,f \in E$ with $e \leq f$, $e \in \chi^{-1}(1)$ implies $f \in \chi^{-1}(1)$.
\item For all $e,f \in E$ with $e,f \in \chi^{-1}(1)$, $ef$ lies in $\chi^{-1}(1)$.
\end{itemize}
Conversely, every filter, i.e., every subset $\cF \in E$ satisfying these three conditions determines a unique $\chi \in \widehat{E}$ with $\chi^{-1}(1) = \cF$. Therefore, we have a one-to-one correspondence between characters $\chi \in \widehat{E}$ and filters.

If $E$ is a semilattice with zero, and $0$ is the distinguished zero element, then we have $\chi(0) = 0$ for all $\chi \in \widehat{E}$. In terms of filters, this amounts to saying that $0$ is never an element of a filter. 
 
As an illustrative example, the reader is encouraged to work out the set of constructible right ideals $\cJ_P$ and the space of characters $\widehat{\cJ_P}$ for the non-abelian free semigroup on two generators $P = \Nz * \Nz$, or in other words, the semilattice of idempotents $E$ and the space of its characters $\widehat{E}$ for the inverse semigroup $S = I_l(\Nz * \Nz)$.

Now assume that we are given a subsemigroup $P$ of a group $G$.  We define
$$I_l(P)\reg \defeq I_l(P) \setminus \gekl{0}$$
if $I_l(P)$ is an inverse semigroup with zero, and $0$ is its distinguished zero element, and
$$I_l(P)\reg \defeq I_l(P)$$
otherwise.

Now it is easy to see that for every partial bijection $s$ in $I_l(P)\reg$, there exists a unique $\sigma(s) \in G$ such that $s$ is of the form 
$$s(x) = \sigma(s) \cdot x \ {\rm for} \ x \in \dom(s).$$
Here we view $P$ as a subset of the group $G$ and make use of multiplication in $G$.

In the alternative picture of $I_l(P)$ as the inverse semigroup of partial isometries on $\ell^2 P$ generated by the isometries $\menge{V_p}{p \in P}$, $I_l(P)\reg$ is given by all non-zero partial isometries in $I_l(P)$. Every element in $I_l(P)\reg$ is of the form
$$
V_{q_1} \dotsm V_{p_n}^*, \,
V_{p_1}^* V_{q_1} \dotsm V_{p_n}^*, \,
V_{q_1} \dotsm V_{p_n}^* V_{q_n}, \, {\rm or} \ 
V_{p_1}^* V_{q_1} \dotsm V_{p_n}^* V_{q_n}.
$$
The map $\sigma$ which we introduced above is then given by
\bglnoz
\sigma(V_{q_1} \dotsm V_{p_n}^*) &=& q_1 \dotsm p_n^{-1} \in G,\\
\sigma(V_{p_1}^* V_{q_1} \dotsm V_{p_n}^*) &=& p_1^{-1} q_1 \dotsm p_n^{-1} \in G,\\
\sigma(V_{q_1} \dotsm V_{p_n}^* V_{q_n}) &=& q_1 \dotsm p_n^{-1} q_n \in G,\\
{\rm or} \ \sigma(V_{p_1}^* V_{q_1} \dotsm V_{p_n}^* V_{q_n}) &=& p_1^{-1} q_1 \dotsm p_n^{-1} q_n \in G.
\eglnoz

To see that $\sigma$ is well-defined, note that, similarly as above, every partial isometry $V \in I_l(P)\reg$ has the property that there exists a unique $g \in G$ such that for every $x \in P$, either $V \delta_x = 0$ or $V \delta_x = \delta_{g \cdot x}$. And $\sigma$ is defined in such a way that $\sigma(V) = g$.

It is easy to see that the map $\sigma: \: I_l(P)\reg \to G$ satisfies
$$\sigma(st) = \sigma(s) \sigma(t)$$
for all $s,t \in I_l(P)\reg$, as long as the product $st$ lies in $I_l(P)\reg$, i.e., is non-zero. Moreover, setting
$$\cJ_P\reg \defeq \cJ_P \setminus \gekl{0}$$
if $\cJ_P$ is a semilattice with zero, and $0$ is the distinguished zero element, and
$$\cJ_P\reg \defeq \cJ_P$$
otherwise, it is also easy to see that
$$\sigma^{-1}(e) = \cJ_P\reg.$$
Here $e$ is the identity in our group $G$.

We formalize this in the next definition: Let $S$ be an inverse semigroup and $E$ the semilattice of idempotents of $S$. We set
$S\reg \defeq S \setminus \gekl{0}$
if $S$ is an inverse semigroup with zero, and $0$ is the distinguished zero element, and
$S\reg \defeq S$
otherwise. Similarly, let
$E\reg \defeq E \setminus \gekl{0}$
if $E$ is a semilattice with zero, and $0$ is the distinguished zero element, and
$E\reg \defeq E$
otherwise. Moreover, let $G$ be a group. 
\bdefin
A map $\sigma: \: S\reg \to G$ is called a partial homomorphism if $\sigma(st) = \sigma(s) \sigma(t)$ for all $s,t \in S\reg$ with $st \in S\reg$.

A map $\sigma: \: S\reg \to G$ is called idempotent pure if $\sigma^{-1}(e) = E\reg$.
\edefin
The existence of an idempotent pure partial homomorphism will allow us to describe C*-algebras attached to inverse semigroups as crossed products of partial dynamical systems later on.

The following is a useful observation which we need later on.
\blemma
\label{dom-sigma==}
Assume that $S$ is an inverse semigroup and $\sigma: \: S\reg \to G$ is an idempotent pure partial homomorphism to a group $G$. Whenever two elements $s$ and $t$ in $S\reg$ satisfy $s^{-1}s = t^{-1}t$ and $\sigma(s) = \sigma(t)$, then we must have $s = t$.
\elemma
\bproof
It is clear that $st^{-1}$ lies in $S\reg$. Since $\sigma(st^{-1}) = e$, we must have $st^{-1} \in E$. Hence
$$
  st^{-1} = ts^{-1}st^{-1} = tt^{-1},
$$
and therefore
$$
  s = s s^{-1} s = s t^{-1} t = tt^{-1}t = t.
$$
\eproof

Let us now explain the construction of reduced and full C*-algebras for inverse semigroups.

Let $S$ be an inverse semigroup, and define $S\reg$ as above. For $s \in S$, define
$$\lambda_s: \: \ell^2 S\reg \to \ell^2 S\reg$$
by setting
$$  
  \lambda_s(\delta_x) \defeq \delta_{sx} \ {\rm if} \ s^{-1} s \geq x x^{-1}, \ {\rm and} \ \lambda_s(\delta_x) \defeq 0 \ {\rm otherwise} .
$$
Note that we require $s^{-1} s \geq x x^{-1}$ because on
$$
  \menge{x \in S}{s^{-1} s \geq x x^{-1}},
$$
the map $x \ma sx$ given by left multiplication with $s$ is injective. This is because we can reconstruct $x$ from $sx$ due to the computation
$$
  x = x x^{-1} x = s^{-1} s x x^{-1} x = s^{-1} (sx).
$$
Therefore, for each $s$, we obtain a partial isometry $\lambda_s$ by our construction. The assignment $s \ma \lambda_s$ is a *-representation of $S$ by partial isometries on $\ell^2 S\reg$. It is called the left regular representation of $S$. The star in *-representation indicates that we have $\lambda_{s^{-1}} = \lambda_s^*$.
\bdefin
We define
$$C^*_{\lambda}(S) \defeq C^*(\menge{\lambda_s}{s \in S}) \subseteq \cL(\ell^2 S\reg).$$
$C^*_{\lambda}(S)$ is called the reduced inverse semigroup C*-algebra of $S$.
\edefin

The full C*-algebra of an inverse semigroup $S$ is given by a universal property.
\bdefin
\label{DEF:C*S}
We define
$$
C^*(S) \defeq C^* \rukl{ \gekl{v_s}_{s \in S} \ \vline \ 
  v_s v_t = v_{st}, \ v_s^* = v_{s^{-1}}, \ v_0 = 0 \ {\rm if} \ 0 \in S
  }.
$$
$C^*(S)$ is the full inverse semigroup C*-algebra of $S$.
\edefin
Here, $0 \in S$ is short for
\setlength{\parindent}{0cm} \setlength{\parskip}{0cm}

\begin{center}
\an{$S$ is an inverse semigroup with zero, and $0$ is the distinguished zero element}.
\end{center}
\setlength{\parindent}{0cm} \setlength{\parskip}{0.5cm}

This means that $C^*(S)$ is uniquely determined by the property that given any C*-algebra $B$ with elements $\menge{w_s}{s \in S}$ satisfying the above relations, i.e.,
$$
  w_s w_t = w_{st}, \ w_s^* = w_{s^{-1}}, \ w_0 = 0 \ {\rm if} \ 0 \in S,
$$
then there exists a unique *-homomorphism from $C^*(S)$ to $B$ sending $v_s$ to $w_s$. 

In other words, $C^*(S)$ is the C*-algebra universal for *-representation of $S$ by partial isometries (in a C*-algebra, or on a Hilbert space). Note that we require that if $0 \in S$, then the zero element of $S$ should be represented by the partial isometry $0$. That is why $v_0 = 0$ in case $0 \in S$. This is different from the definition in \cite[\S~2.1]{Pat99}, where the partial isometry representing $0$ in the full C*-algebra of $S$ is a non-zero, minimal and central projection. We will come back to this difference in the definitions later on.

By construction, there is a canonical *-homomorphism $\lambda: \: C^*(S) \to C^*_{\lambda}(S), \, v_s \ma \lambda_s$. It is called the left regular representation (of $C^*(S)$).

We refer the reader to \cite{Pat99} for more about inverse semigroups and their C*-algebras.

\subsection{Partial dynamical systems}
\label{ss:PDS}

Whenever we have a semigroup embedded into a group, or an inverse semigroup with an idempotent pure partial homomorphism to a group, we can construct a partial dynamical system. Let us first present the general framework.

In the following, our convention will be that all our groups are discrete and countable, and all our topological spaces are locally compact, Hausdorff and second countable. 

\bdefin
Let $G$ be a group with identity $e$, and let $X$ be a topological space.

A partial action $\alpha$ of $G$ on $X$ consists of
\begin{itemize}
\item a collection $\gekl{U_g}_{g \in G}$ of open subspaces $U_g \subseteq X$,
\item a collection $\gekl{\alpha_g}_{g \in G}$ of homeomorphisms $\alpha_g: \: U_{g^{-1}} \to U_g, \, x \ma g.x$ such that
\begin{itemize}
\item $U_e = X$, $\alpha_e = \id_X$;
\item for all $g_1, g_2 \in G$, we have
\bglnoz
&& g_2.(U_{(g_1 g_2)^{-1}} \cap U_{g_2^{-1}}) = U_{g_2} \cap U_{g_1^{-1}},\\
&{\rm and}& (g_1 g_2).x = g_1.(g_2.x) \ {\rm for} \ {\rm all} \ x \in U_{(g_1 g_2)^{-1}} \cap U_{g_2^{-1}}.
\eglnoz
\end{itemize}
\end{itemize}
We call such a triple $(X,G,\alpha)$ a partial dynamical system, and denote it by $\alpha: \: G \curvearrowright X$ or simply $G \curvearrowright X$. 
\edefin

Let $\alpha: \: G \curvearrowright X$ be a partial dynamical system. The dual action $\alpha^*$ of $\alpha$ is the partial action (in the sense of \cite{McCl}) of $G$ on $C_0(X)$ given by $$\alpha^*_g: \: C_0(U_{g^{-1}}) \to C_0(U_g), \, f \ma f(g^{-1}.\sqcup).$$

We set out to describe a canonical partial action attached to a semigroup $P$ embedded into a group $G$. Let $C^*_{\lambda}(P)$ be the reduced semigroup C*-algebra of $P$. It contains a canonical commutative subalgebra $D_{\lambda}(P)$, which is given by 
$$D_{\lambda}(P) \defeq C^*(\menge{1_X}{X \in \cJ_P}) \subseteq C^*_{\lambda}(P).$$
It is clear that $D_{\lambda}(P)$ coincides with $\clspan(\sigma^{-1}(e))$. Recall that the map $\sigma: \: I_l(P)\reg \to G$ is given as follows: Every partial isometry $V \in I_l(P)\reg$ has the property that there exists a unique $g \in G$ such that for every $x \in P$, either $V \delta_x = 0$ or $V \delta_x = \delta_{g \cdot x}$. And $\sigma$ is defined in such a way that $\sigma(V) = g$.

Let us now describe the canonical partial action $G \curvearrowright D_{\lambda}(P)$. We will think of it as a dual action $\alpha^*$. For $g \in G$, let
$$D_{g^{-1}} \defeq \clspan(\menge{V^*V}{V \in I_l(P)\reg, \, \sigma(V) = g}).$$
By construction, we have that $D_e = D_{\lambda}(P)$. Moreover, it is easy to see that $D_{g^{-1}}$ is an ideal of $D_{\lambda}(P)$. Here is the argument: Suppose we are given $V \in I_l(P)\reg$ with $\sigma(V) = g$, and $W \in I_l(P)\reg$ with $\sigma(W) = e$. Then $W$ must be a projection since for every $x \in P$, either $W \delta_x = 0$ or $W \delta_x = \delta_{e \cdot x} = \delta_x$. Moreover, $W$ and $V^*V$ commute as both of these are elements in the commutative C*-algebra $\ell^{\infty}(P)$. Hence $WV^*V$ is non-zero if and only if $V^*VW$ is non-zero, and if that is the case, we obtain
$$
  WV^*V = V^*VW = WV^*VW = (VW)^*(VW).
$$
As $\sigma(VW) = g$, this implies that both $WV^*V$ and $V^*VW$ lie in $D_{g^{-1}}$. Therefore, as we claim, $D_{g^{-1}}$ is an ideal of $D_{\lambda}(P)$.

We then define $\alpha^*_g$ as $\alpha^*_g : \: D_{g^{-1}} \to D_g, \, V^*V \to VV^*$ for $V \in I_V\reg$ with $\sigma(V) = g$. This is well-defined: If we view $\ell^2 P$ as a subspace $\ell^2 G$ and let $\lambda$ be the left regular representation of $G$, then every $V \in I_V\reg$ with $\sigma(V) = g$ satisfies $V = \lambda_g V^*V$. Therefore, $VV^* = \lambda_g V^*V \lambda_g^*$. This shows that $\alpha^*_g$ is just conjugation with the unitary $\lambda_g$. This also explains why $\alpha^*_g$ is an isomorphism.

Of course, we can also describe the dual action $\alpha$. Set
$$\Omega_P \defeq \Spec(D_{\lambda}(P))$$
and for every $g \in G$, let
$$U_{g^{-1}} \defeq \widehat{D_{g^{-1}}}.$$
It is easy to see that
$$U_{g^{-1}} = \menge{\chi \in \Omega_P}{\chi(V^*V) = 1 \ {\rm for} \ {\rm some} \ V \in I_V\reg \ {\rm with} \ \sigma(V) = g}.$$
We then define $\alpha_g$ by setting $\alpha_g(\chi) \defeq \chi \circ \alpha^*_{g^{-1}}$. These $\alpha_g$, $g \in G$, give rise to the canonical partial dynamical system $G \curvearrowright \Omega_P$ attached to a semigroup $P$ embedded into a group $G$.

Our next goal is to describe a canonical partial dynamical system attached to inverse semigroups equipped with a idempotent pure partial homomorphism to a group. Let $S$ be an inverse semigroup and $E$ the semilattice of idempotents of $S$. Let $G$ be a group. Assume that $\sigma$ is a partial homomorphism $S\reg \to G$ which is idempotent pure.

In this situation, we describe a partial dynamical system $G \curvearrowright \widehat{E}$, and we will show later (see Corollary~\ref{C*S=C*ExG}) that the reduced C*-algebra $C^*_{\lambda}(S)$ of $S$ is canonically isomorphic to $C_0(\widehat{E}) \rtimes_r G$. 

Consider the sub-C*-algebra
$$C^*(E) \defeq C^*(\menge{\lambda_e}{e \in E}) \subseteq C^*_{\lambda}(S).$$
As we will see, we have a canonical isomorphism $\Spec(C^*(E)) \cong \widehat{E}$, so that $C_0(\widehat{E}) \cong C^*(E)$.

Now let us describe the partial action $G \curvearrowright C^*(E)$. For $g \in G$, define a sub-C*-algebra of $C^*(E)$ by 
$$C^*(E)_{g^{-1}} \defeq \clspan(\menge{\lambda_{s^{-1}s}}{s \in S\reg, \sigma(s) = g}).$$
As $\sigma$ is idempotent pure, we have $C^*(E)_e = C^*(E)$. For every $g \in G$, we have a C*-isomorphism
$$\alpha^*_g: \: C^*(E)_{g^{-1}} \to C^*(E)_g, \, \lambda_{s^{-1}s} \ma \lambda_{ss^{-1}}.$$

The corresponding dual action is given as follows: We identify $\Spec(C^*(E))$ with $\widehat{E}$. Then, for every $g \in G$, we set
$$U_g = \Spec(C^*(E)_g) \subseteq \widehat{E}.$$
It is easy to see that
$$U_{g^{-1}} = \menge{\chi \in \widehat{E}}{\chi(s^{-1}s) = 1 \ {\rm for} \ {\rm some} \ s \in S\reg \ {\rm with} \ \sigma(s) = g}.$$
For every $g \in G$, the homeomorphism $\alpha_g: \: U_{g^{-1}} \to U_g$ defining the partial dynamical system $G \curvearrowright \widehat{E}$ is given by $\alpha_g(\chi) = \chi \circ \alpha^*_{g^{-1}}$. More concretely, given $\chi \in U_{g^{-1}}$ and $s \in S\reg$ with $\sigma(s) = g$ and $\chi(s^{-1}s) = 1$, we have $\alpha_g(\chi)(e) = \chi(s^{-1}es)$. These $\alpha_g$, $g \in G$, give rise to the canonical partial dynamical system $G \curvearrowright \widehat{E}$ attached to an inverse semigroup $S$ equipped with an idempotent pure partial homomorphism to a group $G$.

At this point, a natural question arises. Assume we are given a semigroup $P$ embedded into a group $G$. We have seen above that this leads to an idempotent pure partial homomorphism on the left inverse hull $I_l(P)$ to our group $G$. How is the partial dynamical system $G \curvearrowright \Omega_P$ related to the partial dynamical system $G \curvearrowright \widehat{\cJ_P}$? We will see the answer in \S~\ref{cropro-GPD-description-C*redP}.

Let us now recall the construction, originally defined in \cite{McCl}, of the reduced and full crossed products $C_0(X) \rtimes_{\alpha^*,r} G$ and $C_0(X) \rtimes_{\alpha^*} G$ attached to our partial dynamical system $\alpha: \: G \curvearrowright X$. We usually omit $\alpha^*$ in our notation for the crossed products for the sake of brevity.

First of all,
$$C_0(X) \rtimes^{\ell^1} G \defeq \menge{\sum_g f_g \delta_g \in \ell^1(G,C_0(X))}{f_g \in C_0(U_g)}$$
becomes a *-algebra under component-wise addition, multiplication given by
$$\rukl{\sum_g f_g \delta_g} \cdot \rukl{\sum_h \ti{f}_h \delta_h} \defeq \sum_{g,h} \alpha^*_g(\alpha^*_{g^{-1}}(f_g) \ti{f}_h) \delta_{gh}$$
and involution 
$$\rukl{\sum_g f_g \delta_g}^* \defeq \sum_g \alpha^*_g(f_{g^{-1}}^*) \delta_g.$$

As in \cite{McCl}, we construct a representation of $C_0(X) \rtimes^{\ell^1} G$. Viewing $X$ as a discrete set, we define $\ell^2 X$ and the representation
$$M: \: C_0(X) \to \cL(\ell^2 X), \, f \ma M(f),$$
where $M(f)$ is the multiplication operator $M(f)(\xi) \defeq f \cdot \xi$ for $\xi \in \ell^2 X$. $M$ is obviously a faithful representation of $C_0(X)$. Every $g \in G$ leads to a twist of $M$, namely
$$M_g: \: C_0(X) \to \cL(\ell^2 X) \ {\rm given} \ {\rm by} \ M_g(f) \xi \defeq f \vert_{U_g}(g. \sqcup) \cdot \xi \vert_{U_{g^{-1}}}.$$
Here we view $f \vert_{U_g}(g. \sqcup)$ as an element in $C_b(U_{g^{-1}})$, and $C_b(U_{g^{-1}})$ acts on $\ell^2 U_{g^{-1}}$ just by multiplication operators. Given $\xi \in \ell^2 X$, we set
$$\xi \vert_{U_{g^{-1}}}(x) \defeq \xi(x) \ {\rm if} \ x \in U_{g^{-1}} \ {\rm and} \ \xi \vert_{U_{g^{-1}}}(x) \defeq 0 \ {\rm if} \ x \notin U_{g^{-1}}.$$
In other words, $\xi \vert_{U_{g^{-1}}}$ is the component of $\xi$ in $\ell^2 U_{g^{-1}}$ with respect to the decomposition
$$\ell^2 X = \ell^2 U_{g^{-1}} \oplus \ell^2 U_{g^{-1}}^c.$$
So we have
$$M_g(f) \xi (x) = f(g.x) \xi(x) \ {\rm if} \ x \in U_{g^{-1}} \ {\rm and} \ M_g(f) \xi (x) = 0 \ {\rm if} \ x \notin U_{g^{-1}}.$$

Consider now the Hilbert space
$$H \defeq \ell^2(G, \ell^2 X) \cong \ell^2 G \otimes \ell^2 X,$$
and define the representation
$$\mu: \: C_0(X) \to \cL(H) \ {\rm given} \ {\rm by} \ \mu(f)(\delta_g \otimes \xi) \defeq \delta_g \otimes M_g(f) \xi.$$
For $g \in G$, let $E_g$ be the orthogonal projection onto $\overline{\mu(C_0(U_{g^{-1}})) H}$. Moreover, let $\lambda$ denote the left regular representation of $G$ on $\ell^2 G$, and set $V_g \defeq (\lambda_g \otimes I) \cdot E_g$. Here $I$ is the identity operator on $H$.

We can now define the representation
$$\mu \times \lambda: \: C_0(X) \rtimes^{\ell^1} G \to \cL(H), \, \sum_g f_g \delta_g \ma \sum_g \mu(f_g) V_g.$$
Following the original definition in \cite{McCl}, we set 
\bdefin
$$
C_0(X) \rtimes_r G \defeq \overline{C_0(X) \rtimes^{\ell^1} G}^{\norm{\cdot}_{\mu \times \lambda}}.
$$
\edefin

To define the full crossed product $C_0(X) \rtimes G$ attached to our partial dynamical system $G \curvearrowright X$, recall that we have already introduced the *-algebra $C_0(X) \rtimes^{\ell^1} G$.

\bdefin
\label{DEF:full-cropro-partial}
Let $C_0(X) \rtimes G$ be the universal enveloping C*-algebra of the *-algebra $C_0(X) \rtimes^{\ell^1} G$.
\edefin

This means that $C_0(X) \rtimes G$ is universal for *-representations of $C_0(X) \rtimes^{\ell^1} G$ as bounded operators on Hilbert spaces or to C*-algebras. To construct this universal C*-algebra, we follow the usual procedure of completing $C_0(X) \rtimes^{\ell^1} G$ with respect to the maximal C*-norm on $C_0(X) \rtimes^{\ell^1} G$. Usually, we only obtain a C*-seminorm and have to divide out vectors with trivial seminorm, but because the *-representation $\mu \times \lambda$ constructed above is faithful, we get a C*-norm. So there is an embedding $C_0(X) \rtimes^{\ell^1} G \into C_0(X) \rtimes G$, and the universal property of $C_0(X) \rtimes G$ means that whenever we have a *-homomorphism $C_0(X) \rtimes^{\ell^1} G \to B$ to some C*-algebra $B$, there is a unique *-homomorphism $C_0(X) \rtimes G \to B$ which makes the diagram
\bgloz
  \xymatrix@C=20mm{
  C_0(X) \rtimes^{\ell^1} G \ar@{^{(}->}[r] \ar[dr] & C_0(X) \rtimes G \ar[d] \\
   & B
  }
\egloz
commutative.

By construction, there is a canonical *-homomorphism $C_0(X) \rtimes G \to C_0(X) \rtimes_r G$ extending the identity on $C_0(X) \rtimes^{\ell^1} G$.

The reader may consult \cites{McCl,Ex4} for more information about partial dynamical systems and their C*-algebras.

\subsection{\'{E}tale groupoids}
\label{ss:GPD}

Groupoids play an important role in operator algebras in general and for our topic of semigroup C*-algebras in particular. This is because many C*-algebras can be written as groupoid C*-algebras. This also applies to many semigroup C*-algebras.

Let us first introduce groupoids. In the language of categories, a groupoid is simply a small category with inverses. Very roughly speaking, this means that a groupoid is a group where multiplication is not globally defined. Roughly speaking, a groupoid $\cG$ is a set, whose elements $\gamma$ are arrows $r(\gamma) \longleftarrow s(\gamma)$. Here $r(\gamma)$ and $s(\gamma)$ are elements in $\cG^{(0)}$, the set of units. $r$ stands for range and $s$ stands for source. For every $u \in \cG^{(0)}$, there is a distinguished arrow $u \overset{\id_u}{\longleftarrow} u$ in our groupoid $\cG$. This allows us to define an embedding $$\cG^{(0)} \into \cG, \, u \ma \id_u,$$
which in turn allows us to view $\cG^{(0)}$ as a subset of $\cG$.

$\cG$ comes with a multiplication
$$
  \menge{(\gamma,\eta) \in \cG \times \cG}{s(\gamma) = r(\eta)} \lori \cG, \, (\gamma, \eta) \ma \gamma \eta.
$$
We think of this multiplication as concatenation of arrows. With this picture in mind, the condition $s(\gamma) = r(\eta)$ makes sense. Also, $\cG$ comes with an inversion
$$
  \cG \to \cG, \, \gamma \to \gamma^{-1}.
$$
We think of this inversion as reversing arrows. The picture of arrows, with concatenation as multiplication and reversing as inversion, leads to obvious axioms, which, once imposed, give rise to the formal definition of a groupoid. Let us present the details.

\bdefin
A groupoid is a set $\cG$, together with a bijective map $\cG \to \cG, \, \gamma \ma \gamma^{-1}$, a subset $\cG * \cG \subseteq \cG \times \cG$, and a map $\cG * \cG \to \cG, \, (\gamma,\eta) \ma \gamma \eta$, such that
\bglnoz
  && (\gamma^{-1})^{-1} = \gamma \ {\rm for} \ {\rm all} \ \gamma \in \cG,\\
  && (\gamma \eta) \zeta = \gamma (\eta \zeta) \ {\rm for} \ {\rm all} \ (\gamma,\eta), \, (\eta, \zeta) \in \cG * \cG,\\
  && \gamma^{-1} \gamma \eta = \eta, \, \gamma \eta \eta^{-1} = \gamma \ {\rm for} \ {\rm all} \ (\gamma,\eta) \in \cG * \cG.
\eglnoz
\edefin
Note that we implicitly impose conditions on $\cG * \cG$ so that these equations make sense. For instance, the second equation implicitly requires that for all $(\gamma,\eta)$ and $(\eta, \zeta)$ in $\cG * \cG$, $((\gamma \eta),\zeta)$ and $(\gamma,(\eta \zeta))$ must lie in $\cG * \cG$ as well.

Elements in $\cG * \cG$ are called composable pairs.

The set of units is now defined by
$$
  \cG^{(0)} \defeq \menge{\gamma^{-1} \gamma}{\gamma \in \cG},
$$
it is also given by
$$
  \cG^{(0)} = \menge{\gamma \gamma^{-1}}{\gamma \in \cG}.
$$
Moreover, we define the source map by setting
$$
  s: \: \cG \to \cG^{(0)}, \, \gamma \ma \gamma^{-1} \gamma
$$
and the range map by setting
$$
  r: \: \cG \to \cG^{(0)}, \, \gamma \ma \gamma \gamma^{-1}.
$$
It is now an immediate consequence of the axioms that
$$
  \cG * \cG = \menge{(\gamma,\eta) \in \cG \times \cG}{s(\gamma) = r(\eta)}.
$$

A groupoid $\cG$ is called a topological groupoid if the set $\cG$ comes with a topology such that multiplication and inversion become continuous maps. A topological groupoid is called \'{e}tale if $r$ and $s$ are local homeomorphisms. A topological groupoid is called locally compact if it is locally compact (and Hausdorff) as a topological space.

As an example, let us describe the partial transformation groupoid attached to the partial dynamical system $\alpha: \: G \curvearrowright X$. It is denoted by $G \mathbin{_{\alpha} \ltimes} X$ and is given by
$$G \mathbin{_{\alpha} \ltimes} X \defeq \menge{(g,x) \in G \times X}{g \in G, x \in U_{g^{-1}}},$$
with source map $s(g,x) = x$, range map $r(g,x) = g.x$, composition
$$(g_1,g_2.x)(g_2,x) = (g_1g_2,x)$$
and inverse
$$(g,x)^{-1} = (g^{-1},g.x).$$
We equip $G \mathbin{_{\alpha} \ltimes} X$ with the subspace topology from $G \times X$. Usually, we write $G \ltimes X$ for $G \mathbin{_{\alpha} \ltimes} X$ if the action $\alpha$ is understood. The unit space of $G \ltimes X$ coincides with $X$. Since $G$ is discrete, $G \ltimes X$ is an \'{e}tale groupoid. Actually, if we set
$$
G_x \defeq \menge{g \in G}{x \in U_{g^{-1}}} \ {\rm and} \ G^x \defeq \menge{g \in G}{x \in U_g}
$$
for $x \in X$, then we have canonical identifications
$$
s^{-1}(x) \cong G_x, \, (g,x) \ma g \ {\rm and} \ r^{-1}(x) \cong G^x, \, (g,g^{-1}.x) \ma g.
$$

Let $\cG$ be an \'{e}tale locally compact groupoid. For $x \in \cG^{(0)}$, let $\cG_x = s^{-1}(x)$ and $\cG^x = r^{-1}(x)$. $C_c(\cG)$ is a *-algebra with respect to the multiplication
$$
  (f*g)(\gamma) = \sum_{\beta \in \cG_{s(\gamma)}} f(\gamma \beta^{-1}) g(\beta)
$$
and the involution
$$f^*(\gamma) = \overline{f(\gamma^{-1})}.$$
For every $x \in \cG^{(0)}$, define a *-representation $\pi_x$ of $C_c(\cG)$ on $\ell^2 \cG_x$ by setting
$$
  \pi_x(f)(\xi)(\gamma) = (f * \xi)(\gamma) = \sum_{\beta \in \cG_x} f(\gamma \beta^{-1}) \xi(\beta).
$$
Alternatively, if we want to highlight why these representations play the role of the left regular representation, attached to left multiplication, we could define $\pi_x$ by setting
$$
  \pi_x(f) \delta_{\gamma} = \sum_{\alpha \in \cG_{r(\gamma)}} f(\alpha) \delta_{\alpha \gamma}.
$$
Here $\menge{\delta_{\gamma}}{\gamma \in \cG_x}$ is the canonical orthonormal basis of $\ell^2 \cG_x$.

With these definitions, we are ready to define groupoid C*-algebras.
\bdefin
Let
$$\norm{f}_{C^*_r(\cG)} \defeq \sup_{x \in \cG^{(0)}} \norm{\pi_x(f)}$$
for $f \in C_c(\cG)$.

We define $C^*_r(\cG) \defeq \overline{C_c(\cG)}^{\norm{\cdot}_{C^*_r(\cG)}}$. 

$C^*_r(\cG)$ is called the reduced groupoid C*-algebra of $\cG$.
\edefin
Alternatively, we could set
$$\pi = \bigoplus_{x \in \cG^{(0)}} \pi_x$$
and
$$C^*_r(\cG) = \overline{\pi(C_c(\cG))} \subseteq \cL(\bigoplus_x \ell^2 \cG_x).$$

Let us now define the full groupoid C*-algebra. Let $\cG$ be an \'{e}tale locally compact groupoid. Then $\cG^{(0)}$ is a clopen subspace of $\cG$. Therefore, we can think of $C_c(\cG^{(0)})$ as a subspace of $C_c(\cG)$ simply by extending functions on $\cG^{(0)}$ by $0$ to functions on $\cG$. This allows us to define the full groupoid C*-algebra.

\bdefin
For $f \in C_c(\cG)$, let
$$\norm{f}_{C^*(\cG)} = \sup_{\pi} \norm{\pi(f)},$$
where the supremum is taken over all *-representations of $C_c(\cG)$ which are bounded on $C_c(\cG^{(0)})$ (with respect to the supremum norm $\norm{\cdot}_{\infty}$).

We then set $C^*(\cG) \defeq \overline{C_c(\cG)}^{\norm{\cdot}_{C^*(\cG)}}$.

$C^*(\cG)$ is called the full groupoid C*-algebra of $\cG$.
\edefin

\bremark
We will only deal with second countable locally compact \'{e}tale groupoids. In that case, \cite[Chapter~II, Theorem~1.21]{Ren1} tells us that every *-representation of $C_c(\cG)$ on a separable Hilbert space is automatically bounded. In other words, the full groupoid C*-algebra of $\cG$ is the universal enveloping C*-algebra of $C_c(\cG)$. This notion has been explained after Definition~\ref{DEF:full-cropro-partial}.
\eremark

By construction, there is a canonical *-homomorphism $C^*(\cG) \to C^*_r(\cG)$ extending the identity on $C_c(\cG)$. It is called the left regular representation.

\subsection{The universal groupoid of an inverse semigroup}

\label{ss:univGPD}

We attach groupoids to inverse semigroups so that full and reduced C*-algebras coincide. The groupoids we construct are basically Paterson's universal groupoid, as in \cite[\S~4.3]{Pat99} or \cite{MS}. There is however a small difference. In case of inverse semigroups with zero, our construction differs from Paterson's because we want the distinguished zero element to be represented by zero in the reduced and full C*-algebras.

Let us first explain our construction. We start with an inverse semigroup $S$ with semilattice of idempotents denoted by $E$. Set
$$
  \Sigma \defeq \menge{(s,\chi) \in S \times \widehat{E}}{\chi(s^{-1}s) = 1}.
$$
Note that in case $0 \in S$, we must have $s \neq 0$ since $\chi(0) = 0$ by our convention.

We introduce an equivalence relation on $\Sigma$. Given $(s,\chi)$ and $(t,\psi)$ in $\Sigma$, we define 
$$
(s,\chi) \sim (t,\psi) \ {\rm if} \ {\rm there} \ {\rm exists} \ e \in E \ {\rm with} \ se = te \ {\rm and} \ \chi(e) = 1.
$$
The equivalence class of $(s,\chi) \in \Sigma$ with respect to $\sim$ is denoted by $[s,\chi]$. We set
$$
  \cG(S) \defeq \Sigma / {}_{\sim}, \, {\rm i.e.}, \, \cG(S) = \menge{[s,\chi]}{(s,\chi) \in \Sigma}.
$$
To define a multiplication on $\cG(S)$, we need to introduce the following notation: Let $s \in S$ and $\chi \in \widehat{E}$ be such that $\chi(s^{-1}s) = 1$. Then we define a new element $s.\chi$ of $\widehat{E}$ by setting
$$(s.\chi)(e) \defeq \chi(s^{-1}es).$$
Then we say that $[t,\psi]$ and $[s,\chi]$ are composable if $\psi = s.\chi$. In that case, we define their product as
$$
  [t,\psi] [s,\chi] \defeq [ts,\chi].
$$
The inverse map is given by
$$
  [s,\chi]^{-1} \defeq [s^{-1},s.\chi].
$$
It is easy to see that multiplication and inverse are well-defined, and they give rise to a groupoid structure on $\cG(S)$.

Moreover, we introduce a topology on $\cG(S)$ by choosing a basis of open subsets. Given $s \in S$ and an open subspace 
$$U \subseteq \menge{\chi \in \widehat{E}}{\chi(s^{-1}s) = 1},$$
we define
$$
  D(s,U) \defeq \menge{[s,\chi]}{\chi \in U}.
$$
We equip $\cG(S)$ with the topology which has as a basis of open subsets
$$
  D(s,U), \ {\rm for} \ s \in S \ {\rm and} \ U \subseteq \menge{\chi \in \widehat{E}}{\chi(s^{-1}s) = 1} \ {\rm open}.
$$
It is easy to check that with this topology, $\cG(S)$ becomes a locally compact \'{e}tale groupoid. In all our examples, $S$ will be countable, in which case $\cG(S)$ will be second countable.

Let us explain the difference between our groupoid $\cG(S)$ and the universal groupoid attached to $S$ in \cite[\S~4.3]{Pat99}. Assume that $S$ is an inverse semigroup with zero, and $0$ is the distinguished zero element. The starting point is that our space $\widehat{E}$ and the space of semi-characters $X$ introduced in \cite[\S~2.1]{Pat99} and \cite[\S~4.3]{Pat99} do not coincide. They are related by
$$
  X = \widehat{E} \sqcup \gekl{\chi_0}.
$$
Here $\chi_0$ is the semi-character on $E$ which sends every element of $E$ to $1$, even $0$. The disjoint union above is not only a disjoint union of sets, but also of topological spaces, i.e., $\chi_0$ is an isolated point in $X$ (it is open and closed).

Now it is easy to see that our $\cG(S)$ is the restriction of the universal groupoid $G_{\mathbf{u}}$ attached to $S$ in \cite[\S~4.3]{Pat99} to $\widehat{E}$. This means that
$$
  \cG(S) = \menge{\gamma \in \cG_{\mathbf{u}}}{r(\gamma) \in \widehat{E}, \, s(\gamma) \in \widehat{E}}.
$$
Actually, the only element in $\cG_{\mathbf{u}}$ which does not have range and source in $\widehat{E}$ is $\chi_0$ itself. It follows that
\bgl
\label{GuGS0}
  G_{\mathbf{u}} = \cG(S) \sqcup \gekl{\chi_0}.
\egl

\subsection{Inverse semigroup C*-algebras as groupoid C*-algebras}

We begin by identifying the full C*-algebras. Given an inverse semigroup $S$ with semilattice of idempotents $E$, let us introduce the notation that for $e \in E$, we write 
$$
  U_e \defeq \menge{\chi \in \widehat{E}}{\chi(e) = 1}.
$$

\btheo
\label{THM:C*S=C*GPD}
For every inverse semigroup $S$, there is a canonical isomorphism
$$
  C^*(S) \overset{\cong}{\lori} C^*(\cG(S))
$$
sending the generator $v_s \in C^*(S)$ to the characteristic function on $D(s,U_{s^{-1}s})$, viewed as an element in $C_c(\cG) \subseteq C^*(\cG)$.
\etheo
Recall that
$$
  D(s,U_{s^{-1}s}) = \menge{[s,\chi]}{\chi \in U_{s^{-1}s}}.
$$
\bproof
If case of inverse semigroups without zero, our theorem is just \cite[Chapter~4, Theorem~4.4.1]{Pat99}.

Now let us assume that $0 \in S$. Then the full C*-algebra attached to $S$ in \cite[\S~2.1]{Pat99} is canonically isomorphic to
$$C^*(S) \oplus \Cz v_0,$$
where $C^*(S)$ is our full inverse semigroup C*-algebra in the sense of Definition~\ref{DEF:C*S}, and $v_0$ is a (non-zero) projection.

For the full groupoid C*-algebra of the universal groupoid $G_{\mathbf{u}}$ attached to $S$ in \cite[\S~4.3]{Pat99}, we get because of \eqref{GuGS0}:
$$
  C^*(G_{\mathbf{u}}) \cong C^*(\cG(S)) \oplus \Cz 1_{\chi_0}.
$$
Here $1_{\chi_0}$ is the characteristic function of the one-point set $\gekl{\chi_0}$, and it is easy to see that $1_{\chi_0}$ is a (non-zero) projection.

With these observations in mind, it is easy to see that the identification in \cite[Chapter~4, Theorem~4.4.1]{Pat99} of the full C*-algebra attached to $S$ in \cite[\S~2.1]{Pat99} with the full groupoid C*-algebra $C^*(G_{\mathbf{u}})$ respects these direct sum decompositions, i.e., it sends $C^*(S)$ in the sense of Definition~\ref{DEF:C*S} to $C^*(\cG(S))$. Finally, it is also easy to see that the identification we get in this way really sends $v_s \in C^*(S)$ to the characteristic function on $D(s,U_{s^{-1}s})$.
\eproof

Next, we identify the reduced C*-algebras.
\btheo
\label{THM:C*S=C*GPS_red}
For every inverse semigroup $S$, there is a canonical isomorphism
$$
  C^*_{\lambda}(S) \overset{\cong}{\lori} C^*_r(\cG(S))
$$
sending the generator $\lambda_s \in C^*_{\lambda}(S)$ to the characteristic function on $D(s,U_{s^{-1}s})$, viewed as an element in $C_c(\cG) \subseteq C^*_r(\cG)$.
\etheo
We could give a proof of this result in complete analogy to the case of the full C*-algebras, using \cite[Chapter~4, Theorem~4.4.2]{Pat99} instead of \cite[Chapter~4, Theorem~4.4.1]{Pat99}. Instead, since all these C*-algebras are defined using concrete representations, we give a concrete proof identifying certain representations.
\bproof
For $e \in E\reg$, define
$$
  S\reg_e \defeq \menge{x \in S\reg}{x^{-1}x = e}.
$$
It is then easy to see that
$$
  S\reg = \bigsqcup_{e \in E\reg} S\reg_e.
$$
This yields the direct sum decomposition
$$
  \ell^2 S\reg = \bigoplus_{e \in E\reg} \ell^2 S\reg_e.
$$
The left regular representation of $S$ respects this direct sum decomposition. This is because given $s \in S$ and $x \in S\reg_e$ with $s^{-1}s \geq xx^{-1}$, we have that $sx \in S\reg_e$ since
$$
  (sx)^{-1}(sx) = x^{-1} (s^{-1}s) x = x^{-1} (s^{-1}s xx^{-1}) x = x^{-1} (xx^{-1}) x = x^{-1}x = e.
$$
Therefore, for every $s \in S$, we have
$$
  \lambda_s = \bigoplus_{e \in E\reg} \lambda_s \big\vert_{\ell^2 S\reg_e}.
$$
Now define for every $e \in E\reg$ the character $\chi_e \in \widehat{E}$ by setting
\bglnoz
  \chi_e(f) &=& 1 \ {\rm if} \ e \leq f,\\
  \chi_e(f) &=& 0 \ {\rm if} \ e \nleq f.
\eglnoz
The map
$$
  S\reg_e \lori \cG(S)_{\chi_e}, \, x \ma [x,\chi_e]
$$
is surjective as every $(x,\chi_e) \in \Sigma$ is equivalent to $(xe,\chi_e)$, and $xe$ lies in $S\reg_e$ as $\chi_e(x^{-1}x) = 1$ implies $e \leq x^{-1}x$. It is also injective as $[x,\chi_e] = [y,\chi_e]$ for $x,y \in S\reg_e$ implies that $xf = yf$ for some $f \in E\reg$ with $e \leq f$, and thus $x = y$. Therefore, the map above is a bijection. It induces a unitary
$$
  U: \: \ell^2 S\reg_e \overset{\cong}{\lori} \ell^2 \cG(S)_{\chi_e}, \, \delta_x \ma \delta_{[x,\chi_e]}.
$$
Now let $1_{D(s,U_{s^{-1}s})}$ be the characteristic function on $D(s,U_{s^{-1}s})$, viewed as an element in $C_c(\cG)$. Then we have
\bgl
\label{Ul=piU}
  U \circ \lambda_s \big\vert_{\ell^2 S\reg_e} = \pi_{\chi_e}(1_{D(s,U_{s^{-1}s})}) \circ U.
\egl
This is because
$$
  (U \circ \lambda_s \big\vert_{\ell^2 S\reg_e})(\delta_x) = U(\delta_{sx}) = [sx,\chi_e]
$$
and
$$
  (\pi_{\chi_e} \circ 1_{D(s,U_{s^{-1}s})} \circ U)(\delta_x) = \pi_{\chi_e}(1_{D(s,U_{s^{-1}s})})([x,\chi_e]) = [sx,\chi_e]
$$
if $s^{-1}s \geq xx^{-1}$, and both sides of \eqref{Ul=piU} are zero if $s^{-1}s \ngeq xx^{-1}$.

Hence it follows that the left regular representation of $C^*(S)$ is unitarily equivalent to
$$
  \bigoplus_{e \in E\reg} \pi_{\chi_e}
$$
under the isomorphism from Theorem~\ref{THM:C*S=C*GPD}.

Thus, all we have to show in order to conclude our proof is that
\bgl
\label{sup=sup}
  \sup_{\chi \in \widehat{E}} \norm{\pi_{\chi}(f)} = \sup_{e \in E\reg} \norm{\pi_{\chi_e}(f)},
\egl
for all $f \in C_c(\cG(S))$. To show this, we first need to observe that
$$
  \menge{\chi_e}{e \in E\reg}
$$
is dense in $\widehat{E}$. This is because a basis of open subsets for the topology of $\widehat{E}$ are given by
$$
  U(e;e_1,\dotsc,e_n) \defeq \menge{\chi \in \widehat{E}}{\chi(e) = 1; \; \chi(e_1) = \dotso = \chi(e_n) = 0},
$$
for $e, e_1, \dotsc, e_n \in E\reg$ with $e_i \nleq e$. It is then clear that $\chi_e$ lies in $U(e;e_1,\dotsc,e_n)$. 

Because of density, \eqref{sup=sup} follows from \cite[Chapter~3, Proposition~3.1.2]{Pat99}.
\eproof

\bremark
\label{isom-leftreg}
It is clear that the explicit isomorphisms provided by Theorem~\ref{THM:C*S=C*GPD} and Theorem~\ref{THM:C*S=C*GPS_red} give rise to a commutative diagram
\bgloz
  \xymatrix@C=20mm{
  C^*(S) \ar[r] \ar[d]^{\cong} & C^*_{\lambda}(S) \ar[d]^{\cong} \\
  C^*(\cG(S)) \ar[r] & C^*_r(\cG(S))
  }
\egloz
where the horizontal arrows are the left regular representations and the vertical arrows are the identifications provided by Theorem~\ref{THM:C*S=C*GPD} and Theorem~\ref{THM:C*S=C*GPS_red}.
\eremark

\subsection{C*-algebras of partial dynamical systems as C*-algebras of partial transformation groupoids}

Our goal is to identify the full and reduced crossed products attached to partial dynamical systems with full and reduced groupoid C*-algebras for the corresponding partial transformation groupoids.

Given a partial dynamical system $G \curvearrowright X$, we have constructed its partial transformation groupoid $G \ltimes X$ in \S~\ref{ss:GPD}. 

The following result is \cite[Theorem~3.3]{Aba2}:
\btheo
\label{THM:C*GPD=C*PDS}
The canonical homomorphism
$$
  C_c(G \ltimes X) \to C_0(X) \rtimes^{\ell^1} G, \, \theta \ma \sum_g \theta(g,g^{-1}.\sqcup) \delta_g,
$$
where $\theta(g,g^{-1}.\sqcup)$ is the function $U_{g^{-1}} \to \Cz, \, x \ma \theta(g,g^{-1}.x)$, extends to an isomorphism
$$C^*(G \ltimes X) \overset{\cong}{\lori} C_0(X) \rtimes G.$$
\etheo
Here we use the same notation for partial dynamical systems and their crossed products as in \S~\ref{ss:PDS}.

Let us now identify reduced crossed products.
\btheo
\label{THM:C*GPD=C*PDS_red}
The canonical homomorphism
\bgl
\label{cG->}
  C_c(G \ltimes X) \to C_0(X) \rtimes^{\ell^1} G, \, \theta \ma \sum_g \theta(g,g^{-1}.\sqcup) \delta_g,
\egl
where $\theta(g,g^{-1}.\sqcup)$ is the function
$$U_{g^{-1}} \to \Cz, \, x \ma \theta(g,g^{-1}.x),$$
extends to an isomorphism
$$C^*_r(G \ltimes X) \overset{\cong}{\lori} C_0(X) \rtimes_r G.$$
\etheo
We include a proof of this result. It is taken from \cite{Li-PTGPD}.
\bproof
We use the same notation as in the construction of the reduced crossed product in \S~\ref{ss:PDS}. As above, let $\mu \times \lambda$ be the representation $C_0(X) \rtimes^{\ell^1} G \to \cL(H)$ which we used to define $C_0(X) \rtimes_r G$. Our first observation is
\bgl
\label{nondeg}
  \overline{\img(\mu \times \lambda)(H)}
  = \bigoplus_{h \in G} \delta_h \otimes \ell^2 U_{h^{-1}}.
\egl
To see this, observe that for all $g \in G$,
$$\img(E_g) \subseteq \bigoplus_h \delta_h \otimes \ell^2(U_{h^{-1}} \cap U_{(gh)^{-1}}).$$
This holds since for
$$x \notin h^{-1}.(U_h \cap U_{g^{-1}}) = U_{(gh)^{-1}} \cap U_{h^{-1}},$$
$f \vert_{U_h} (h.x) = 0$ for $f \in C_0(U_{g^{-1}})$. Therefore,
$$\pi(C_0(U_{g^{-1}}))(\delta_h \otimes \ell^2 X) \subseteq \delta_h \otimes \ell^2(U_{h^{-1}} \cap U_{(gh)^{-1}}).$$
Hence
$$\img(E_g) \subseteq \bigoplus_h \delta_h \otimes \ell^2(U_{h^{-1}} \cap U_{(gh)^{-1}}),$$
and thus,
$$\img(V_g) \subseteq \bigoplus_h \delta_{gh} \otimes \ell^2(U_{h^{-1}} \cap U_{(gh)^{-1}}) \subseteq \bigoplus_h \delta_h \otimes \ell^2 U_{h^{-1}}.$$
This shows \an{$\subseteq$} in \eqref{nondeg}. For \an{$\supseteq$}, note that for $f \in C_0(X)$,
$$(\mu \times \lambda)(f \delta_e) = \mu(f) E_e,$$
and for $\xi \in \ell^2 U_{h^{-1}}$, 
$$\mu(f) E_e (\delta_h \otimes \xi) = \delta_h \otimes f \vert_{U_h}(h. \sqcup) \xi.$$
So $(\mu \times \lambda)(f \delta_e)(H)$ contains $\delta_h \otimes f \cdot \xi$ for all $f \in C_0(U_{h^{-1}})$ and $\xi \in \ell^2 U_{h^{-1}}$, hence also $\delta_h \otimes \ell^2 U_{h^{-1}}$. This proves \an{$\supseteq$}.

For $x \in X$, let $G_x = \menge{g \in G}{x \in U_{g^{-1}}}$ as before. Our second observation is that for every $x \in X$, the subspace $H_x \defeq \ell^2 G_x \otimes \delta_x$ is $(\mu \times \lambda)$-invariant. It is clear that $\mu(f)$ leaves $H_x$ invariant for all $f \in C_0(X)$. For $g, h \in G$,
$$
E_g(\delta_h \otimes \delta_x) = \delta_h \otimes \delta_x
$$
if $x \in U_{h^{-1}} \cap U_{(gh)^{-1}}$, and if that is the case, then
$$V_g(\delta_h \otimes \delta_x) = \delta_{gh} \otimes \delta_x \in H_x.$$

Therefore,
$$
H = \rukl{\bigoplus_{x \in X} H_x} \oplus (\mu \times \lambda)(C_0(X) \rtimes^{\ell^1} G)(H)^{\perp}
$$
is a decomposition of $H$ into $\mu \times \lambda$-invariant subspaces. For $x \in X$, set
$$\rho_x \defeq (\mu \times \lambda) \vert_{H_x}.$$
Then
$$C_0(X) \rtimes_r G = \overline{C_0(X) \rtimes^{\ell^1} G}^{\norm{\cdot}_{\bigoplus_x \rho_x}}.$$

Moreover, we have for $x \in U_{h^{-1}}$,
\bgln
  && \rho_x \rukl{\sum_g f_g \delta_g}(\delta_h \otimes \delta_x)
  = \sum_g \mu(f_g) V_g (\delta_h \otimes \delta_x) \nonumber \\
  &=& \sum_{g: \: x \in U_{(gh)^{-1}}} \mu(f_g) (\delta_{gh} \otimes \delta_x)
  = \sum_{g: \: x \in U_{(gh)^{-1}}} \delta_{gh} \otimes f_g(gh.x) \delta_x \nonumber \\
\label{rho_x}
  &=& \sum_{k \in G_x} \delta_k \otimes f_{kh^{-1}}(k.x) \delta_x. 
\egln
Let us compare this construction with the construction of the reduced groupoid C*-algebra of $G \ltimes X$. Obviously, \eqref{cG->} is an embedding of $C_c(G \ltimes X)$ as a subalgebra which is $\norm{\cdot}_{\ell^1}$-dense in $C_0(X) \rtimes^{\ell^1} G$. Therefore,
$$C_0(X) \rtimes_r G = \overline{C_c(G \ltimes X)}^{\norm{\cdot}_{\bigoplus_x \rho_x}}.$$

Now, to construct the reduced groupoid C*-algebra $C^*_r(G \ltimes X)$, we follow our explanations in \S~\ref{ss:GPD} and construct for every $x \in X$ the representation
$$\pi_x: \: C_c(G \ltimes X) \to \cL(\ell^2(s^{-1}(x)))$$
by setting
$$\pi_x(\theta)(\xi)(\zeta) \defeq \sum_{\eta \, \in \, s^{-1}(x)} \theta(\zeta \eta^{-1}) \xi(\eta).$$
In our case, using $s^{-1}(x) = G_x \times \gekl{x}$, we obtain for $\xi = \delta_h \otimes \delta_x$ with $h \in G_x$:
$$
\pi_x(\theta)(\delta_h \otimes \delta_x)(k,x) = \theta((k.x)(h,x)^{-1}) = \theta(kh^{-1},h.x).$$
Thus,
\bgl
\label{pi_x}
  \pi_x(\theta)(\delta_h \otimes \delta_x)(k,x) = \sum_{k \in G_x} \theta(kh^{-1},h.x) \delta_k \otimes \delta_x.
\egl
By definition,
$$C^*_r(G \ltimes X) = \overline{C_c(G \ltimes X)}^{\norm{\cdot}_{\bigoplus_x \pi_x}}.$$
Therefore, in order to show that $\norm{\cdot}_{\bigoplus_x \rho_x}$ and $\norm{\cdot}_{\bigoplus_x \rho_x}$ coincide on $C_c(G \ltimes X)$, it suffices to show that for every $x \in X$, $\pi_x$ and the restriction of $\rho_x$ to $C_c(G \ltimes X)$ are unitarily equivalent. Given $x \in X$, using $s^{-1}(x) = G_x \times \gekl{x}$, we obtain the canonical unitary
$$\ell^2(s^{-1}(x)) \cong H_x = \ell^2(G_x) \otimes \delta_x,$$
so that we may think of both $\rho_x$ and $\pi_x$ as representations on $\ell^2(G_x) \otimes \delta_x$. We then have for $x \in X$, $\theta \in C_c(G \ltimes X)$ and $h \in G_x$:
\bglnoz
  \rho_x(\theta)(\delta_h \otimes \delta_x)
  &\overset{\eqref{cG->}}{=}&
  \rho_x(\sum_g \theta(g,g^{-1}. \sqcup) \delta_g)(\delta_h \otimes \delta_x)\\
  &\overset{\eqref{rho_x}}{=}&
  \sum_{k \in G_x} \delta_k \otimes \theta(kh^{-1},h.x) \delta_x
  \overset{\eqref{pi_x}}{=} \pi_x(\theta)(\delta_h \otimes \delta_x).
\eglnoz
This yields the canonical identification
$$C_0(X) \rtimes_r G \cong C^*_r(G \ltimes X),$$
as desired.
\eproof

\subsection{The case of inverse semigroups admitting an idempotent pure partial homomorphism to a group}

We would like to show that in the case of inverse semigroups which admit an idempotent pure partial homomorphism to a group, all our constructions above coincide.

Let $S$ be an inverse semigroup and $E$ the semilattice of idempotents of $S$. Let $G$ be a group. Assume that $\sigma$ is a partial homomorphism $S\reg \to G$ which is idempotent pure.

In this situation, we constructed a partial dynamical system $G \curvearrowright \widehat{E}$ in \S~\ref{ss:PDS}. Our first observation is that the partial transformation groupoid of $G \curvearrowright \widehat{E}$ can be canonically identified with the groupoid $\cG(S)$ we attached to $S$ in \S~\ref{ss:univGPD}.

\blemma
\label{GS=GxhatE}
In the situation described above, we have a canonical identification
$$
  \cG(S) \overset{\cong}{\lori} G \ltimes \widehat{E}, \, [s,\chi] \ma (\sigma(s),\chi).
$$
of topological groupoids.
\elemma
\bproof
We use the notations from \S~\ref{ss:PDS} and \S~\ref{ss:univGPD}.

To see that the mapping $[s,\chi] \ma (\sigma(s),\chi)$ is well-defined, suppose that $(s,\chi)$ and $(t,\chi)$ in $\Sigma$ are equivalent. Then there exists $e \in E\reg$ such that $se = te$, and $se$ (or $te$) cannot be zero in case $0 \in S$. Therefore,
$$
  \sigma(s) = \sigma(se) = \sigma(te) = \sigma(t).
$$
To see that $[s,\chi] \ma (\sigma(s),\chi)$ is a morphism of groupoids, note that $[s,\chi]^{-1} = [s^{-1},s.\chi]$ is sent to $(\sigma(s^{-1}),s.\chi) = (\sigma(s),\chi)^{-1}$. Hence our mapping respects inverses. For multiplication, observe that
$$
  s([s,\chi]) = \chi = s(\sigma(s),\chi)
$$
and
$$
  r([s,\chi]) = s.\chi = \sigma(s).\chi = r(\sigma(s),\chi).
$$
Moreover, $[t,s.\chi] \cdot [s,\chi] = [ts,\chi]$ is mapped to $[\sigma(ts),\chi] = [\sigma(t),s.\chi] \cdot [\sigma(s),\chi]$. Hence it follows that our mapping is a groupoid morphism.

We now set out to construct an inverse. Define the map  
$$
  G \ltimes \widehat{E} \lori \cG(S), \, (g,\chi) \ma [s,\chi]
$$
where for every $g$ in $G$, we choose $s \in S$ with $\sigma(s) = g$ and $\chi(s^{-1}s) = 1$. This is well-defined: Given $t \in S$ with $\sigma(t) = g$ and $\chi(t^{-1}t) = 1$, set $e \defeq s^{-1}s t^{-1}t$. Then $\chi(e) = 1$. Moreover, $se = s t^{-1}t$ and $te = t s^{-1}s$. As $\sigma(se) = \sigma(s) = g = \sigma(t) = \sigma(te)$ and $(se)^{-1}(se) = e = (te)^{-1}(te)$, we deduce by Lemma~\ref{dom-sigma==} that $se = te$. Hence $(s,\chi) \sim (t,\chi)$.

It is easy to see that we have just constructed the inverse of 
$$
  \cG(S) \lori G \ltimes \widehat{E}, \, [s,\chi] \ma (\sigma(s),\chi).
$$
Moreover, it is also easy to see that both our mappings are open, so that they give rise to the desired identification of topological groupoids.
\eproof

Combining Theorem~\ref{THM:C*S=C*GPD} with Theorem~\ref{THM:C*GPD=C*PDS} and Theorem~\ref{THM:C*S=C*GPS_red} with Theorem~\ref{THM:C*GPD=C*PDS_red}, we obtain the following
\bcor
\label{C*S=C*ExG}
Let $S$ be an inverse semigroup and $E$ the semilattice of idempotents of $S$. Let $G$ be a group. Assume that $\sigma$ is a partial homomorphism $S\reg \to G$ which is idempotent pure.

In this situation, we have canonical isomorphisms
$$C^*(S) \to C^*(E) \rtimes G, \, v_s \ma \lambda_{ss^{-1}} \delta_{\sigma(s)}$$
and
$$C^*_{\lambda}(S) \to C^*(E) \rtimes_r G, \, \lambda_s \ma \lambda_{ss^{-1}} V_{\sigma(s)}.$$
\ecor

\section{Amenability and nuclearity}

Amenability is an important structural property for groups and groupoids, while nuclearity plays a crucial role in the structure theory for C*-algebras, in particular in the classification program. In the case of groups and groupoids, it is known that amenability and nuclearity of C*-algebras are closely related. Moreover, there are further alternative ways to characterize amenability in terms of C*-algebras. Our goal now is to explain to what extent analogous results hold true in the semigroup context.

\subsection{Groups and groupoids}

Let us start by reviewing the case of groups and groupoids.

Let $G$ be a discrete group. We recall three conditions.

\bdefin
Our group $G$ is said to be amenable if there exists a left invariant state on $\ell^{\infty}(G)$.
\edefin
This means that we require the existence of a state $\mu: \: \ell^{\infty}(G) \to \Cz$ with the property that $\mu(f(s \sqcup)) = \mu(f)$ for every $f \in \ell^{\infty}(G)$ and $s \in G$. Here $f(s \sqcup)$ is the function $G \to \Cz, \, x \ma f(sx)$.

\bdefin
Our group $G$ is said to satisfy Reiter's condition if there exists a net $(\theta_i)_i$ of probability measures on $G$ such that 
$$\lim_{i \to \infty} \norm{\theta_i - g \theta_i} = 0$$ 
for all $g \in G$.
\edefin
Here $g \theta$ is the pushforward of $\theta$ under
$$G \cong G, \, x \ma gx.$$

\bdefin
Our group $G$ is said to satisfy F{\o}lner's condition if for every finite subset $E \subseteq G$ and every $\varepsilon > 0$, there exists a non-empty finite subset $F \subseteq G$ with
$$
  \abs{(sF) \triangle F} / \abs{F} < \varepsilon
$$
for all $s \in E$.
\edefin
Here $sF = \menge{sx}{x \in F}$, and $\triangle$ stands for symmetric difference.

It turns out that a group is amenable if and only if it satisfies Reiter's condition if and only if it satisfies F{\o}lner's condition. We refer the reader to \cite[Chapter~2, \S~6]{BO} for more details.

All abelian, nilpotent and solvable groups are amenable, to mention some examples. Non-abelian free groups are not amenable.

We now turn to groupoids.

\bdefin
An \'{e}tale locally compact groupoid $\cG$ is amenable if there is a net $(\theta_i)_i$ of continuous systems of probability measures $\theta_i = (\theta_i^x)_{x \in \cG^{(0)}}$ with
$$\lim_{i \to \infty} \norm{\theta_i^{r(\gamma)} - \gamma \theta_i^{s(\gamma)}} = 0 \ {\rm for} \ {\rm all} \ \gamma \in \cG. $$
\edefin
Here $\theta^x$ is a probability measure on $\cG$ with support contained in $\cG^x$. \an{Continuous} means that for every $f \in C_c(\cG)$, the function
$$\cG^{(0)} \to \Cz, \, x \ma \int f {\rm d} \theta^x$$
is continuous. As above, $\gamma \theta$ is the pushforward of $\theta$ under
$$\cG^{s(\gamma)} \to \cG^{r(\gamma)}, \, \eta \ma \gamma \eta.$$
Note that what we call amenability of groupoids is really Reiter's condition for groupoids. Moreover, we may require that the convergence in our definition happens uniform on compact subsets of $\cG$. This is because of \cite{Rnew}.

For instance, if $G$ is an amenable group, and $G \curvearrowright \Omega$ is a partial dynamical system on a locally compact Hausdorff space $\Omega$, then the partial transformation groupoid $G \ltimes \Omega$ is amenable by \cite[Theorem~20.7 and Theorem~25.10]{Ex4}. But we can get amenable partial transformation groupoids even if $G$ is not amenable.

Let us now introduce nuclearity for C*-algebras.
\bdefin
A C*-algebra $A$ is nuclear if there exists a net of contractive completely positive maps $\varphi_i: \: A \to F_i$ and $\psi_i: \: F_i \to A$, where $F_i$ are finite dimensional C*-algebras, such that
$$\lim_{i \to \infty} \norm{\psi_i \circ \varphi_i (a) - a} = 0$$
for all $a \in A$.
\edefin
For instance, all commutative C*-algebras are nuclear, and all finite dimensional C*-algebras are nuclear.

The reader may find more about nuclearity for C*-algebras for example in \cite[Chapter~2]{BO}.

Let us now relate amenability and nuclearity. Let us start with the case of groups.

Recall that the full group C*-algebra $C^*(G)$ of a discrete group $G$ is the C*-algebra universal for unitary representations of $G$. This means that $C^*(G)$ is generated by unitaries $u_g$, $g \in G$, satisfying
$$
  u_{gh} = u_g u_h \ {\rm for} \ {\rm all} \ g, h \in G,
$$
and whenever we find unitaries $v_g$, $g \in G$, in another C*-algebra $B$ satisfying
$$
  v_{gh} = v_g v_h \ {\rm for} \ {\rm all} \ g, h \in G,
$$
then there exists a (unique) *-homomorphism $C^*(G) \to B$ sending $u_g$ to $v_g$. 

The reduced group C*-algebra $C^*_{\lambda}(G)$ of a discrete group $G$ is the C*-algebra generated by the left regular representations of $G$. The left regular representation is exactly what we get when we apply the construction at the beginning of \S~\ref{leftregrep} to $G$. Therefore, $C^*_{\lambda}(G)$ is the C*-algebra we get when we apply Definition~\ref{sgpCSTAR} to $G$ in place of $P$.

By construction, we have a canonical *-homomorphism
$$
  \lambda: \: C^*(G) \to C^*_{\lambda}(G), \, u_g \to \lambda_g.
$$
It is called the left regular representation (of $C^*(G)$).

Here are a couple of C*-algebraic characterizations of amenability for groups. We refer the reader to \cite[Chapter~2, \S~6]{BO} for details and proofs.
\btheo
\label{gp-am-nuc}
Let $G$ be a discrete group. The following are equivalent:
\begin{itemize}
\item $G$ is amenable.
\item $C^*(G)$ is nuclear.
\item $C^*_{\lambda}(G)$ is nuclear.
\item The left regular representation $\lambda: \: C^*(G) \to C^*_{\lambda}(G)$ is an isomorphism.
\item There exists a character on $C^*_{\lambda}(G)$.
\end{itemize}
\etheo
Here, by a character on a unital C*-algebra $A$, we simply mean a unital *-homomorphism from $A$ to $\Cz$.

We now turn to groupoids and C*-algebraic characterizations of amenability for them. We already introduced full and reduced groupoid C*-algebras in \S~\ref{ss:GPD}. We also introduced the left regular representation (of the full groupoid C*-algebra)
$$\lambda: \: C^*(\cG) \to C^*_r(\cG).$$

\btheo
\label{gpd-am-nuc}
Let $\cG$ be an \'{e}tale locally compact groupoid. Consider the statements
\begin{enumerate}
\item[(i)] $\cG$ is amenable.
\item[(ii)] $C^*(\cG)$ is nuclear.
\item[(iii)] $C^*_{\lambda}(\cG)$ is nuclear.
\item[(iv)] $\lambda: \: C^*(\cG) \to C^*_{\lambda}(\cG)$ is an isomorphism.
\end{enumerate}
Then (i) $\LRarr$ (ii) $\LRarr$ (iii) $\Rarr$ (iv).
\etheo

We refer to \cite[Chapter~5, \S~6]{BO} and \cite{ADR} for more details.

It was an open question whether statement (iv) implies the other statements. But Rufus Willett gave a counterexample in \cite{Wil}. There are, however, results saying that statement (iv) does imply the other statements for particular classes of groupoids. For instance, we mention \cite{Mat}.

\subsection{Amenability for semigroups}
\label{ss:am_P}

Let us now turn to amenability for semigroups. As in the group case, we have the following definitions:
\bdefin
A discrete semigroup $P$ is called left amenable if there exists a left invariant mean on $\ell^{\infty}(P)$, i.e. a state $\mu$ on $\ell^{\infty}(P)$ such that for every $p \in P$ and $f \in \ell^{\infty}(P)$, $\mu(f(p \sqcup)) = \mu(f)$.
\edefin
Here $f(p \sqcup)$ is the function $P \to \Cz, \, x \ma f(px)$.

For instance, every abelian semigroup is left amenable.

\bdefin
A discrete semigroup $P$ is said to satisfy Reiter's condition if there is a net $(\theta_i)_i$ of probability measures on $P$ with the property that
\bgloz
  \lim_i \norm{\theta_i - p \theta_i} = 0 \ {\rm for} \ {\rm all} \ p \in P.
\egloz
\edefin
Here $p \theta$ is the pushforward of $\theta$ under $P \to P, \, x \ma px$.

\bdefin
A discrete semigroup $P$ satisfies the strong F{\o}lner condition if for every finite subset $E \subseteq P$ and every $\varepsilon > 0$, there exists a non-empty finite subset $F \subseteq P$ such that
$$\abs{(pF) \triangle F} / \abs{F} < \varepsilon$$
for all $p \in C$.
\edefin
Here $pF = \menge{px}{x \in F}$ and $\triangle$ stands for symmetric difference.

As in the group case, a discrete left cancellative semigroup is left amenable if and only if it satisfies Reiter's condition if and only if it satisfies the strong F{\o}lner condition. The reader may consult \cite{LiSG} for a proof, and we also refer to \cite{Pat} for more details.

Our goal now is to find the analogues of Theorem~\ref{gp-am-nuc} and Theorem~\ref{gpd-am-nuc} in the context of semigroups and their C*-algebras. The motivation is to understand and explain -- in a conceptual way -- the following two observations:

Let $P = \Nz \times \Nz$, the universal monoid generated by two commuting elements. This is an abelian semigroup, so it is left amenable. So far, we have not discussed the question how to construct full semigroup C*-algebras. But a natural candidate for the full semigroup C*-algebra of $\Nz \times \Nz$ would be 
$$
  C^* \rukl{v_a, v_b \ \vline \ v_a^* v_a = 1, \ v_b^* v_b = 1, \ v_a v_b = v_b v_a}.
$$
In other words, this is the universal C*-algebra generated by two commmuting isometries. It is the C*-algebra universal for isometric representations of our semigroup. This is a very natural candidate for the full semigroup C*-algebra. But Murphy showed that this C*-algebra is not nuclear in \cite[Theorem~6.2]{Mur4}. 

Next, consider $P = \Nz * \Nz$, the non-abelian free monoid on two generators. As in the group case, non-abelian free semigroups are examples of semigroups which are not left amenable. But it is easy to see that $C^*_{\lambda}(\Nz * \Nz)$ is generated as a C*-algebra by two isometries $V_a$ and $V_b$ with orthogonal range projections, i.e.,
$$
  (V_a V_a^*) \cdot (V_b V_b^*) = 0.
$$ 
Therefore, $C^*_{\lambda}(\Nz * \Nz)$ is isomorphic to the canonical extension of the Cuntz algebra $\cO_2$, as introduced in \cite[\S~3]{CuCMP}. It fits into an exact sequence
$$
  0 \to \cK \to C^*_{\lambda}(\Nz * \Nz) \to \cO_2 \to 0,
$$
where $\cK$ is the C*-algebra of compact operators on a infinite dimensional and separable Hilbert space. Hence it follows that $C^*_{\lambda}(\Nz * \Nz)$ is nuclear. Moreover, $C^*_{\lambda}(\Nz * \Nz)$ can be described as a universal C*-algebra, because 
$$
  C^*_{\lambda}(\Nz * \Nz) \cong C^* \rukl{v_a, v_b \ \vline \ v_a^* v_a = 1, \ v_b^* v_b = 1, \ v_a v_a^* v_b v_b^* = 0}.
$$
So this is a hint that for the semigroup $\Nz * \Nz$, the full and reduced semigroup C*-algebras are isomorphic. But, as we remarked above, $\Nz * \Nz$ is not left amenable.

Our goal now is to explain these phenomena, to clarify the relation between amenability and nuclearity, and to obtain analogues of Theorem~\ref{gp-am-nuc} and Theorem~\ref{gpd-am-nuc} in the context of semigroups. The first step for us will be to find a systematic and reasonable way to define full semigroup C*-algebras. It turns out that left inverse hulls attached to left cancellative semigroups, as introduced in \S~\ref{ss:ISGP}, give rise to an approach to this problem. However, before we come to the construction of full semigroup C*-algebras, we first need to compare the reduced C*-algebras of left cancellative semigroups and their left inverse hulls.

\subsection{Comparing reduced C*-algebras for left cancellative semigroups and their left inverse hulls}
\label{Comparing-I-P}

Let $P$ be a left cancellative semigroup and $I_l(P)$ the left inverse hull attached to $P$, as in \S~\ref{ss:ISGP}. As we explained in \S~\ref{ss:ISGP}, we have a canonical embedding of $P$ into $I_l(P)$, denoted by
$$P \into I_l(P), \, p \ma p.$$
It gives rise to the isometry 
$$\bI: \: \ell^2 P \to \ell^2 S\reg, \, \delta_p \ma \delta_p.$$
Thus, we may think of $\ell^2 P$ as a subspace of $\ell^2 S\reg$.

The following observation appears in \cite[\S~3.2]{Nor1}.
\blemma
\label{S-->P}
Assume that $P$ is a left cancellative semigroup with left inverse hull $I_l(P)$. Then the subspace $\ell^2 P$ of $\ell^2 I_l(P)\reg$ is invariant under $C^*_{\lambda}(I_l(P))$. Moreover, we obtain a well-defined surjective *-homomorphism
$$
  C^*_{\lambda}(I_l(P)) \to C^*_{\lambda}(P), \, T \ma \bI^* T \bI
$$
sending $\lambda_p$ to $V_p$ for every $p \in P$.
\elemma
\bproof
We first claim that every $s \in I_l(P)$ has the following property: 
\bgl
\label{s(xr)=s(x)r}
{\rm For} \ {\rm every} \ x \in \dom(s) \ {\rm and} \ {\rm every} \ r \in P, \ xr \ {\rm lies} \ {\rm in} \ \dom(s), \ {\rm and} \ s(xr) = s(x)r.
\egl

To prove our claim, first observe that for every $p \in P$, the partial bijection $p \in I_l(P)$ certainly has this property, as it is just given by left multiplication with $p$. Moreover, $p^{-1}$ is the partial bijection
$$pP \to P, \, px \ma x.$$
Certainly, for every $px \in pP$ and every $r \in P$, $pxr$ lies in $pP$, and
$$p^{-1}(pxr) = xr = p^{-1}(px)r.$$
Hence $p^{-1}$ has the desired property as well. To conclude the proof of our claim, suppose that $s, t \in I_l(P)$ both have the desired property. Choose $x \in \dom(st)$. Then for every $r \in P$, $xr$ lies in $\dom(t)$, and $t(xr) = t(x)r$. Since $t(x)$ lies in $\dom(s)$, $t(x)r$ lies in $\dom(s)$ as well. The conclusion is that $xr$ lies in $\dom(st)$, and we have
$$
  (st)(xr) = s(t(x)r) = s(t(x))r = (st)(x)r.
$$
As every element in $I_l(P)$ is a finite product of partial bijections in 
$$\menge{p}{p \in P} \cup \menge{p^{-1}}{p \in P},$$ 
this proves our claim.

The second step is to show that for every $s \in I_l(P)$ and $x \in P$ with $s^{-1}s \geq pp^{-1}$, we must have $sx = s(x) \in P$. This is because we have, for every $y \in P$:
$$
  (sx)(y) = s(x(y)) = s(xy) = s(x)y = (s(x))(y).
$$
Here we used our first claim from above.

Now let $s \in I_l(P)$ be arbitrary. We want to show that $\lambda_s(\ell^2 P) \subseteq \ell^2 P$. Given $x \in P$, we have $\lambda_s(\delta_x) = 0$ if $s^{-1}s \ngeq pp^{-1}$. If $s^{-1}s \geq pp^{-1}$, then what we showed in the second step implies that $\lambda_s(\delta_x) = \delta_{s(x)}$ lies in $\ell^2 P$. As $s$ was arbitrary, this shows that
$$C^*_{\lambda}(I_l(P))(\ell^2 P) \subseteq \ell^2 P.$$

Therefore, every $T \in C^*_{\lambda}(I_l(P))$ satisfies $T \bI \bI^* = \bI \bI^* T \bI \bI^*$, and since $C^*_{\lambda}(I_l(P))$ is *-invariant, we even obtain that every $T \in C^*_{\lambda}(I_l(P))$ satisfies $T \bI \bI^* = \bI \bI^* T$. This shows that the map
$$
  C^*_{\lambda}(I_l(P)) \to \cL(\ell^2 P), \, T \ma \bI^* T \bI
$$
is a *-homomorphism. Its image is $C^*_{\lambda}(P)$ because we have, for $p \in P$ and $x \in P$:
$$
  \lambda_p(\delta_x) = \delta_{px} = V_p(\delta_x),
$$
so that $\bI^* \lambda_p \bI = V_p$ for all $p \in P$.
\eproof

Recall that we denote the semilattice of idempotents in $I_l(P)$ by $\cJ_P$, and we identified this semilattice with the constructible right ideals of $P$ (see \S~\ref{ss:ISGP}). Moreover, we also introduced in \S~\ref{ss:PDS} the sub-C*-algebra of $C^*_{\lambda}(I_l(P))$ generated by $\cJ_P$:
$$
  C^*(\cJ_P) = C^*(\menge{\lambda_X}{X \in \cJ_P}).
$$
It is easy to see that for every $X \in \cJ_P$, we get
$$
  \bI^* \lambda_X \bI = 1_X,
$$
where $1_X$ is the characteristic function of $X$, viewed as an element in $\ell^{\infty}(P)$. 

Hence, restricting the *-homomorphism
$$
  C^*_{\lambda}(I_l(P)) \to C^*_{\lambda}(P)
$$
from Lemma~\ref{S-->P} to $C^*(\cJ_P)$, we obtain a *-homomorphism from $C^*(\cJ_P)$ onto the sub-C*-algebra $D_{\lambda}(P) = C^*(\menge{1_X}{X \in \cJ_P})$ of $C^*_{\lambda}(P)$, which is generated by $\menge{1_X}{X \in \cJ_P}$,
$$
  C^*(\cJ_P) \onto D_{\lambda}(P), \, T \ma \bI^* T \bI.
$$
Obviously, if the *-homomorphism from Lemma~\ref{S-->P} is an isomorphism, then its restriction to $C^*(\cJ_P)$ must be an isomorphism (onto its image) as well. Let us now discuss a situation when the converse holds.

We need the following
\blemma
\label{f.c.e.}
Let $X$ be a set. There exists a faithful conditional expectation
$$\Theta_X: \: \cL(\ell^2 X) \onto \ell^{\infty}(X)$$
such that, for every $T \in \cL(\ell^2 X)$, we have
\bgl
\label{char_THETA}
  \spkl{\Theta_X(T)\delta_x,\delta_y} = \delta_{x,y} \spkl{T \delta_x,\delta_y}
\egl
for all $x, y \in X$.
\elemma
\bproof
Let $e_{x,x}$ be the rank one projection onto $\Cz \delta_x \subseteq \ell^2 X$, given by
$$e_{x,x}(\xi) = \spkl{\xi,\delta_x}\delta_x \ {\rm for} \ {\rm all} \ \xi \in \ell^2 X.$$
Consider the linear map
\bgl
\label{map_THETA}
  \lspan(\menge{\delta_x}{x \in X}) \to \lspan(\menge{\delta_x}{x \in X}), \, \sum_x \alpha_x \delta_x \ma \sum_x \alpha_x (e_{x,x} \circ T)(\delta_x).
\egl
We have
\bglnoz
  \norm{\sum_x \alpha_x (e_{x,x} \circ T)(\delta_x)}^2
  &=& \spkl{\sum_x \alpha_x (e_{x,x} \circ T)(\delta_x),\sum_x \alpha_x (e_{x,x} \circ T)(\delta_x)}\\
  &=& \sum_x \abs{\alpha_x}^2 \spkl{(e_{x,x} \circ T)(\delta_x),(e_{x,x} \circ T)(\delta_x)}\\
  &\leq& \norm{T}^2 \sum_x \abs{\alpha_x}^2 = \norm{T}^2 \norm{\sum_x \alpha_x \delta_x}^2
\eglnoz
So the linear map in \eqref{map_THETA} extends to a bounded linear operator $\ell^2 X \to \ell^2 X$, which we denote by $\Theta_X(T)$. Our computation shows that
$$
  \norm{\Theta_X(T)} \leq \norm{T}.
$$
By definition,
$$\Theta_X(T)(\delta_x) = \spkl{T\delta_x,\delta_x} \delta_x.$$
This shows that $\Theta_X(T)$ lies in $\ell^{\infty}(X)$. It also shows that $\Theta_X(T)$ satisfies \eqref{char_THETA}.

Moreover, by construction, $\Theta_X(T) = T$ for all $T \in \ell^{\infty}(X)$. Therefore, the map
$$
  \Theta_X: \: \cL(\ell^2 X) \to \ell^{\infty}(X), \, T \ma \Theta_X(T)
$$
is a projection of norm $1$. Hence it follows by \cite[Theorem~II.6.10.2]{Bla06} that $\Theta_X$ is a conditional expectation.

Finally, $\Theta_X$ is faithful because given $T \in \cL(\ell^2 X)$, $\Theta_X(T^*T) = 0$ implies that
$$
  0 = \spkl{T^*T \delta_x,\delta_x} = \norm{T \delta_x}^2,
$$
so that $T \delta_x = 0$ for all $x \in X$, and hence $T = 0$.
\eproof

Applying Lemma~\ref{f.c.e.} to $X = I_l(P)\reg$ and $X = P$, we obtain faithful conditional expectations
$$
  \Theta_{I_l(P)}: \: \cL(\ell^2 I_l(P)\reg) \onto \ell^{\infty}(I_l(P)\reg)
$$
and
$$
  \Theta_P: \: \cL(\ell^2 P) \onto \ell^{\infty}(P).
$$
They fit into the following commutative diagram:
\bgl
\label{Theta-I*I}
  \xymatrix@C=20mm{
  \cL(\ell^2 I_l(P)\reg) \ar[r]^{\bI^* \, \sqcup \, \bI} \ar[d]_{\Theta_{I_l(P)}} & \cL(\ell^2 P) \ar[d]^{\Theta_P} \\
  \cL(\ell^2 I_l(P)\reg) \ar[r]^{\bI^* \, \sqcup \, \bI} & \ell^{\infty}(P)
  }
\egl
Here $\bI^* \sqcup \bI$ is our notation for the map sending $T$ to $\bI^* T \bI$. Commutativity of the diagram above follows from the following computation:
$$
  \Theta_P(\bI^* T \bI) \, \delta_x = \spkl{\bI^* T \bI \delta_x,\delta_x} \delta_x 
  = \spkl{T \delta_x,\delta_x} \delta_x = (\bI^* \Theta_{I_l(P)}(T) \bI) \, \delta_x.
$$

This leads us to
\bcor
\label{S-->P__E}
Assume that
\bgl
\label{Theta=cJ}
  \Theta_{I_l(P)}(C^*_{\lambda}(I_l(P))) = C^*(\cJ_P).
\egl
Then the *-homomorphism
$$
  C^*_{\lambda}(I_l(P)) \to C^*_{\lambda}(P), \, T \ma \bI^* T \bI
$$
from Lemma~\ref{S-->P} is an isomorphism if and only if its restriction to $C^*(\cJ_P)$,
$$
  C^*(\cJ_P) \onto D_{\lambda}(P), \, T \ma \bI^* T \bI,
$$
is an isomorphism.
\ecor
\bproof
Take the commutative diagram \eqref{Theta-I*I} and restrict the upper left corner to $$C^*_{\lambda}(I_l(P)) \subseteq \cL(\ell^2 I_l(P)\reg).$$
As $\bI^* C^*_{\lambda}(I_l(P)) \bI = C^*_{\lambda}(P)$ by Lemma~\ref{S-->P}, and because of \eqref{Theta=cJ}, we obtain the commutative diagram
\bgl
\label{Theta-I*I_v2}
  \xymatrix@C=36mm{
  C^*_{\lambda}(I_l(P)) \ar[r]^{\bI^* \, \sqcup \, \bI} \ar[d]_{\Theta_{I_l(P)}} & C^*_{\lambda}(P) \ar[d]^{\Theta_P} \\
  C^*(\cJ_P) \ar[r]^{\bI^* \, \sqcup \, \bI} & D_{\lambda}(P)
  }
\egl
As the vertical arrows are faithful, it is now easy to see that if the lower horizontal arrow is faithful, the upper horizontal arrow has to be faithful as well. This proves our corollary.
\eproof

\bremark
\label{C*E=cap}
The condition \eqref{Theta=cJ}, i.e.,
$$
  \Theta_{I_l(P)}(C^*_{\lambda}(I_l(P))) = C^*(\cJ_P),
$$
implies that
$$
  C^*(\cJ_P) = C^*_{\lambda}(I_l(P)) \cap \ell^{\infty}(I_l(P)\reg),
$$
and
$$
  D_{\lambda}(P) = C^*_{\lambda}(P) \cap \ell^{\infty}(P).
$$
This is because we always have
\bgl
\label{CCTheta_IP}
  C^*(\cJ_P) \subseteq C^*_{\lambda}(I_l(P)) \cap \ell^{\infty}(I_l(P)\reg) \subseteq \Theta_{I_l(P)}(C^*_{\lambda}(I_l(P))),
\egl
and
\bgl
\label{CCTheta_P}
  D_{\lambda}(P) \subseteq C^*_{\lambda}(P) \cap \ell^{\infty}(P) \subseteq \Theta_P(C^*_{\lambda}(P)),
\egl
and \eqref{Theta=cJ} implies that all these inclusions are equalities in \eqref{CCTheta_IP}, and also in \eqref{CCTheta_P} because
\bglnoz
  && D_{\lambda}(P)\\
  &=& \bI^* \, C^*(\cJ_P) \, \bI \overset{\eqref{Theta=cJ}}{=} \bI^* \, \Theta_{I_l(P)}(C^*_{\lambda}(I_l(P))) \, \bI = \Theta_P(\bI^* \, C^*_{\lambda}(I_l(P)) \, \bI)\\
  &=& \Theta_P(C^*_{\lambda}(P)).
\eglnoz
Here we used commutativity of the diagram in \eqref{Theta-I*I_v2} and Lemma~\ref{S-->P}.
\eremark

It remains to find out when condition~\eqref{Theta=cJ} holds. We follow \cite[\S~3.2]{Nor1}. Let us introduce the following
\bdefin
An inverse semigroup $S$ is called E*-unitary if for every $s \in S$, we must have $s \in E$ if there exists $x \in S\reg$ with $sx = x$. 
\edefin

\bremark
\label{stronglyE*->E*}
If there exists an idempotent pure partial homomorphism $\sigma: \: S\reg \to G$ to some group $G$, then $S$ is E*-unitary. This is because if we are given $s \in S$, and there exists $x \in S\reg$ with $sx = x$, then $\sigma(x) = \sigma(s) \sigma(x)$, so that $\sigma(s) = e$, where $e$ is the identity element in $G$. Since $\sigma$ is idempotent pure, $s$ must lie in $E$. 
\eremark

Now we apply Lemma~\ref{f.c.e.} to $X = S\reg$. Then we get a faithful conditional expectation
$$
  \Theta_{S\reg}: \cL(\ell^2 S\reg) \onto \ell^{\infty}(S\reg),
$$
and we may apply it to elements in $C^*_{\lambda}(S)$.

\blemma
In the situation above, our inverse semigroup $S$ is E*-unitary if and only if for every $s \in S$, we always have
$$\Theta_{S\reg}(\lambda_s) = 0$$
or
$$s \in E \ {\rm and} \ \Theta_{S\reg}(\lambda_s) = \lambda_s.$$
\elemma
\bproof
For \an{$\Rarr$}, assume that $\Theta_{S\reg}(\lambda_s) \neq 0$. This is equivalent to saying that there exists $x \in S\reg$ with $sx = x$. But since $S$ is E*-unitary, this implies $s \in E$. And since $\lambda_s$ lies in $\ell^{\infty}(S)$ for all $s \in E$, we must have $\Theta_{S\reg}(\lambda_s) = \lambda_s$.

Conversely, for \an{$\Larr$}, take $s \in S$ and suppose that there is $x \in S\reg$ with $sx = x$. Then $sxx^{-1} = xx^{-1}$, so that $sxx^{-1}$ is idempotent, and we conclude that
$$
  s^{-1}sxx^{-1} = (xx^{-1}s^{-1})(sxx^{-1}) = sxx^{-1} = xx^{-1},
$$
i.e., $s^{-1}s \geq xx^{-1}$. Hence
$$
  \lambda_s (\delta_x) = \delta_{sx} = \delta_x.
$$
Hence it follows that $\Theta_{S\reg}(\lambda_s) \neq 0$, and this implies, by assumption, that $s$ lies in $E$.
\eproof

In particular, we can draw the following conclusion
\bcor
\label{E*->Theta}
If $S$ is an E*-unitary inverse semigroup, then $\Theta_{S\reg}(C^*_{\lambda}(S)) = C^*(E)$.
\ecor

Combining Corollary~\ref{S-->P__E}, Remark~\ref{C*E=cap}, Corollary~\ref{E*->Theta}, Remark~\ref{stronglyE*->E*} and the observation that $I_l(P)$ admits an idempotent pure partial homomorphism to a group if $P$ embeds into a group (see \S~\ref{ss:ISGP}), we obtain
\bcor
\label{PG_SP}
Assume that $P$ is a semigroup which embeds into a group $G$. Then condition~\eqref{Theta=cJ} holds, i.e., 
$$
  \Theta_{I_l(P)}(C^*_{\lambda}(I_l(P))) = C^*(\cJ_P),
$$
and the *-homomorphism
$$
  C^*_{\lambda}(I_l(P)) \to C^*_{\lambda}(P), \, T \ma \bI^* T \bI
$$
from Lemma~\ref{S-->P} is an isomorphism if and only if its restriction to $C^*(\cJ_P)$,
$$
  C^*(\cJ_P) \onto D_{\lambda}(P), \, T \ma \bI^* T \bI,
$$
is an isomorphism.

Moreover,
$$
  C^*(\cJ_P) = C^*_{\lambda}(I_l(P)) \cap \ell^{\infty}(I_l(P)\reg),
$$
and
\bgl
\label{diagonal_P}
  D_{\lambda}(P) = C^*_{\lambda}(P) \cap \ell^{\infty}(P).
\egl
\ecor

Corollary~\ref{PG_SP} prompts the question when the *-homomorphism
$$
  C^*(\cJ_P) \onto D_{\lambda}(P), \, T \ma \bI^* T \bI,
$$
is an isomorphism. Note that both $C^*(\cJ_P)$ and $D_{\lambda}(P)$ are generated by a family of commuting projections, closed under multiplication, and our *-homomorphism sends generator to generator, i.e., $\lambda_X$ to $1_X$ for all $X \in \cJ_P$. Let us now investigate when such a *-homomorphism is an isomorphism.

\subsection{C*-algebras generated by semigroups of projections}

We basically follow \cite[\S~2.6]{LiSG} in this subsection.

If we think of elements of an inverse semigroup as partial isometries on a Hilbert space, then the semilattice of idempotents is a family of commuting projections, closed under multiplication, or in other words, a semigroup of projections.

Let us consider the general setting of a semilattice $E$ of idempotents, i.e., $E$ is an abelian semigroup consisting of idempotents. Suppose that $D$ is a C*-algebra generated by a multiplicatively closed family $\menge{d_e}{e \in E}$ of projections such that
$$
  E \to D, \, e \ma d_e
$$
is a semigroup homomorphism.

We make the following easy observation
\blemma
\label{Vd_e}
For every finite subset $F$ of $E$, there exists a projection in $D$, denoted by $\bigvee_{f \in F} d_f$, which is the smallest projection dominating all the projections $d_f$, $f \in F$.

Moreover, with $E(F)$ denoting the subsemigroup of $E$ generated by $F$, $\bigvee_{f \in F} d_f$ lies in
$$
  \lspan(\menge{d_e}{e \in E(F)}).
$$
\elemma
Just to be clear, the projection $\bigvee_{f \in F} d_f$ is uniquely characterized by 
$$
  d_f \leq \bigvee_{f \in F} d_f \ {\rm for} \ {\rm all} \ f \in F,
$$
and whenever a projection $d \in D$ satisfies
$$
  d_f \leq d \ {\rm for} \ {\rm all} \ f \in F,
$$
then we must have
$$
  \bigvee_{f \in F} d_f \leq d.
$$
\bproof
We proceed inductively on the cardinality of $F$. The case $\abs{F} = 1$ is trivial. Now assume that our claim holds for a finite subset $F$, and take an arbitrary element $\ti{f} \in E$. We want to check our claim for $F \cup \gekl{\ti{f}}$. Consider the element
\bgl
\label{Vdfvde}
  \bigvee_{f \in F} d_f + d_{\ti{f}} - \rukl{\bigvee_{f \in F} d_f} \cdot d_{\ti{f}}.
\egl
It is easy to see that this is projection in $D$, which dominates all the $d_f$, $f \in F$, as well as $d_{\ti{f}}$. Moreover, if $d$ is a projection in $D$ which dominates all the $d_f$, $f \in F$, and also $d_{\ti{f}}$, then $d$ obviously also dominates the projection in \eqref{Vdfvde}. Furthermore, since $\bigvee_{f \in F} d_f$ lies in
$$
  \lspan(\menge{d_e}{e \in E(F)})
$$
by induction hypothesis, the projection in \eqref{Vdfvde} lies in 
$$
  \lspan(\menge{d_e}{e \in E(F \cup \gekl{\ti{f}})}).
$$
\eproof

As above, let $E$ be a semilattice of idempotents. Suppose that $D$ is a C*-algebra generated by projections $\menge{d_e}{e \in E}$ such that $d_0 = 0$ if $0 \in E$ and $d_{ef} = d_e d_f$ for all $e,f \in E$. We prove the following result about *-homomorphisms out of $D$.
\bprop
\label{D-->B}
Let $B$ be a C*-algebra containing a semigroup of projections $\menge{b_e}{e \in E}$ such that $b_0 = 0$ if $0 \in E$ and $b_{ef} = b_e b_f$ for all $e,f \in E$.

There exists a *-homomorphism $D \to B$ sending $d_e$ to $b_e$ for all $e \in E$ if and only of for every $e \in E$ and every finite subset $F \subseteq E$ such that $f \lneq e$ for all $f \in F$, the equation 
$$
  d_e = \bigvee_{f \in F} d_f \ {\rm in} \ D
$$
implies that
$$
  b_e = \bigvee_{f \in F} b_f \ {\rm in} \ B.
$$

In that case, the kernel of the *-homomorphism
$$
  D \to B, \, d_e \to b_e
$$
is generated by
$$
  \menge{d_e - \bigvee_{f \in F} d_f \in D}{e \in E, \, F \subseteq \menge{f \in E}{f \lneq e} \ {\rm finite}, \,
  b_e = \bigvee_{f \in F} b_f \ {\rm in} \ B}.
$$
\eprop
\bproof
Let us start with the first part. Our condition is certainly a necessary condition for the existence of a *-homomorphism $D \to B, \, d_e \to b_e$. To prove that it is also sufficient, write $E$ as an increasing union of finite subsemigroups $E_i$, i.e.,
$$
  E = \bigcup_i E_i.
$$
Let $D_i \defeq C^*(\menge{d_e}{e \in E_i})$. Obviously,
$$
  D = \overline{\bigcup_i D_i}.
$$
For every $e \in E_i$, let $F_e \defeq \menge{f \in E_i}{f \lneq e}$. Then, by Lemma~\ref{Vd_e},
$$
  d_e - \bigvee_{f \in F_e} d_f
$$
is a projection in $D_i$. It is easy to see that
$$
  \menge{d_e - \bigvee_{f \in F_e} d_f}{e \in E_i}
$$
is a family of pairwise orthogonal projections which generates $D_i$. Moreover, it is also easy to see that
$$
  \menge{b_e - \bigvee_{f \in F_e} b_f}{e \in E_i}
$$
is a family of pairwise orthogonal projections in $B$. Hence it follows that there exists a *-homomorphism $D_i \to B$ sending
$$
  d_e - \bigvee_{f \in F_e} d_f 
$$
to
$$
  b_e - \bigvee_{f \in F_e} b_f
$$
for all $e \in E_i$ if and only if 
$$
  d_e - \bigvee_{f \in F_e} d_f = 0 \ {\rm in} \ D_i
$$
implies
$$
  b_e - \bigvee_{f \in F_e} b_f = 0 \ {\rm in} \ B,
$$
for all $e \in E_i$. But this is precisely the condition in the first part of our proposition. Moreover, it is easy to see that the *-homomorphism $D_i \to B$ we just constructed sends $d_e$ to $b_e$ for all $e \in E_i$. Hence these *-homomorphisms, taken together for all $i$, are compatible and give rise to the desired *-homomorphism from $D = \overline{\bigcup_i D_i}$ to $B$.

For the second part of the proposition, let $I$ be the ideal of $D$ generated by
$$
  \menge{d_e - \bigvee_{f \in F} d_f \in D}{F \subseteq E \ {\rm finite}, \,
  b_e = \bigvee_{f \in F} b_f \ {\rm in} \ B}.
$$
Obviously, $I$ is contained in the kernel of $D \to B, \, d_e \ma b_e$. It remains to show that the induced *-homomorphism $D/I \to B$ is injective. With the $D_i$s as above, set $I_i \defeq I \cap D_i$. Obviously, we have
$$
  I = \overline{\bigcup_i I_i} \ \ {\rm and} \ \ D/I = \overline{\bigcup_i D_i/I_i}.
$$
Hence it suffices to prove that the restriction $D_i/I_i \to B$ is injective, or in other words, that the *-homomorphism $D_i \to B$ we constructed above has kernel equal to $I_i$. But we have seen that
$$
  \menge{d_e - \bigvee_{f \in F_e} d_f}{e \in E_i}
$$
is a family of pairwise orthogonal projections which generates $D_i$. So the kernel is generated by those projections
$$
  d_e - \bigvee_{f \in F_e} d_f
$$
for which we have
$$
  b_e - \bigvee_{f \in F_e} b_f = 0 \ {\rm in} \ B.
$$
Therefore, the kernel is $I_i$, as required.
\eproof

As before, let $D$ be a C*-algebra generated by a semigroup $\menge{d_e}{e \in E}$ of projections such that $d_0 = 0$ if $0 \in E$ and $d_{ef} = d_e d_f$ for all $e,f \in E$. We set $E\reg \defeq E$ if $E$ is a semilattice without zero, and $E\reg \defeq E \setminus \gekl{0}$ if $0 \in E$.
\bprop
\label{When-D-universal}
The following are equivalent:
\begin{itemize}
\item[(i)] Our C*-algebra $D$ is universal for representations of $E$ by projections, i.e., we have an isomorphism
$$
  D \overset{\cong}{\lori} C^*(\menge{v_e}{e \in E} \ \vert \ v_e^* = v_e = v_e^2, \, v_0 = 0 \ {\rm if} \ 0 \in E, \, v_{ef} = v_e v_f)
$$
sending $d_e$ to $v_e$.
\item[(ii)] For every $e \in E$ and every finite subset $F \subseteq E$ with $f \lneq e$ for all $f \in F$, we have
$$\bigvee_{f \in F} d_f \lneq d_e.$$
\item[(iii)] The projections $\menge{d_e}{e \in E\reg}$ are linearly independent in $D$.
\end{itemize}
\eprop
\bproof
Obviously, (iii) implies (ii).

Moreover, (ii) implies (i) by Proposition~\ref{D-->B}, because if (ii) holds, we can never have
$$
  d_e = \bigvee_{f \in F} d_f \ {\rm in} \ D
$$
for any finite subset $F \subseteq E$ with $f \lneq e$ for all $f \in F$.

It remains to prove that (i) implies (iii). First of all, consider the left regular representation $\lambda$ on $\ell^2 E\reg$ as in \S~\ref{ss:ISGP}. It is given by $\lambda_e \delta_x = \delta_x$ if $e \geq x$ and $\lambda_e \delta_x = 0$ if $e \ngeq x$. By universal property of $D$, there is a *-homomorphism $D \to \cL(\ell^2 E\reg)$ sending $d_e$ to $\lambda_e$. But it is easy to see that $\lambda_e = \lambda_f$ if and only if $e = f$. Hence it follows that $d_e = d_f$ if and only if $e = f$.

Furthermore, again by universal property of $D$, there exists a *-homomorphism
$$
  D \to D \otimes D, \, d_e \ma d_e \otimes d_e.
$$
Let
$$\cD = \lspan(\menge{d_e}{e \in E}) \subseteq D.$$
Restricting the *-homomorphism $D \to D \otimes D$ from above to $\cD$, we obtain a homomorphism $\Delta: \: \cD \to \cD \odot \cD$ which is determined by $d_e \ma d_e \otimes d_e$ for every $e \in E$.

We now deduce from the existence of such a homomorphism $\Delta$ that $\menge{d_e}{e \in E\reg}$ is a $\Cz$-basis of $D$. As $\menge{d_e}{e \in E\reg}$ generates $\cD$ as a $\Cz$-vector space, we can always find a subset $\cS$ of $E\reg$ such that $\menge{d_e}{e \in \cS}$ is a $\Cz$-basis for $\cD$. It then follows that $\menge{d_e \otimes d_f}{e, f \in \cS}$ is a $\Cz$-basis of $\cD \odot \cD$.

Now take $e \in E\reg$. We can find finitely many $e_i \in \cS$ and $\alpha_i \in \Cz$ with $d_e = \sum_i \alpha_i d_{e_i}$. Applying $\Delta$ yields
\bgloz
  \sum_{i,j} \alpha_i \alpha_j d_{e_i} \otimes d_{e_j} = d_e \otimes d_e = \Delta(d_e) 
  = \sum_i \alpha_i \Delta(d_{e_i}) = \sum_i \alpha_i d_{e_i} \otimes d_{e_i}. 
\egloz
Hence it follows that among the $\alpha_i$s, there can only be one non-zero coefficient which must be $1$. The corresponding vector $d_{e_i}$ must then coincide with $d_e$. This implies $e = e_i \in \cS$, i.e. $\menge{d_e}{e \in E\reg}$ is a $\Cz$-basis of $\cD$. This proves (iii).
\eproof

Now let $S$ be an inverse semigroup with semilattice of idempotents $E$, and let $C^*_{\lambda}(S)$ be its reduced C*-algebra. Recall that we defined
$$
  C^*(E) \defeq \menge{\lambda_e}{e \in E}.
$$
\blemma
\label{C*E-universal}
The C*-algebra $C^*(E)$ is universal for representations of $E$ by projections.
\elemma
\bproof
By Proposition~\ref{When-D-universal}, all we have to show is that for every $e \in E$ and every finite subset $F \subseteq E$ with $f \lneq e$ for all $f \in F$, we have
$$\bigvee_{f \in F} \lambda_f \lneq \lambda_e.$$
But this follows from $\lambda_f (\delta_e) = 0$ for all $f \in E$ with $f \lneq e$, while $\lambda_e (\delta_e) = \delta_e$ for all $e \in E\reg$.
\eproof

It turns out that $C^*(E)$ can be identified with the corresponding sub-C*-algebra of the full C*-algebra of $S$.
\bcor
\label{C*E=C*E}
We have an isomorphism
$$
  C^*(E) \overset{\cong}{\lori} C^*(\menge{v_e}{e \in E}) \subseteq C^*(S)
$$
sending $\lambda_e$ to $v_e$ for all $e \in E$.
\ecor
\bproof
By Lemma~\ref{C*E-universal}, there is a *-homomorphism
$$
  C^*(E) \overset{\cong}{\lori} C^*(\menge{v_e}{e \in E}) \subseteq C^*(S)
$$
sending $\lambda_e$ to $v_e$ for all $e \in E$. It is an isomorphism because the inverse is given by restricting the left regular representation $C^*(S) \to C^*_{\lambda}(S)$ to $C^*(\menge{v_e}{e \in E}) \subseteq C^*(S)$.
\eproof
This justifies why we denote the sub-C*-algebra $C^*(\menge{\lambda_e}{e \in E})$ of $C^*_{\lambda}(S)$ by $C^*(E)$.

\bcor
We have a canonical identification $\widehat{E} \cong \Spec(C^*(E))$.
\ecor
\bproof
This is because by universal property of $C^*(E)$ (see Lemma~\ref{C*E-universal}), there is a one-to-one correspondence between non-zero *-homomorphisms $C^*(E) \to \Cz$ and non-zero semigroup homomorphisms $E \to \gekl{0,1}$ (sending $0$ to $0$ if $0 \in E$).
\eproof

Now suppose that we have an inverse semigroup $S$ with semilattice of idempotents $E$, and that we have a surjective *-homomorphism $C^*(E) \to D$ sending $\lambda_e \to d_e$. Then $D$ is a commutative C*-algebra, and we can describe its spectrum as follows:
\bcor
\label{SpecD}
Viewing $\Spec(D)$ as a closed subspace of $\widehat{E}$, $\Spec(D)$ is given by the subspace of all $\chi \in \widehat{E}$ with the property that whenever we have $e \in E$ with $\chi(e) = 1$ and a finite subset $F \subseteq E$ with $f \lneq e$ for every $f \in F$ satisfying $d_e = \bigvee_{f \in F} d_f$ in $D$, then we must have $\chi(f) = 1$ for some $f \in F$.
\ecor
\bproof
This is an immediate consequence of Proposition~\ref{D-->B}.
\eproof

Now let us suppose that we have a left cancellative semigroup $P$. We now apply Corollary~\ref{SpecD} and Proposition~\ref{When-D-universal} to the situation where $S = I_l(P)$, $E = \cJ_P$ and $D = D_{\lambda}(P) \subseteq C^*_{\lambda}(P)$. First, we make the following easy observation:
\blemma
Suppose that we are given finitely many $X_i \in \cJ_P$. Then we have
$$
  \bigvee_i 1_{X_i} = 1_{\bigcup_i X_i} \ {\rm in} \ D_{\lambda}(P) \subseteq \ell^{\infty}(P).
$$
\elemma

The following follows immediately from Corollary~\ref{SpecD}:
\bcor
\label{Omega_P-in-hatE}
The spectrum $\Omega_P = \Spec(D_{\lambda}(P))$ is given by the closed subspace of $\widehat{\cJ_P}$ consisting of all $\chi \in \widehat{\cJ_P}$ with the property that for all $X \in \cJ_P$ with $\chi(X) = 1$ and all $X_1, \dotsc, X_n \in \cJ_P$ with $X = \bigcup_{i=1}^n X_i$ in $P$, we must have $\chi(X_i) = 1$ for some $1 \leq i \leq n$.
\ecor

Proposition~\ref{When-D-universal} yields in our situation:
\bcor
\label{When-C*E=D}
The following are equivalent:
\begin{itemize}
\item We have an isomorphism
$$
  D_{\lambda}(P) \overset{\cong}{\lori} 
  C^* \rukl{\menge{v_X}{X \in \cJ_P} \ \vline 
  \begin{array}{c}
  v_X^* = v_X = v_X^2,
  \\
  v_0 = 0 \ {\rm if} \ 0 \in \cJ_P,
  \\
  v_{X \cap Y} = v_X v_Y)
  \end{array}
  }, \, 1_X \ma v_X.
$$
\item We have an isomorphism
$$
  C^*(E) \overset{\cong}{\lori} D_{\lambda}(P), \, \lambda_X \ma 1_X.
$$
\item For every $X \in \cJ_P$ and all $X_1, \dotsc, X_n \in \cJ_P$, 
$$
  X = \bigcup_{i=1}^n X_i
$$
implies that $X = X_i$ for some $1 \leq i \leq n$.
\item The projections $\menge{1_X}{X \in \cJ_P\reg}$ are linearly independent in $D_{\lambda}(P)$.
\end{itemize}
\ecor

\subsection{The independence condition}
\label{independence}

Corollary~\ref{When-C*E=D} justifies the following
\bdefin
\label{Def:ind}
We say that our left cancellative semigroup $P$ satisfies the independence condition (or simply independence) if for every $X \in \cJ_P$ and all $X_1, \dotsc, X_n \in \cJ_P$, 
$$
  X = \bigcup_{i=1}^n X_i
$$
implies that $X = X_i$ for some $1 \leq i \leq n$.
\edefin

Let us now discuss examples of left cancellative semigroups which satisfy independence, and also some examples which do not. We start with the following
\blemma
Suppose that $P$ is a left cancellative semigroup with identity $e$. If every non-empty constructible right ideal of $P$ is principal, i.e.,
$$\cJ_P\reg = \menge{pP}{p \in P},$$
then $P$ satisfies independence.
\elemma
\bproof
Suppose that
$$
  pP = \bigcup_{i=1}^n p_i P
$$
for some $p, p_1, \dotsc, p_n \in P$. Then, since $P$ has an identity, the element $p$ lies in $pP$, hence we must have $p \in p_i P$ for some $1 \leq i \leq n$. But then, since $p_i P$ is a right ideal, we conclude that $pP \subseteq p_iP$. Hence it follows that $pP = p_iP$, since we always have $pP \supseteq p_iP$.
\eproof

When are all non-empty constructible right ideals principal? Here is a necessary and sufficient condition:
\blemma
\label{pPCAPqP}
For a left cancellative semigroup $P$ (with or without identity), we have 
$$\cJ_P\reg = \menge{pP}{p \in P}$$
if and only if the following criterion holds:

For all $p, q \in P$ with $pP \cap qP \neq \emptyset$, there exists $r \in P$ with $pP \cap qP = rP$.
\elemma
\bproof
Our criterion is certainly necessary, since $\cJ_P$ is a semilattice, hence closed under intersections. To show that our condition is also sufficient, we first observe that $\cJ_P$ can be characterized as the smallest family of subsets of $P$ containing $P$ itself and closed under left multiplication, i.e.,
$$
  X \in \cJ_P, \, p \in P \ \Rarr \ p(X) \in cJ_P,
$$
as well as pre-images under left multiplication, i.e.,
$$
  X \in \cJ_P, \, q \in P \ \Rarr \ q^{-1}(X) \in cJ_P.
$$
Now $\menge{pP}{p \in P}$ is obviously closed under left multiplication. Hence it suffices to prove that principal right ideals are also closed under pre-images under left multiplication, up to $\emptyset$. Take $p, q \in P$. We always have
$$
  q^{-1}(pP) = q^{-1}(pP \cap qP).
$$
Therefore, if $pP \cap qP = \emptyset$, then $q^{-1}(pP) = \emptyset$. If $pP \cap qP \neq \emptyset$, then by our criterion, there exists $r \in P$ with $pP \cap qP = rP$. As $rP \subseteq qP$, we must have $r \in qP$, so that we can write $r = qx$ for some $x \in P$. Therefore, we conclude that
$$
  q^{-1}(pP) = q^{-1}(pP \cap qP) = q^{-1}(rP) = q^{-1}(qxP) = xP.
$$
\eproof

For instance, positive cones in totally ordered groups (as in \S~\ref{ex:totalorder}) always satisfy independence. This is because if $P$ is such a positive cone, then for $p, q \in P$, we have $pP \cap qP = pP$ if $p \geq q$ and $pP \cap qP = qP$ if $p \leq q$. Hence, all constructible right ideals are principal by Lemma~\ref{pPCAPqP}.

Moreover, right-angled Artin monoids (see \S~\ref{ex:presentations}) satisfy independence. Actually, all non-empty constructible right ideals are principal, because the criterion of Lemma~\ref{pPCAPqP} is true. This will come out of our general discussion of graph products in \S~\ref{sec:GraphProducts}.

To discuss more examples, let us explain a general method for verifying the criterion in Lemma~\ref{pPCAPqP}. This is based on \cite{Deh}.

Suppose that we are given a monoid $P$ defined by a presentation, i.e., generators $\Sigma$ and relations $R$, so that $P = \spkl{\Sigma \, \vert \, R}^+$. Assume that all the relations in $R$ are of the form $w_1 = w_2$, where $w_1$ and $w_2$ are formal words in $\Sigma$. Now we introduce formal symbols
$$
  \menge{\sigma^{-1}}{\sigma \in \Sigma} \eqdef \Sigma^{-1},
$$
and look at formal words in $\Sigma$ and $\Sigma^{-1}$. For two such words $w$ and $w'$, we write $w \curvearrowright_R w'$ if $w$ can be transformed into $w'$ be finitely many of the following two possible steps:
\begin{itemize}
\item Delete $\sigma^{-1} \sigma$.
\item Replace $\sigma_i^{-1} \sigma_j$ by $uv^{-1}$ if $\sigma_i u = \sigma_j v$ is a relation in $R$.
\end{itemize}
We then say that our presentation $(\Sigma,R)$ is complete for $\curvearrowright_R$ if for two formal words $u$ and $v$ in $\Sigma$, we have
$$
  u^{-1}v \curvearrowright_R \varepsilon \ \ {\rm (where} \ \varepsilon \ {\rm is \ the \ empty \ word)}
$$
if and only if $u$ and $v$ define the same element in our monoid $P = \spkl{\Sigma \, \vert \, R}^+$.

There are criteria on $(\Sigma,R)$ which ensure completeness for $\curvearrowright_R$ (see \cite{Deh}).

If completeness for $\curvearrowright_R$ is given, then we can read of properties of our monoid $P = \spkl{\Sigma \, \vert \, R}^+$ from the presentation $(\Sigma,R)$. We refer the reader to \cite{Deh} for a general and more complete discussion. For our purposes, the following observation is important: If $(\Sigma,R)$ is complete for $\curvearrowright_R$, then $P = \spkl{\Sigma \, \vert \, R}^+$ has the property that

\begin{center}
for all $p, q \in P$ with $pP \cap qP \neq \emptyset$, there exists $r \in P$ with $pP \cap qP = rP$
\setlength{\parindent}{0cm} \setlength{\parskip}{0cm}

if and only if

for all $\sigma_i, \sigma_j \in \Sigma$, there is at most one relation of the form $\sigma_i u = \sigma_j v$ in $R$.
\setlength{\parindent}{0cm} \setlength{\parskip}{0.5cm}
\end{center}

Coming back to examples, it turns out that the presentations for Artin monoids, discussed in \S~\ref{ex:presentations}, are complete for $\curvearrowright_R$. Also, the presentations for Baumslag-Solitar monoids $B_{k,l}^+$, for $k,l \geq 1$, are complete for $\curvearrowright_R$. Furthermore, the presentation for the Thompson monoid $F^+$ is complete for $\curvearrowright_R$.

Following our discussion above, it is now easy to see that for Artin monoids, the Baumslag-Solitar monoids $B_{k,l}^+$, for $k,l \geq 1$, and the Thompson monoid $F^+$, all non-empty constructible right ideals are principal. In particular, all these examples satisfy independence.

For semigroups coming from rings, we have the following result:
\blemma
\label{MnREG}
Let $R$ be a principal ideal domain. For both semigroups $M_n\reg(R)$ and $M_n(R) \rtimes M_n\reg(R)$, every non-empty constructible right ideal is principal. 
\elemma

For the proof, we need the following
\blemma
\label{int-principal}
For every $a$, $c$ in $M_n\reg(R)$, there exists $x \in M_n\reg(R)$ such that
\bgloz
  a M_n(R) \cap c M_n(R) = x M_n(R) \ {\rm and} \ a M_n\reg(R) \cap c M_n\reg(R) = x M_n\reg(R).
\egloz
\elemma
\bproof
For brevity, we write $M$ for $M_n(R)$ and $M\reg$ for $M_n\reg(R)$.

We will use the observation that for every $z \in M\reg$, there exist $u$ and $v$ in $GL_n(R)$ such that $uzv$ is a diagonal matrix (see for instance \cite{Kap}).

To prove our lemma, let us first of all define $x$. Let $\ti{c} \in M\reg$ satisfy $c \ti{c} = \ti{c} c = \det(c) \cdot 1_n$ ($1_n$ is the identity matrix). Choose $u$ and $v$ in $GL_n(R)$ with
$$\ti{c} a = u \cdot \text{diag}(\alpha_1, \dotsc, \alpha_n) \cdot v,$$
where $\text{diag}(\alpha_1, \dotsc, \alpha_n)$ is the diagonal matrix with $\alpha_1$, ..., $\alpha_n$ on the diagonal. For all $1 \leq i \leq n$, set $\beta_i \defeq \lcm(\alpha_i,\det(c))$ and $\gamma_i \defeq \det(c)^{-1} \beta_i$. Then our claim is that we can choose $x$ as $x = c \cdot u \cdot \text{diag}(\gamma_1, \dotsc, \gamma_n)$. In the following, we verify our claim:
\bglnoz
  aM \cap cM &=& \ti{c}^{-1} (\ti{c} a M \cap (\det(c) \cdot 1_n)M) \\
  &=& \ti{c}^{-1} ((u \cdot \text{diag}(\alpha_1, \dotsc, \alpha_n) \cdot v)M  \cap (\det(c) \cdot 1_n)M) \\
  &=& \ti{c}^{-1} u (\text{diag}(\alpha_1, \dotsc, \alpha_n)M  \cap (\det(c) \cdot 1_n)M) \\
  &=& \ti{c}^{-1} \cdot u \cdot \text{diag}(\beta_1, \dotsc, \beta_n) M \\
  &=& \ti{c}^{-1} (\det(c) \cdot 1_n) \cdot u \cdot \text{diag}(\gamma_1, \dotsc, \gamma_n) M \\
  &=& c \cdot u \cdot \text{diag}(\gamma_1, \dotsc, \gamma_n) M.
\eglnoz
Thus we have shown $aM \cap cM = xM$. Exactly the same computation shows that $aM\reg \cap cM\reg = x M\reg$.
\eproof

\bproof[Proof of Lemma~\ref{MnREG}]
For $M_n\reg(R)$, our claim is certainly a consequence of the Lemma~\ref{int-principal}. For $M_n(R) \rtimes M_n\reg(R)$, first note that given $(b,a)$ and $(d,c)$ in $M_n(R) \rtimes M_n\reg(R)$, we have
\bglnoz
  && (b,a)(M_n(R) \rtimes M_n\reg(R)) = (b + a M_n(R)) \times (a M_n\reg(R)), \\
  && (d,c)(M_n(R) \rtimes M_n\reg(R)) = (d + c M_n(R)) \times (c M_n\reg(R)).
\eglnoz
Moreover, the intersection
$$(b + a M_n(R)) \cap (d + c M_n(R))$$
is either empty or of the form
$$y + (a M_n(R) \cap c M_n(R))$$
for some $y \in M_n(R)$. Now Lemma~\ref{int-principal} provides an element $x \in M_n\reg(R)$ with
$$
a M_n(R) \cap c M_n(R) = x M_n(R) \ {\rm and} \ a M_n\reg(R) \cap c M_n\reg(R) = x M_n\reg(R).
$$
Thus either
$$(b,a)(M_n(R) \rtimes M_n\reg(R)) \cap (d,c)(M_n(R) \rtimes M_n\reg(R))$$
is empty or we obtain
\bgloz
  (b,a)(M_n(R) \rtimes M_n\reg(R)) \cap (d,c)(M_n(R) \rtimes M_n\reg(R)) = (y,x)(M_n(R) \rtimes M_n\reg(R)).
\egloz
\eproof

In general, however, given an integral domain $R$, the semigroups $R\reg$ and $R \rtimes R\reg$ do not have the property that all non-empty constructible right ideals are principal. For example, just take a number field with non-trivial class number, and let $R$ be its ring of algebraic integers. The property that all non-empty constructible right ideals are principal, for $R\reg$ or $R \rtimes R\reg$, translates to the property of the ring $R$ of being a principal ideal domain. But this is not the case if the class number is bigger than $1$. However, for all rings of algebraic integers, and more generally, for all Krull rings $R$, the semigroups $R\reg$ and $R \rtimes R\reg$ do satisfy independence.

Let $R$ be an integral domain. Recall that we introduced the set $\cI(R)$ of constructible ideals in \S~\ref{ss:Krull}. It is now easy to see that
$$
  \cJ_{R\reg} = \menge{I\reg}{I \in \cI(R)}
$$
and
$$
  \cJ_{R \rtimes R\reg} = \menge{(r+I) \times I\reg}{r \in R, a, I \in \cI(R)},
$$
where $I\reg = I \setminus \gekl{0}$.

Let us make the following observation about the relationship between the independence condition for multiplicative semigroups and $ax+b$-semigroups:
\blemma
\label{independence:R<->ax+b}
Let $R$ be an integral domain. Then $R\reg$ satisfies independence if and only if $R \rtimes R\reg$ satisfies independence.
\elemma
\bproof
If $\cJ_{R \rtimes R\reg}$ is not independent, then we have a non-trivial equation of the form
$$
  (r+I) \rtimes I\reg = \bigcup_{i=1}^n (r_i + I_i) \times I_i\reg \text{ with } (r_i + I_i) \times I_i\reg \subsetneq (r+I) \rtimes I\reg.
$$
It is clear that
$$(r_i + I_i) \times I_i\reg \subsetneq (r+I) \rtimes I\reg$$
implies that $I_i \subsetneq I$, for all $1 \leq i \leq n$. Projecting onto the second coordinate of $R \times R\reg$, we obtain
$$I\reg = \bigcup_{i=1}^n I_i\reg.$$
This means that $R\reg$ does not satisfy independence.

Conversely, assume that $R\reg$ does not satisfy independence, so that we have a non-trivial equation of the form
$$I\reg = \bigcup_{i=1}^n I_i\reg$$
with $I_i\reg \subsetneq I\reg$. Hence it follows that
$$I = \bigcup_{i=1}^n I_i,$$
and $I_i \subsetneq I$ for all $1 \leq i \leq n$. By \cite[Theorem~18]{Gott}, we may assume without loss of generality that
$$[I:I_i] < \infty \ {\rm for} \ {\rm all} \ 1 \leq i \leq n.$$
But then we have
$$
  I \times I\reg = \bigcup_{i=1}^n \bigcup_{r + I_i \in I / I_i} (r+I_i) \times I_i\reg.
$$
This shows that $\cJ_{R \rtimes R\reg}$ does not satisfy independence.
\eproof

\blemma
\label{Krull-independence}
For a Krull ring $R$, both semigroups $R\reg$ and $R \rtimes R\reg$ satisfy independence.
\elemma
\bproof
We use the same notations as in \S~\ref{ss:Krull}.

Let $Q$ be the quotient field of $R$, and let $I$, $I_1$, ..., $I_n$ be ideals in $\cI(R)$ with $I_i \subsetneq I$ for all $1 \leq i \leq n$. Then for every $1 \leq i \leq n$, there exists $\mfp_i \in \cP(R)$ with
$$v_{\mfp_i}(I_i) > v_{\mfp_i}(I).$$
By Proposition~\ref{approx}, there exists $x \in Q\reg$ with
$$v_{\mfp_i}(x) = v_{\mfp_i}(I) \ {\rm for} \ {\rm all} \ 1 \leq i \leq n$$
and
$$v_{\mfp}(x) \geq v_{\mfp}(I) \ {\rm for} \ {\rm all} \ \mfp \in \cP(R) \setminus \gekl{\mfp_1, \dotsc, \mfp_r}.$$
Thus $x$ lies in $I$, but does not lie in $I_i$ for any $1 \leq i \leq n$. Therefore,
$$\bigcup_{i=1}^n I_i \subsetneq I,$$
and thus
$$\bigcup_{i=1}^n I_i\reg \subsetneq I\reg.$$
This shows that $R\reg$ satisfies independence. By Lemma~\ref{independence:R<->ax+b}, $R \rtimes R\reg$ must satisfy independence as well.
\eproof

Let us present an example of a semigroup coming from a ring which does not satisfy independence. Consider the ring $R \defeq \Zz[i \sqrt{3}]$. Its quotient field is given by $Q = \Qz[i \sqrt{3}]$. $R$ is not integrally closed in $Q$. Let $\alpha \defeq \halb (1 + i \sqrt{3})$. $\alpha$ is a primitive sixth root of unity. It is clear that $\alpha \notin R$. But $2 \alpha = 1 + i \sqrt{3}$ lies in $R$.

The integral closure of $R$ is given by $\bar{R} \defeq \Zz[\alpha]$. We claim that
$$2 \bar{R} = 2^{-1} (2 \alpha R) = 2^{-1} (1 + i \sqrt{3}) R.$$
To prove \an{$\subseteq$}, observe that $\bar{R} = \Zz \cdot 1 + \Zz \cdot \alpha$. Now $$2 \cdot (2 \cdot 1) = 4 = (1 + i \sqrt{3}) \cdot (1 - i \sqrt{3}) \in (1 + i \sqrt{3})R,$$
and
$$2 \cdot (2 \alpha) = 2 \cdot (1 + i \sqrt{3}) \in (1 + i \sqrt{3})R.$$
For \an{$\supseteq$}, let $x = m + n \cdot i \sqrt{3}$ be in $R$ such that $2x \in 2 \alpha R$. As
$$2 \alpha R = (1 + i \sqrt{3}) R = \Zz \cdot (1 + i \sqrt{3}) + \Zz \cdot ((1 + i \sqrt{3}) i \sqrt{3}) = \Zz \cdot (1 + i \sqrt{3}) + \Zz \cdot (-3 + i \sqrt{3}),$$
there exist $k, l \in \Zz$ with
$$
  2x = 2m + 2n \cdot i \sqrt{3} = k (1 + i \sqrt{3}) + l (-3 + i \sqrt{3}) = (k-3l) + (k+l)(i \sqrt{3}),
$$
so that $2m = k-3l$ and $2n = k+l$. It follows that $2n = 2m + 4l$, and thus $n = m + 2l$ or $m = n - 2l$. We conclude that
$$
  x = -2l + n \cdot (1 + i \sqrt{3}) \in 2 \bar{R}.
$$
This shows that $2 \bar{R} = 2^{-1} (1 + i \sqrt{3}) R$. Hence it follows that $2\bar{R}$ is a constructible (ring-theoretic) ideal of $R$.

We have $\bar{R} = R \cup \alpha R \cup \alpha^2 R$ in $Q$. This is because
$$R = \Zz + \Zz (2 \alpha), \ \alpha R = \Zz \alpha + \Zz (2 \alpha^2) = \Zz \alpha + \Zz (2 \alpha - 2) \ {\rm and} \ \alpha^2 R = \Zz (\alpha - 1) + \Zz 2.$$
Now take $x = m + n \alpha \in \bar{R}$ with $m, n \in \Zz$. If $n$ is even, then $x$ is contained in $R$. If $n$ is odd and $m$ is even, then write $l = \tfrac{m}{2}$. We have $$x = (n+m) \cdot \alpha + (-l) \cdot (2 \alpha - 2) \in \alpha R.$$
Finally, if $n$ is odd and $m$ is odd, we write $k = \tfrac{m+n}{2}$. Then
$$x = n \cdot (\alpha - 1) + k \cdot 2 \in \alpha^2 R.$$
This shows $\bar{R} = R \cup \alpha R \cup \alpha^2 R$. Therefore,
$$2 \bar{R} = 2R \cup 2 \alpha R \cup 2 \alpha^2 R = 2 R \cup (1 + i \sqrt{3}) R \cup (-1 + i \sqrt{3}) R.$$
But $2R \subsetneq 2 \bar{R}$, $(1 + i \sqrt{3})R \subsetneq 2 \bar{R}$ and $(-1 + i \sqrt{3})R \subsetneq 2 \bar{R}$. This means that $R\reg$ does not satisfy independence. By Lemma~\ref{independence:R<->ax+b}, $R \rtimes R\reg$ does not satisfy independence, either.

Let us present another example of a left cancellative semigroup not satisfying independence. Consider $P = \Nz \setminus \gekl{1}$. Clearly, $P$ is a semigroup under addition. We have the following constructible right ideals
$$
  2+P = \gekl{2,4,5,6,\dotsc} \ {\rm and} \ 3+P = \gekl{3,5,6,7,\dotsc}.
$$
Hence
$$
  5+\Nz = \gekl{5,6,7,\dotsc} = (2+P) \cap (3+P)
$$
is also a constructible right ideal of $P$. Moreover, it is clear that
$$
  5+\Nz = (5+P) \cup (6+P).
$$
But since $5+P \subsetneq 5+\Nz$ and $6+P \subsetneq 5+\Nz$, it follows that $P$ does not satisfy independence.

A similar argument shows that for every numerical semigroup of the form $\Nz \setminus F$, where $F$ is a non-empty finite subset of $\Nz$ such that $\Nz \setminus F$ is still closed under addition, the independence condition does not hold. The reader may also compare \cite{CunToricVar} for more examples of a similar kind (which are two-dimensional versions), where the independence condition typically fails.

Now let us come back to the comparison of reduced C*-algebras for left cancellative semigroups and their left inverse hulls. Combining Corollary~\ref{PG_SP} and Proposition~\ref{When-C*E=D}, we get
\bprop
\label{PG_SP_ind}
Let $P$ be a subsemigroup of a group. The *-homomorphism
$$
  C^*_{\lambda}(I_l(P)) \to C^*_{\lambda}(P), \, \lambda_p \ma V_p
$$
is an isomorphism if and only if $P$ satisfies independence.
\eprop

\subsection{Construction of full semigroup C*-algebras}

Proposition~\ref{PG_SP_ind} explains when we can identify $C^*_{\lambda}(I_l(P))$ and $C^*_{\lambda}(P)$ in a canonical way, in case $P$ embeds into a group. Motivated by this result, we construct full semigroup C*-algebras.

\bdefin
Let $P$ be a left cancellative semigroup, and $I_l(P)$ its left inverse hull. We define the full semigroup C*-algebra of $P$ as the full inverse semigroup C*-algebra of $I_l(P)$, i.e.,
$$
  C^*(P) \defeq C^*(I_l(P)).
$$
\edefin
Recall that $C^*(I_l(P))$ is the C*-algebra universal for *-representations of the inverse semigroup $I_l(P)$ by partial isometries (see \S~\ref{ss:ISGP}).

As we saw in \S~\ref{ss:ISGP}, there is a canonical *-homomorphism
$$
  C^*(I_l(P)) \to C^*_{\lambda}(I_l(P)), \, v_p \ma \lambda_p.
$$
Composing with the *-homomorphism
$$
  C^*_{\lambda}(I_l(P)) \to C^*_{\lambda}(P), \, \lambda_p \ma V_p,
$$
we obtain a canonical *-homomorphism
$$
  C^*(P) \to C^*_{\lambda}(P), \, v_p \ma V_p.
$$
We call it the left regular representation of $C^*(P)$.

\bremark
\label{lrrISO->ind}
It is clear that if the left regular representation of $C^*(P)$ is an isomorphism, then $P$ must satisfy independence. This is because the restriction of $C^*(P) \to C^*_{\lambda}(P)$ to $C^*(\menge{v_X}{X \in \cJ_P})$ is the composition
$$
  C^*(\menge{v_X}{X \in \cJ_P}) \to C^*(E) \to D_{\lambda}(P),
$$
and we know that the first *-homomorphism is always an isomorphism (see Corollary~\ref{C*E=C*E}), while the second one is an isomorphism if and only if $P$ satisfies independence (see Corollary~\ref{When-C*E=D}).
\eremark

Given a concrete left cancellative semigroup $P$, it is usually possible to find a natural and simple presentation for $C^*(P)$ as a universal C*-algebra generated by isometries and projections, subject to relations. Let us discuss some examples.

For the example $P = \Nz$, the full semigroup C*-algebra $C^*(\Nz)$ is the universal unital C*-algebra generated by one isometry,
$$
  C^*(\Nz) \cong C^*(v \ \vert \ v^*v = 1).
$$

For $P = \Nz \times \Nz$, $C^*(\Nz \times \Nz)$ is the universal unital C*-algebra generated by two isometries which *-commute, i.e., 
$$
  C^*(\Nz \times \Nz) \cong C^*(v_a, v_b \ \vert \ v_a^*v_a = 1 = v_b^* v_b, \, v_a v_b = v_b v_a, \, v_a^* v_b = v_b v_a^*).
$$
Note that this C*-algebra is a quotient of
$$
  C^*(v_a, v_b \ \vert \ v_a^*v_a = 1 = v_b^* v_b, \, v_a v_b = v_b v_a).
$$
As we remarked in \S~\ref{ss:am_P}, the latter C*-algebra is not nuclear by \cite[Theorem~6.2]{Mur4}. However, as we will see in \S~\ref{ss:am_P-C*}, this quotient, and hence $C^*(\Nz \times \Nz)$, is nuclear.

For the non-abelian free monoid on two generators $P = \Nz * \Nz$, $C^*(\Nz * \Nz)$ is the universal unital C*-algebra generated by two isometries with orthogonal range projections, i.e.,
$$
  C^*(\Nz * \Nz) \cong C^*(v_a, v_b \ \vert \ v_a^*v_a = 1 = v_b^* v_b, \, v_a v_a^* v_b v_b^* = 0).
$$

More generally, for a right-angled Artin monoid $P$, a natural and simple presentation for $C^*(P)$ has been established in \cite{CrLa1} (see also \cite{ELR}).

Let us also mention that for a class of left cancellative semigroups, full semigroup C*-algebras can be identified in a canonical way with semigroup crossed products by endomorphisms. Let $P$ be a left cancellative semigroup with constructible right ideals $\cJ_P$. We then have a natural action $\alpha$ of $P$ by endomorphisms on
$$
  D(P) \defeq C^*(\menge{v_X}{X \in \cJ_P}) \subseteq C^*(P),
$$
where $p \in P$ acts by the endomorphism
$$
  \alpha_p: \: D(P) \to D(P), \, v_X \ma v_{pX}.
$$
If $P$ is right reversible, i.e., $Pp \cap Pq \neq \emptyset$ for all $p, q \in P$, or if every non-empty constructible right ideal of $P$ is principal, i.e., $\cJ_P\reg = \menge{pP}{p \in P}$, then we have a canonical isomorphism
$$
  C^*(P) \cong D(P) \rtimes_{\alpha} P.
$$
We refer to \cite[\S~3]{LiSG} for more details. Writing out the definition of the crossed 
product, we get the following presentation:
\bgloz
  C^*(P) \cong C^* \rukl{\menge{e_X}{X \in \cJ_P} \cup \menge{v_p}{p \in P} \vline 
  \begin{array}{c}
  e_X^* = e_X = e_X^2; \, v_p^*v_p = 1; \\ 
  e_{\emptyset} = 0 \ {\rm if} \ \emptyset \in \cJ_P , \, e_P = 1,
  \\ 
  e_{X \cap Y} = e_X \cdot e_Y;
  \\
  \, v_{pq} = v_p v_q;
  \\
  v_p e_X v_p^* = e_{pX}
  \end{array}
  }
\egloz

In particular, for an integral domain $R$, we obtain the following presentation for the full semigroup C*-algebra of $R \rtimes R\reg$:
\bglnoz
  && C^*(R \rtimes R\reg)\\ 
  &\cong& C^* \rukl{
  \begin{array}{c}
  \menge{e_I}{I \in \cI(R)} \\
  \cup \menge{u^b}{b \in R} \\
  \cup \menge{s_a}{a \in R\reg} 
  \end{array}  
  \vline 
  \begin{array}{c}
  e_I^* = e_I = e_I^2; \\
  u^b (u^b)^* = 1 = (u^b)^* u^b; \, v_a^* v_a = 1 \\ 
  e_R = 1, \, e_{I \cap J} = e_I \cdot e_J;
  \\
  s_{ac} = s_a s_c, \, u^{b+d} = u^b u^d, \, s_a u^b = u^{ab} s_a;
  \\
  s_a e_I s_a^* = e_{aI};
  \\
  u^b e_I = e_I u^b \ {\rm if} \ b \in I, \, e_I u^b e_I = 0 \ {\rm if} \ b \notin I
  \end{array}
  }
\eglnoz
We refer to \cite[\S~2]{CDL} as well as \cite[\S~2.4]{LiSG}.

In order to explain how this definition of full semigroup C*-algebras is related to previous constructions in the literature, we mention first of all that our definition generalizes Nica's construction in the quasi-lattice ordered case \cite{Nica}. Moreover, in the case of $ax+b$-semigroups over rings of algebraic integers (or more generally Dedekind domains), our definition includes the construction in \cite{CDL}. In the case of subsemigroups of groups, our definition coincides with the construction, denoted by $C^*_s(P)$, in \cite[Definition~3.2]{LiSG}. Last but not least, we point out that in comparison with another construction in \cite[Definition~2.2]{LiSG}, our definition is always a quotient of the construction in \cite[Definition~2.2]{LiSG}, and in certain cases (see \cite[\S~3.1]{LiSG} for details), our definition is actually isomorphic to the construction in \cite[Definition~2.2]{LiSG}.

\subsection{Crossed product and groupoid C*-algebra descriptions of reduced semigroup C*-algebras}
\label{cropro-GPD-description-C*redP}

We now specialize to the case where our semigroup $P$ embeds into a group $G$. To explain the connection between amenability and nuclearity, we would like to write the reduced C*-algebra $C^*_{\lambda}(P)$ of $P$ as a reduced crossed product attached to a partial dynamical system, and hence as a reduced groupoid C*-algebra. Let us start with the underlying partial dynamical system.

We already saw that $\Omega_P = \Spec(D_{\lambda}(P))$ may be identified with the subspace of $\widehat{\cJ_P}$ given by the characters $\chi$ with the property that for all $X, X_1, \dotsc, X_n$ in $\cJ_P$ with $X = \bigcup_{i=1}^n X_i$, $\chi(X) = 1$ implies that $\chi(X_i) = 1$ for some $1 \leq i \leq n$ (see Corollary~\ref{Omega_P-in-hatE}).

Moreover, we introduced the partial dynamical system $G \curvearrowright \widehat{\cJ_P}$ in \S~\ref{ss:ISGP}. It is given as follows: Every $g \in G$ acts on
$$
U_{g^{-1}} = \menge{\chi \in \widehat{\cJ_P}}{\chi(s^{-1}s) = 1 \ {\rm for} \ {\rm some} \ s \in I_l(P)\reg \ {\rm with} \ \sigma(s) = g},
$$
and for $\chi \in U_{g^{-1}}$, $g.\chi = \chi(s^{-1} \sqcup s)$ where $s \in I_l(P)\reg$ is an element satisfying $\chi(s^{-1}s) = 1$ and $\sigma(s) = g$.

We now claim:
\blemma
$\Omega_P$ is an $G$-invariant subspace of $\widehat{\cJ_P}$.
\elemma
\bproof
Take $g \in G$ and $\chi \in U_{g^{-1}} \cap \Omega_P$, and suppose that $s \in I_l(P)\reg$ satisfies $\chi(s^{-1}s) = 1$ and $\sigma(s) = g$. We have to show that $g.\chi = \chi(s^{-1} \sqcup s)$ lies in $\Omega_P$. 

Suppose that $X, X_1, \dotsc, X_n$ in $\cJ_P$ satisfy $X = \bigcup_{i=1}^n X_i$. Then, identifying $s^{-1}s$ with $\dom(s)$, we have
$$
  s^{-1}Xs = (g^{-1}X) \cap \dom(s) = \bigcup_{i=1}^n (g^{-1}X_i) \cap \dom(s) = \bigcup_{i=1}^n s^{-1}X_is.
$$
Hence, if $g.\chi(X) = 1$, then $\chi(s^{-1}Xs) = 1$, and hence $g.\chi(X_i) = \chi(s^{-1}X_is) = 1$ for some $1 \leq i \leq n$. This shows that $g.\chi$ lies in $\Omega_P$.
\eproof

Hence we obtain a partial dynamical system $G \curvearrowright \Omega_P$ by restricting $G \curvearrowright \widehat{\cJ_P}$ to $\Omega_P$. A moment's thought shows that this partial dynamical system coincides with the one introduced in \S~\ref{ss:PDS}.

If our group $G$ were exact, then this observation, together with Corollary~\ref{C*S=C*ExG}, would immediately imply that $C^*_{\lambda}(P) \cong C(\Omega_P) \rtimes_r G$ with respect to the $G$-action $G \curvearrowright \Omega_P$. However, it turns out that we do not need exactness here.

\btheo
\label{P-DG}
There is a canonical isomorphism $C^*_{\lambda}(P) \cong C(\Omega_P) \rtimes_r G$ determined by $V_p \ma W_p$. Here $W_g$ denote the canonical partial isometries in $C(\Omega_P) \rtimes_r G$.
\etheo
\bproof
We work with the dual action $G \curvearrowright D_{\lambda}(P)$ as described in \S~\ref{ss:PDS}. Our strategy is to describe both $C^*_{\lambda}(P)$ and $D_{\lambda}(P) \rtimes_r G$ as reduced (cross sectional) algebras of Fell bundles, and then to identify the underlying Fell bundles.

Let us start with $C^*_{\lambda}(P)$. As in \S~\ref{ss:ISGP}, we think of $I_l(P)$ as partial isometries. Recall that we defined the partial homomorphism $\sigma: \: I_l(P)\reg \to G$ in \S~\ref{ss:ISGP}. Now we set
$$B_g \defeq \clspan(\sigma^{-1}(g))$$
for every $g \in G$. We want to see that $(B_g)_{g \in G}$ is a grading for $C^*_{\lambda}(P)$, in the sense of \cite[Definition~3.1]{Ex1}. Conditions (i) and (ii) are obviously satisfied. For (iii), we use the faithful conditional expectation $\Theta_P: \: C^*_{\lambda}(P) \onto D_{\lambda}(P) = B_e$ from \S~\ref{Comparing-I-P}. Given a finite sum
$$x = \sum_g x_g \in C^*_{\lambda}(P)$$
of elements $x_g \in B_g$ such that $x = 0$, we conclude that
$$0 = x^* x = \sum_{g,h} x_g^* x_h,$$
and hence
$$0 = \Theta_P(x^* x) = \sum_g x_g^* x_g.$$
Here we used that $\Theta_P \vert_{B_g} = 0$ if $g \neq e$. This implies that $x_g = 0$ for all $g$. Therefore, the subspaces $B_g$ are independent. It is clear that the linear span of all the $B_g$ is dense in $C^*_{\lambda}(P)$. This proves (iii). If we let $\cB$ be the Fell bundle given by $(B_g)_{g \in G}$, then \cite[Proposition~3.7]{Ex1} implies $C^*_{\lambda}(P) \cong C^*_r(\cB)$ because $\Theta_P: \: C^*_{\lambda}(P) \onto D_{\lambda}(P) = B_e$ is a faithful conditional expectation satisfying $\Theta_P \vert_{B_e} = \id_{B_e}$ and $\Theta_P \vert_{B_g} = 0$ if $g \neq e$.

Let us also describe $D_{\lambda}(P) \rtimes_r G$ as a reduced algebra of a Fell bundle. We denote by $W_g$ the partial isometry in $D_{\lambda}(P) \rtimes_r G$ corresponding to $g \in G$, and we set $B'_g \defeq D_g W_g$. Recall that we defined
$$D_{g^{-1}} = \clspan(\menge{V^*V}{V \in I_l(P)\reg, \, \sigma(V) = g})$$
in \S~\ref{ss:PDS}. It is easy to check that $(B'_g)_{g \in G}$ satisfy (i), (ii) and (iii) in \cite[Definition~3.1]{Ex1}. Moreover, $B'_e = D_e = D_{\lambda}(P)$, and it follows immediately from the construction of the reduced partial crossed product that there is a faithful conditional expectation $D_{\lambda}(P) \rtimes_r G \onto D_{\lambda}(P) = B'_e$ which is identity on $B'_e$ and $0$ on $B'_g$ for $g \neq e$. Hence if we let $\cB'$ be the Fell bundle given by $(B'_g)_{g \in G}$, then \cite[Proposition~3.7]{Ex1} implies $D_{\lambda}(P) \rtimes_r G \cong C^*_r(\cB')$.

To identify $C^*_{\lambda}(P)$ and $D_{\lambda}(P) \rtimes_r G$, it now remains to identify $\cB$ with $\cB'$. We claim that the map
$$\lspan(\menge{V}{\sigma(V) = g}) \to \lspan(\menge{VV^* W_g}{\sigma(V) = g}), \, \sum_i \alpha_i V_i \ma \sum_i \alpha_i V_i V_i^* W_g$$
is well-defined and extends to an isometric isomorphism $B_g \to B'_g$, for all $g \in G$.

All we have to show is that our map is isometric. We have
$$
\norm{\sum_i \alpha_i V_i}^2 = \norm{\sum_{i,j} \alpha_i \overline{\alpha_j} V_i V_j^*}_{D_{\lambda}(P)}$$
and
$$
\norm{\sum_i \alpha_i V_i V_i^* W_g}^2 = \norm{\sum_{i,j} \alpha_i \overline{\alpha_j} V_i V_i^* V_j V_j^*}_{D_{\lambda}(P)}.
$$
Since $V_i = V_i V_i^* \lambda_g$ and $V_j^* = \lambda_{g^{-1}} V_j V_j^*$, we have 
$$V_i V_j^* = V_i V_i^* \lambda_g \lambda_{g^{-1}} V_j V_j^* = V_i V_i^* V_j V_j^*.$$
Hence, indeed,
$$\norm{\sum_i \alpha_i V_i}^2 = \norm{\sum_i \alpha_i V_i V_i^* W_g}^2,$$
and we are done.

All in all, we have proven that
$$C^*_{\lambda}(P) \cong C^*_r(\cB) \cong C^*_r(\cB') \cong D_{\lambda}(P) \rtimes_r G.$$
Our isomorphism sends $V_p$ to $V_p V_p^* W_p$, but a straightforward computation shows that actually, $V_p V_p^* W_p = W_p$ for all $p \in P$. Thus the isomorphism we constructed is given by $V_p \ma W_p$ for all $p \in P$.
\eproof

In particular, in combination with Theorem~\ref{THM:C*GPD=C*PDS_red}, we get an isomorphism
\bgl
\label{P-->GxOmega}
  C^*_{\lambda}(P) \overset{\cong}{\lori} C^*_r(G \ltimes \Omega_P), \, V_p \ma 1_{\gekl{p} \times \Omega_P}.
\egl
Together with Remark~\ref{isom-leftreg} and Lemma~\ref{GS=GxhatE}, we see that we obtain a commutative diagram
\bgl
\label{PPPGPDGPDGPD}
  \xymatrix@C=20mm{
  C^*(P) = C^*(I_l(P)) \ar[r]^{\cong} \ar[d] & C^*(G \ltimes \widehat{\cJ_P}) \ar[d] \\
  C^*_{\lambda}(I_l(P)) \ar[r]^{\cong} \ar[d] & C^*_r(G \ltimes \widehat{\cJ_P}) \ar[d] \\
  C^*_{\lambda}(P) \ar[r]^{\cong} & C^*_r(G \ltimes \Omega_P)
  }
\egl
Here the upper left vertical arrow is the left regular representation of $C^*(I_l(P))$. The lower left vertical arrow is the *-homomorphism provided by Lemma~\ref{S-->P}. The upper right vertical arrow is the left regular representation of $C^*(G \ltimes \widehat{\cJ_P})$. The lower right vertical arrow is the canonical projection map; it corresponds to the canonical map $C(\widehat{\cJ_P}) \rtimes_r G \onto C(\Omega) \rtimes_r G$ under the identification from Theorem~\ref{THM:C*GPD=C*PDS_red}. The first horizontal arrow is the identifications from Theorem~\ref{THM:C*GPD=C*PDS}. The second horizontal arrow is the isomorphism from Theorem~\ref{THM:C*GPD=C*PDS_red}. For both of these horizontal arrows, we also need Lemma~\ref{GS=GxhatE}. The third horizontal arrow is provided by the isomorphism \eqref{P-->GxOmega}.

Now we are ready to discuss the relationship between amenability and nuclearity and thereby explain the strange phenomena mentioned at the beginning of \S~\ref{ss:am_P}.

\subsection{Amenability of semigroups in terms of C*-algebras}
\label{ss:am_P-C*}

Let us start by explaining how to characterize amenability of semigroups in terms of their C*-algebras.
\btheo
\label{THM:am_P-C*}
Let $P$ be a cancellative semigroup, i.e., $P$ is both left and right cancellative. Assume that $P$ satisfies the independence condition. Then the following are equivalent:
\begin{itemize}
\item[1)] $P$ is left amenable.
\item[2)] $C^*(P)$ is nuclear and there is a character on $C^*(P)$.
\item[3)] $C^*_{\lambda}(P)$ is nuclear and there is a character on $C^*(P)$. 
\item[4)] The left regular representation $C^*(P) \to C^*_{\lambda}(P)$ is an isomorphism and there is a character on $C^*(P)$. 
\item[5)] There is a character on $C^*_{\lambda}(P)$.
\end{itemize}
\etheo
By a character, we mean a unital *-homomorphism to $\Cz$.

For the proof, we need the following
\blemma
\label{char-leftrev}
Let $P$ be a left cancellative semigroup. The following are equivalent:
\begin{itemize}
\item[1.] There is a character on $C^*(P)$.
\item[2.] $P$ is left reversible, i.e., $pP \cap qP \neq \emptyset$ for all $p, q \in P$.
\item[3.] $I_l(P)$ does not contain $\emptyset \to \emptyset$, the partial bijection which is nowhere defined.
\end{itemize}
\elemma
Recall that in the convention we introduced in \S~\ref{ss:ISGP}, if $\emptyset \to \emptyset$ lies in $I_l(P)$, then we say that $I_l(P)$ is an inverse semigroup with zero, and let $\emptyset \to \emptyset$ be its distinguished zero element, which we denote by $0$.
\bproof
1. $\Rarr$ 2.: If $\chi$ is a character on $C^*(P)$, then for every $p, q \in P$, we have
$$
  \chi(1_{pP \cap qP}) = \chi(1_{pP}) \chi(1_{qP}) = \chi(V_p V_p^*) \chi(V_q V_q^*) = \abs{\chi(V_p)}^2 \abs{\chi(V_q)}^2 = 1.
$$
Hence $pP \cap qP \neq \emptyset$.

2. $\Rarr$ 3.: Every partial bijection in $I_l(P)$ is a finite product of elements in
$$
\menge{p}{p \in P} \cup \menge{q^{-1}}{q \in P}.
$$
Hence, by an inductive argument, it suffices to show that if $s \in I_l(P)$ is not $\emptyset \to \emptyset$, then for all $p, q \in P$, $ps$ and $q^{-1}s$ are not $\emptyset \to \emptyset$. For $ps$, this is clear. For $q^{-1}s$, choose $x \in \dom(s)$. Then $xP \subseteq \dom(s)$ and $s(xr) = s(x)r$ for all $r \in P$ by property~\eqref{s(xr)=s(x)r}. As $P$ is left reversible, there exists $y \in P$ with $y \in qP \cap s(x)P$. Hence $y = s(x)r = qz$ for some $r, z \in P$. Therefore,
$$
  (q^{-1}s)(xz) = q^{-1}(s(xr)) = q^{-1}(s(x)r) = q^{-1}(qz) = z.
$$
Hence $q^{-1}s$ is not $\emptyset \to \emptyset$, as desired.

3. $\Rarr$ 1.: Since $I_l(P)$ does not contain $\emptyset \to \emptyset$, we have by definition that
$$
  C^*(P) = C^*(I_l(P)) = C^*(\menge{v_s}{s \in I_l(P)} \vert v_{st} = v_s v_t, \, v_{s^{-1}} = v_s^*).
$$
Obviously, by universal property, we obtain a character $C^*(P) \to \Cz, \, v_s \to 1$.
\eproof

\bproof[Proof of Theorem~\ref{THM:am_P-C*}]
1) $\Rarr$ 2): If $P$ is left amenable, then there exists a left invariant state $\mu$ on $\ell^{\infty}(P)$ by definition. Hence, for every $p \in P$, we have
$$
  \mu(1_{pP}) = \mu(1_{pP}(p \sqcup)) = \mu(1_P) = 1.
$$
Now, if there were $p, q \in P$ with $pP \cap qP = \emptyset$, then $1_{pP} + 1_{qP}$ would be a projection in $\ell^{\infty}(P)$ with $1_{pP} + 1_{qP} \leq 1_P$, so that
$$
  1 = \mu(1_P) \geq \mu(1_{pP} + 1_{qP}) = \mu(1_{pP}) + \mu(1_{qP}) = 1+1 = 2.
$$
This is a contradiction. Therefore, $P$ must be left reversible. By Lemma~\ref{char-leftrev}, it follows that $C^*(P)$ has a character.

In addition, by our discussion of group embeddability in \S~\ref{PinG}, we see that $P$ embeds into its group $G$ of right quotients. Moreover, as $P$ is left amenable, $G$ must be amenable by \cite[Proposition~(1.27)]{Pat}. Hence, statement 2) follows from Theorem~\ref{THM:nuc-am} (see also Corollary~\ref{Cor:PinGam}).

2) $\Rarr$ 3) is obvious.

3) $\Rarr$ 4) follows again from Theorem~\ref{THM:nuc-am}.

4) $\Rarr$ 5) is obvious.

5) $\Rarr$ 1): We follow \cite[\S~4.2]{LiSG}. Let $\chi: \: C^*_{\lambda}(P) \to \Cz$ be a non-zero character. Viewing $\chi$ as a state, we can extend it by the theorem of Hahn-Banach to a state on $\cL(\ell^2(P))$. We then restrict the extension to $\ell^{\infty}(P) \subseteq \cL(\ell^2(P))$ and call this restriction $\mu$. The point is that by construction, $\mu \vert_{C^*_{\lambda}(P)} = \chi$ is multiplicative, hence $C^*_{\lambda}(P)$ is in the multiplicative domain of $\mu$. Thus we obtain for every $f \in \ell^{\infty}(P)$ and $p \in P$
\bgloz
  \mu(f(p \sqcup)) = \mu(V_p^* f V_p) = \mu(V_p^*) \mu(f) \mu(V_p) = \mu(V_p)^* \mu(V_p) \mu(f) = \mu(f).
\egloz
Thus $\mu$ is a left invariant mean on $\ell^{\infty}(P)$. This shows \an{5) $\Rarr$ 1)}.
\eproof

Theorem~\ref{THM:am_P-C*} tells us that for the example $P = \Nz \times \Nz$ discussed in \S~\ref{ss:am_P}, our definition of full semigroup C*-algebras leads to a full C*-algebra $C^*(\Nz \times \Nz)$ which is nuclear and whose left regular representation is an isomorphism. This explains and resolves the strange phenomenon described in \S~\ref{ss:am_P}.

At the same time, we see why it is not a contradiction that $\Nz * \Nz$ is not amenable while its C*-algebra behaves like those of amenable semigroups. The point is that there is no character on $C^*(\Nz * \Nz)$ because $\Nz * \Nz$ is not left reversible.

However, we still need an explanation why the semigroup C*-algebra of $\Nz * \Nz$ behaves like those of amenable semigroups. This leads us to our next result.

\subsection{Nuclearity of semigroup C*-algebras and the connection to amenability}

\btheo
\label{THM:nuc-am}
Let $P$ be a semigroup which embeds into a group $G$. Consider
\begin{itemize}
\item[(i)] $C^*(P)$ is nuclear.
\item[(ii)] $C^*_{\lambda}(P)$ is nuclear.
\item[(iii)] $G \ltimes \Omega_P$ is amenable.
\item[(iv)] The left regular representation $C^*(P) \to C^*_{\lambda}(P)$ is an isomorphism.
\end{itemize}
We always have (i) $\Rarr$ (ii) $\LRarr$ (iii), and (iv) implies that $P$ satisfies independence.

If $P$ satisfies independence, then we also have (iii) $\Rarr$ (i) and (iii) $\Rarr$ (iv).
\etheo
Note that the \'{e}tale locally compact groupoid $G \ltimes \Omega_P$ really only depends on $P$, not on the embedding $P \into G$. This follows from Lemma~\ref{GS=GxhatE}.

\bproof
The first claim follows from the description of $C^*(P) = C^*(I_l(P))$ as a full groupoid C*-algebra (see Theorem~\ref{THM:C*S=C*GPD}), the description of $C^*_{\lambda}(P)$ as a reduced groupoid C*-algebra (see \S~\ref{cropro-GPD-description-C*redP} and the isomorphism~\eqref{P-->GxOmega}), the commutative diagram~\eqref{PPPGPDGPDGPD}, and Theorem~\ref{gpd-am-nuc}. That (iv) implies that $P$ satisfies independence was explained in Remark~\ref{lrrISO->ind}.

The second claim follows from the observation that if $P$ satisfies independence, then $\Omega_P = \widehat{\cJ_P}$ (see Corollary~\ref{Omega_P-in-hatE} and equation~\eqref{diagonal_P}), so that the partial dynamical systems $G \curvearrowright \Omega_P$ and $G \curvearrowright \widehat{\cJ_P}$, and hence their partial transformation groupoids coincide, and Theorem~\ref{gpd-am-nuc}.
\eproof

\bcor
\label{Cor:PinGam}
If $P$ is a subsemigroup of an amenable group $G$, then statements (i), (ii) and (iii) from Theorem~\ref{THM:nuc-am} hold, and (iv) holds if and only if $P$ satisfies independence.
\ecor
\bproof
This is because if $G$ is amenable, the partial transformation groupoid $G \ltimes \Omega_P$ is amenable by \cite[Theorem~20.7 and Theorem~25.10]{Ex4}.
\eproof

This explains the second strange phenomenon mentioned at the beginning of \S~\ref{ss:am_P}, that the semigroup C*-algebra of $\Nz * \Nz$ behaves like those of amenable semigroups. The underlying reason is that $\Nz * \Nz$ embeds into an amenable group: Let $\Fz_2$ be the free group on two generators. By \cite{Hoch}, we have an embedding $\Nz * \Nz \into \Fz_2 / \Fz_2''$, where $\Fz_2''$ is the second commutator subgroup of $\Fz_2$. But $\Fz_2 / \Fz_2''$ is solvable, in particular amenable. Moreover, $\Nz * \Nz$ satisfies independence (see \S~\ref{independence}). This is why statements (i) to (iv) from Theorem~\ref{THM:nuc-am} are all true for the semigroup $P = \Nz * \Nz$.

\bremark
If we modify the definition of full semigroup C*-algebras, then we can get the same results as in Theorem~\ref{THM:am_P-C*}, Theorem~\ref{THM:nuc-am} and Corollary~\ref{Cor:PinGam} without having to mention the independence condition. Simply define $C^*(P)$ as the full groupoid C*-algebra of the restriction
$$
  \cG(I_l(P)) \vert \Omega_P = \menge{\gamma \in \cG(I_l(P))}{r(\gamma), s(\gamma) \in \Omega_P}
$$
of the universal groupoid $\cG(I_l(P))$ of $I_l(P)$ to $\Omega_P$. This means that we would set
$$
  C^*(P) \defeq C^*(\cG(I_l(P)) \vert \Omega_P).
$$
Then, in Theorem~\ref{THM:nuc-am}, we would have (i) $\LRarr$ (ii) $\LRarr$ (iii), and all these statements imply (iv). Corollary~\ref{Cor:PinGam} would say that statements (i) to (iv) from Theorem~\ref{THM:nuc-am} hold whenever $G$ is amenable. Moreover, Theorem~\ref{THM:am_P-C*} would be true without the assumption that $P$ satisfies independence.

We have chosen not to follow this route and keep the definition of full semigroup C*-algebras as full C*-algebras of left inverse hulls because the C*-algebras $C^*(I_l(P))$ usually have a nicer presentation, i.e., a nicer and simpler description as universal C*-algebras given by generators and relations. Moreover, in the case of semigroups embeddable into groups, we know that these two definitions of full semigroup C*-algebras differ precisely by the (failure of the) independence condition.
\eremark

\section{Topological freeness, boundary quotients, and C*-simplicity}
\label{sec:TFBQsimp}

Given a semigroup $P$ which embeds into a group $G$, we have constructed a partial dynamical system $G \curvearrowright \Omega_P$ and identified the reduced semigroup C*-algebra $C^*_{\lambda}(P)$ with the reduced crossed product $C(\Omega_P) \rtimes_r G$. Let us now present a criterion for topological freeness of $G \curvearrowright \Omega_P$. First recall (compare \cite{ELQ} and \cite{Li-PTGPD}) that a partial dynamical system $G \curvearrowright X$ is called topologically free if for every $e \neq g \in G$,
$$
  \menge{x \in U_{g^{-1}}}{g.x \neq x}
$$
is dense in $U_{g^{-1}}$. Here, we use the same notation as in \S~\ref{ss:PDS}.

We first need the following observation: Let $P$ be a monoid. For $p \in P$, let $\chi_{pP} \in \widehat{\cJ_P}$ be defined by $\chi_{pP}(X) = 1$ if and only if $pP \subseteq X$, for $X \in \cJ_P$. Since $P$ is a monoid, $\chi_{pP}$ lies in $\Omega_P$ for all $p \in P$.
\blemma
\label{chi_P:dense}
The subset $\menge{\chi_{pP}}{p \in P}$ is dense in $\Omega_P$.
\elemma
\bproof
Basic open sets in $\Omega_P$ are of the form
$$U(X;X_1, \dotsc, X_n) = \menge{\chi \in \Omega_P}{\chi(X) = 1, \, \chi(X_i) = 0 \ {\rm for} \ {\rm all} \ 1 \leq i \leq n}.$$
Here $X, X_1, \dotsc, X_n$ are constructible ideals of $P$. Clearly, $U(X;X_1, \dotsc, X_n)$ is empty if $X = \bigcup_{i=1}^n X_i$. Thus, for a non-empty basic open set $U(X;X_1, \dotsc, X_n)$, we may choose $p \in X$ such that $p \notin \bigcup_{i=1}^n X_i$, and then $\chi_{pP} \in U(X;X_1, \dotsc, X_n)$.
\eproof

\btheo
\label{P*=e->TF}
Let $P$ be a monoid with identity $e$ which embeds into a group $G$. If $P$ has trivial units $P^* = \gekl{e}$, then $G \curvearrowright \Omega_P$ is topologically free.
\etheo
\bproof
For $p \in P$, let $\chi_{pP} \in \widehat{\cJ_P}$ be defined as in Lemma~\ref{chi_P:dense}, i.e., $\chi_{pP}(X) = 1$ if and only if $pP \subseteq X$, for $X \in \cJ_P$. Assume that $g \in G$ satisfies $g.\chi_{pP} = \chi_{pP}$ for some $p \in P$. This equality only makes sense if $\chi_p \in U_{g^{-1}}$, i.e., there exists $s \in I_l(P)$ with $\sigma(s) = g$ and $\chi_p(s^{-1}s) = 1$. The latter condition is equivalent to $pP \subseteq \dom(s)$. Then
$$
  g.\chi_{pP}(X) = \chi_{pP}(s^{-1}Xs) = \chi_{pP}(s^{-1}(X \cap \img(s))) = \chi_{pP}(g^{-1}(X \cap \img(s))).
$$
So for $X \in \cJ_P$,
$$g.\chi_{pP}(X) = 1$$
if and only if
$$pP \subseteq g^{-1}(X \cap \img(s)) = g^{-1}X \cap \dom(s).$$
But since $pP \subseteq \dom(s)$ holds, we have that $g.\chi_{pP}(X) = 1$ if and only if $pP \subseteq g^{-1}X$ if and only if $gpP \subseteq X$. Therefore, $\chi_{pP} = g.\chi_{pP}$ means that for $X \in \cJ_P$, we have $pP \subseteq X$ if and only if $gpP \subseteq X$. Note that $gpP = s(pP)$ lies in $\cJ_P$. Hence, for $X = pP$, we obtain $gpP \subseteq pP$, and for $X = gpP$, we get $pP \subseteq gpP$. Hence there exist $x,y \in P$ with
$$gp = px \ {\rm and} \ p = gpy.$$
So $p = gpy = pxy$ and $gp = px = gpyx$. Thus $xy = yx = e$. Hence $x,y \in P^*$. Since $P^* = \gekl{e}$ by assumption, we must have $x = y = e$, and hence $gp = p$. This implies $g = e$. In other words, for every $e \neq g \in G$, we have $g.\chi_{pP} \neq \chi_{pP}$ for all $p \in P$ such that $\chi_{pP} \in U_{g^{-1}}$. Hence it follows that 
$$
  \menge{\chi \in U_{g^{-1}}}{g.\chi \neq \chi} \ {\rm contains} \ \menge{\chi_{pP} \in U_{g^{-1}}}{p \in P},$$
and the latter set is dense in $U_{g^{-1}}$ as $\menge{\chi_{pP}}{p \in P}$ is dense in $\Omega_P$.
\eproof
Note that $G \curvearrowright \Omega_P$ can be topologically free if $P^* \neq \gekl{e}$. For instance, partial dynamical systems attached to $ax+b$-semigroups over rings of algebraic integers in number fields are shown to be topologically free in \cite{EL}. A generalization of this result is obtained in \cite[Proposition~5.8]{Li5}.

By \cite[Theorem~2.6]{ELQ} and because of Theorem~\ref{P-DG}, we obtain the following
\bcor
Suppose that $P$ is a monoid with trivial units which embeds into a group. Let $I$ be an ideal of $C^*_{\lambda}(P)$.

If $I \cap D_{\lambda}(P) = (0)$, then $I = (0)$.

In other words, a representation of $C^*_{\lambda}(P)$ is faithful if and only if it is faithful on $D_{\lambda}(P)$.
\ecor

Let us now discuss boundary quotients. We start with general inverse semigroups (with or without zero). In many situations, we are not only interested in the reduced C*-algebra of an inverse semigroup, but also in its boundary quotient. This is a notion going back to Exel (see \cites{Ex2,Ex3,Ex4,EGS}). Let us recall the construction. Given a semilattice $E$, let $\widehat{E}_{\max}$ be the subset of $\widehat{E}$ consisting of those $\chi \in \widehat{E}$ such that $\menge{e \in E}{\chi(e) = 1}$ is maximal among all characters $\chi \in \widehat{E}$. Note that if $E$ is a semilattice without zero, then $\widehat{E}_{\max}$ consists of only one element, namely the character $\chi$ satisfying $\chi(e) = 1$ for all $e \in E$. For later purposes, we make the following observation:
\blemma
\label{chimax=0}
Let $E$ be a semilattice with zero, and let $0$ be its distinguished zero element. Suppose that $\chi \in \widehat{E}_{\max}$ satisfies $\chi(e) = 0$ for some $e \in E\reg$. Then there exists $f \in E\reg$ with $\chi(f) = 1$ and $ef = 0$.
\elemma
\bproof
If every $f \in E\reg$ with $\chi(f) = 1$ satisfies $ef \neq 0$, then we can define a filter $\cF$ by defining, for every $\ti{f} \in E\reg$,
$$
  \ti{f} \in \cF \ {\rm if} \ {\rm there} \ {\rm exists} \ f \in E\reg \ {\rm } \ {\rm with} \ \chi(f) = 1 \ {\rm and} \ ef \leq \ti{f}.
$$
It is obvious that $\cF$ is a filter, so that there exists a character $\chi_F \in \widehat{E}$ with $\chi_F^{-1} = \cF$. By construction,
$$
  \menge{f \in E\reg}{\chi(f) = 1} \subseteq \menge{f \in E\reg}{\chi_F(f) = 1},
$$
but $\chi_F(e) = 1$ while $\chi(e) = 0$. This contradicts maximality of
$$
  \menge{f \in E}{\chi(f) = 1}.
$$
\eproof
We define
$$\partial \widehat{E} \defeq \overline{\widehat{E}_{\max}} \subseteq \widehat{E}.$$
Now let $E$ be the semilattice of idempotents in an inverse semigroup $S$. As $\partial \widehat{E} \subseteq \widehat{E}$ is closed, we obtain a short exact sequence
$$0 \to I \to C_0(\widehat{E}) \to C_0(\partial \widehat{E}) \to 0.$$
Now there are two options. We could view $I$ as a subset of $C^*_{\lambda}(S)$ and form the ideal $\spkl{I}$ of $C^*_{\lambda}(S)$ generated by $I$. The boundary quotient in Exel's sense (see \cites{Ex2,Ex3,Ex4,EGS}) is given by
$$\partial C^*_{\lambda}(S) \defeq C^*_{\lambda}(S) / \spkl{I}.$$
Alternatively, we could take the universal groupoid $\cG(S)$ of our inverse semigroup, form its restriction to $\partial \widehat{E}$,
$$
  \cG(S) \, \vert \, \partial \widehat{E} \defeq \menge{\gamma \in \cG(S)}{r(\gamma), s(\gamma) \in \partial \widehat{E}},
$$
and form the reduced groupoid C*-algebra
$$
  C^*_r(\cG(S) \, \vert \, \partial \widehat{E}).
$$
As the canonical homomorphism
$$
  C^*_{\lambda}(S) \cong C^*_r(\cG(S)) \onto C^*_r(\cG(S) \, \vert \, \partial \widehat{E})
$$
contains $\spkl{I}$ in its kernel, we obtain canonical projections
$$
  C^*_{\lambda}(S) \onto C^*_{\lambda}(S) / \spkl{I} \onto C^*_r(\cG(S) \, \vert \, \partial \widehat{E}).
$$
Under an exactness assumption, the second *-homomorphism actually becomes an isomorphism, so that our two alternatives for the boundary quotient coincide. For our purposes, it is more convenient to work with $C^*_r(\cG(S) \, \vert \, \partial \widehat{E})$ because it is, by its very definition, a reduced groupoid C*-algebra, so that groupoid techniques apply.

Now let us assume that our inverse semigroup $S$ admits an idempotent pure partial homomorphism $\sigma: \: S\reg \to G$ to a group $G$. In that situation, we can define the partial dynamical system $G \curvearrowright \widehat{E}$ (see \S~\ref{ss:PDS}) and identify $\cG(S)$ with the partial transformation groupoid $G \ltimes \widehat{E}$ (see Lemma~\ref{GS=GxhatE}). We have the following
\blemma
\label{bd_Ginv}
Let $S$ be an inverse semigroup with an idempotent pure partial homomorphism to a group $G$. Let $G \curvearrowright \widehat{E}$ be its partial dynamical system. Then $\partial \widehat{E}$ is $G$-invariant.
\elemma
\bproof
Let us first show that for every $g \in G$,
$$g.(U_{g^{-1}} \cap \widehat{E}_{\max}) \subseteq U_g \cap \widehat{E}_{\max}.$$
Take $\chi \in \widehat{E}_{\max}$ with $\chi(s^{-1}s) = 1$ for some $s \in S$ with $\sigma(s) = g$. Then $g.\chi(e) = \chi(s^{-1}es)$. Assume that $g.\chi \notin \widehat{E}_{\max}$. This means that there is $\psi \in \widehat{E}_{\max}$ such that $\psi(e) = 1$ for all $e \in E$ with $g.\chi(e) = 1$, and there exists $f \in E$ with $\psi(f) = 1$ but $\chi(s^{-1}fs) = 0$. Then $\psi \in U_g$ since $g.\chi(ss^{-1}) = 1$, which implies $\psi(ss^{-1}) = 1$. Consider $g^{-1}.\psi$ given by $g^{-1}.\psi(e) = \psi(ses^{-1})$. Then for every $e \in E$, $\chi(e) = 1$ implies $\chi(s^{-1}ses^{-1}s) = 1$, hence $\chi(s^{-1}(ses^{-1})s) = 1$, so that $g^{-1}.\psi(e) = \psi(ses^{-1}) = 1$. But $\chi(s^{-1}fs) = 0$ and $g^{-1}.\psi(s^{-1}fs) = \psi(ss^{-1}fss^{-1}) = \psi(f) = 1$. This contradicts $\chi \in \widehat{E}_{\max}$. Hence $g.(U_{g^{-1}} \cap \widehat{E}_{\max}) \subseteq U_g \cap \widehat{E}_{\max}$.

To see that
$$g.(U_{g^{-1}} \cap \partial \widehat{E}) \subseteq U_g \cap \partial \widehat{E},$$
let $\chi \in U_{g^{-1}} \cap \partial \widehat{E}$ and choose a net $(\chi_i)_i$ in $\widehat{E}_{\max}$ with $\lim_i \chi_i = \chi$. As $U_{g^{-1}}$ is open, we may assume that all the $\chi_i$ lie in $U_{g^{-1}}$. Then $g.\chi_i \in \widehat{E}_{\max}$, and $\lim_i g.\chi_i = g.\chi$. This implies $g.\chi \in \partial \widehat{E}$.
\eproof

\bcor
\label{bdquot=cropro}
In the situation of Lemma~\ref{bd_Ginv}, we have canonical isomorphisms
$$
  \cG(S) \, \vert \, \partial \widehat{E} \cong G \ltimes \partial \widehat{E}
$$
and
$$
  C^*_r(\cG(S) \, \vert \, \partial \widehat{E}) \cong C_0(\partial \widehat{E}) \rtimes_r G.
$$
\ecor
\bproof
The first identification follows immediately from Lemma~\ref{bd_Ginv}, while the second one is a consequence of the first one and Theorem~\ref{THM:C*GPD=C*PDS_red}.
\eproof

Let us now specialize to the case where $S$ is the left inverse hull of a left cancellative semigroup $P$. First, we observe the following:
\blemma
\label{bdOmegaInOmega}
We have $\partial \widehat{\cJ_P} \subseteq \Omega_P$.
\elemma
\bproof
Let $X, X_1, \dotsc, X_n \in \cJ_P$ satisfy $X = \bigcup_{i=1}^n X_i$. Then for $\chi \in (\widehat{\cJ_P})_{\max}$, $\chi(X_i) = 0$ implies that there exists $X_i' \in \cJ$ with $\chi(X_i') = 1$ and $X_i \cap X_i' = \emptyset$ (see Lemma~\ref{chimax=0}). Thus if $\chi(X_i) = 0$ for all $1 \leq i \leq n$, then let $X_i'$, $1 \leq i \leq n$ be as above. Then for $X' = \bigcap_{i=1}^n X_i'$, $\chi(X') = 1$ and $X \cap X' = \emptyset$. Thus $\chi(X) = 0$. This shows $(\widehat{\cJ_P})_{\max} \subseteq \Omega_P$. As $\Omega_P$ is closed, we conclude that $\partial \widehat{\cJ_P} \subseteq \Omega_P$.
\eproof

\bdefin
We write $\partial \Omega_P \defeq \partial \widehat{\cJ_P}$.
\edefin

For simplicity, let us now restrict to semigroups which embed into groups.
\bdefin
We call $C^*_r(\cG(I_l(P)) \, \vert \, \partial \Omega_P)$ the boundary quotient of $C^*_{\lambda}(P)$, and denote it by $\partial C^*_{\lambda}(P)$.
\edefin

Note that by Corollary~\ref{bdquot=cropro}, given a semigroup $P$ embedded into a group $G$, we have a canonical isomorphism
$$
  \partial C^*_{\lambda}(P) \cong C(\partial \Omega_P) \rtimes_r G.
$$

Let us discuss some examples. Assume that our semigroup $P$ is cancellative, and that it is left reversible, i.e., $pP \cap qP \neq \emptyset$ for all $p, q \in P$. This is for instance the case for positive cones in totally ordered groups. Given such a semigroup, we know because of Lemma~\ref{char-leftrev} that $\cJ_P$ is a semilattice without zero, so that $(\widehat{\cJ_P})_{\max}$ degenerates to a point. Therefore, $\partial \Omega_P$ degenerates to a point. Hence it follows that the boundary quotient $\partial C^*_{\lambda}(P)$ coincides with the reduced group C*-algebra of the group of right quotients of $P$.

For the non-abelian free monoid $\Nz * \Nz$ on two generators, the boundary quotient $\partial C^*_{\lambda}(\Nz * \Nz)$ is canonically isomorphic to the Cuntz algebra $\cO_2$. More generally, boundary quotients for right-angled Artin monoids are worked out and studied in \cite{CrLa2}.

Given an integral domain $R$, the boundary quotient $\partial C^*_{\lambda}(R \rtimes R\reg)$ of the $ax+b$-semigroup over $R$ is canonically isomorphic to the ring C*-algebra $\fA_r[R]$ of $R$ (see \cites{CuLi,CuLiX,LiRing}). It is given as follows:

Consider the Hilbert space $\ell^2 R$ with canonical orthonormal basis $\menge{\delta_x}{x \in R}$. For every $a \in R\reg$, define $S_a (\delta_x) \defeq \delta_{ax}$, and for every $b \in R$, define $U^b(\delta_x) \defeq \delta_{b+x}$. Then the ring C*-algebra of $R$ is the C*-algebra generated by these two families of operators, i.e,
$$
  \fA_r[R] \defeq C^*(\menge{S_a}{a \in R\reg} \cup \menge{U^b}{b \in R}) \subseteq \cL(\ell^2 R).
$$
We refer to \cites{CuLi,CuLiX,LiRing} and also \cite[\S~8.3]{LiNuc} for details.

Let us now establish structural properties for boundary quotients. From now on, let us suppose that our semigroup $P$ embeds into a group $G$.

\blemma
\label{Lem:bd-min}
$\partial \Omega_P$ is the minimal non-empty closed $G$-invariant subspace of $\widehat{\cJ_P}$.
\elemma
\bproof
Let $C \subseteq \widehat{\cJ_P}$ be non-empty, closed and $G$-invariant. Let $\chi \in(\widehat{\cJ_P})_{\max}$ be arbitrary, and choose $X \in \cJ_P$ with $\chi(X) = 1$. Choose $p \in X$ and $\chi \in C$. As $U_{p^{-1}} = \widehat{\cJ_P}$, we can form $p.\chi$, and we know that $p.\chi \in C$. We have $p.\chi(pP) = \chi(P) = 1$, so that $p.\chi(X) = 1$ as $p \in X$ implies $pP \subseteq X$ ($X$ is a right ideal). Set $\chi_X \defeq p.\chi$. Consider the net $(\chi_X)_X$ indexed by $X \in \cJ$ with $\chi(X) = 1$, ordered by inclusion. Passing to a convergent subnet if necessary, we may assume that $\lim_X \chi_X$ exists. But it is clear because of $\chi \in (\widehat{\cJ_P})_{\max}$ that $\lim_X \chi_X = \chi$. As $\chi_X \in C$ for all $X$, we deduce that $\chi \in C$. Thus $(\widehat{\cJ_P})_{\max} \subseteq C$, and hence $\partial \Omega_P \subseteq C$.
\eproof

In particular, $\partial \Omega_P$ is the minimal non-empty closed $G$-invariant subspace of $\Omega_P$. Another immediate consequence is
\bcor
\label{Cor:bdGPD-min}
The transformation groupoid $G \ltimes \partial \Omega_P$ is minimal.
\ecor

To discuss topological freeness of $G \curvearrowright \partial \Omega_P$, let
$$G_0 = \menge{g \in G}{X \cap g P \neq \emptyset \neq X \cap g^{-1} P \ {\rm for} \ {\rm all} \ \emptyset \neq X \in \cJ_P},$$
as in \cite[\S~7.3]{LiNuc}. Clearly,
$$
  G_0 = \menge{g \in G}{pP \cap g P \neq \emptyset \neq pP \cap g^{-1} P \ {\rm for} \ {\rm all} \ p \in P}.
$$
Furthermore, we have the following
\blemma
$G_0$ is a subgroup of $G$.
\elemma
\bproof
Take $g_1$, $g_2$ in $G_0$. Then for all $\emptyset \neq X \in \cJ$, we have
\bglnoz
  && ((g_1 g_2) P) \cap X = g_1 ((g_2 P) \cap (g_1^{-1} X))\\
  &\supseteq& g_1 ((g_2 P) \cap (g_1^{-1} X)) \cap (g_1 P) = g_1 ((g_2 P) \cap ((g_1^{-1} X) \cap P)).
\eglnoz
Now
$$(g_1^{-1} X) \cap P = g_1^{-1} (X \cap (g_1 P)) \neq \emptyset.$$
Thus there exists $x \in P$ such that $x \in (g_1^{-1} X) \cap P$. Hence $xP \subseteq (g_1^{-1} X) \cap P$. Thus
$$
  \emptyset \neq g_1 ((g_2 P) \cap (xP)) \subseteq ((g_1 g_2) P) \cap X.
$$
\eproof

\bprop
\label{Prop_bd-tf}
$G \curvearrowright \partial \Omega_P$ is topologically free if and only if $G_0 \curvearrowright \partial \Omega_P$ is topologically free. 
\eprop
\bproof
\an{$\Rarr$} is clear. For \an{$\Larr$}, assume that $G_0 \curvearrowright \partial \Omega_P$ is topologically free, and suppose that $G \curvearrowright \partial \Omega_P$ is not topologically free, i.e., there exists $g \in G$ and $U \subseteq U_{g^{-1}} \cap \partial \Omega_P$ such that $g.\chi = \chi$ for all $\chi \in U$. As $\overline{(\widehat{\cJ_P})_{\max}} = \partial \Omega_P$, we can find $\chi \in U_{g^{-1}} \cap (\widehat{\cJ_P})_{\max}$ with $g.\chi = \chi$.

For every $X \in \cJ_P$ with $\chi(X) = 1$, choose $x \in X$ and $\psi_X \in (\widehat{\cJ_P})_{\max}$ with $\psi_X(xP) = 1$, so that $\psi_X(X) = 1$. Consider the net $(\psi_X)_X$ indexed by $X \in \cJ_P$ with $\chi(X) = 1$, ordered by inclusion. Passing to a convergent subnet if necessary, we may assume that $\lim_X \psi_X = \chi$. As $U$ is open, we may assume that $\psi_X \in U$ for all $X$. Then $\psi_X(xP) = 1$ implies that $\psi_X \in U_x \cap U$.

Hence for sufficiently small $X \in \cJ_P$ with $\chi(X) = 1$, there exists $x \in X$ such that $x^{-1}.(U_x \cap U)$ is a non-empty open subset of $\partial \Omega_P$. We conclude that $(x^{-1}gx).\psi = \psi$ for all $\psi \in x^{-1}.(U_x \cap U)$. This implies that $x^{-1}gx \notin G_0$ as $G_0 \curvearrowright \partial \Omega_P$ is topologically free. So there exists $p \in P$ with
$$pP \cap x^{-1}gx P = \emptyset \ {\rm or} \ pP \cap x^{-1}g^{-1}x P = \emptyset.$$
Let $\chi_X \in (\widehat{\cJ_P})_{\max}$ satisfy $\chi_X(xpP) = 1$. If $pP \cap x^{-1}gx P = \emptyset$, then 
$$
xpP \cap gx P = \emptyset, \ {\rm so} \ {\rm that} \ xpP \cap g^{-1}xp P = \emptyset.
$$
Hence $g.\chi_X \neq \chi_X$ if $\chi_X \in U_{g^{-1}}$. If $pP \cap x^{-1}g^{-1}x P = \emptyset$, then 
$$
xpP \cap g^{-1}x P = \emptyset, \ {\rm so} \ {\rm that} \ xpP \cap g^{-1}xp P = \emptyset.$$
Again, $g.\chi_X \neq \chi_X$ if $\chi_X \in U_{g^{-1}}$.

For every sufficiently small $X \in \cJ_P$ with $\chi(X) = 1$, we can find $x \in X$ and $\chi_X$ as above. Hence we can consider the net $(\chi_X)_X$ as above, and assume after passing to a convergent subnet that $\lim_X \chi_X = \chi$. As $\chi \in U \subseteq U_{g^{-1}} \cap \partial \Omega_P$, it follows that $\chi_X \in U \subseteq U_{g^{-1}} \cap \partial \Omega_P$ for sufficiently small $X$. So we obtain $g.\chi_X \neq \chi_X$, although $g$ acts trivially on $U$. This is a contradiction.
\eproof

\bcor
\label{Cor:bd-tf->simple}
If $G_0 \curvearrowright \partial \Omega_P$ is topologically free, then $\partial C^*_{\lambda}(P)$ is simple.
\ecor
\bproof
This follows from Lemma~\ref{Lem:bd-min}, Proposition~\ref{Prop_bd-tf} and \cite[Chapter~II, Proposition~4.6]{Ren1}. 
\eproof

We present a situation where Corollary~\ref{Cor:bd-tf->simple} applies. Recall that we introduced the notion of \an{completeness for $\curvearrowright_R$} for presentations after Lemma~\ref{pPCAPqP}. Moreover, a pair $P \subseteq G$ consisting of a monoid $P$ embedded into a group $G$ is called quasi-lattice ordered (see \cite{Nica}) if $P$ has trivial units $P^* = \gekl{e}$ and for every $g \in G$ with $gP \cap P \neq \emptyset$, we can find an element $p \in P$ such that $gP \cap P = pP$.
\btheo
\label{THM_sgp-rep->G0=trivial}
Let $P = \spkl{\Sigma, R}^+$ be a monoid given by a presentation $(\Sigma,R)$ which is complete for $\curvearrowright_R$, in the sense of \cite{Deh}. Assume that for all $u \in \Sigma$, there is $v \in \Sigma$ such that there is no relation of the form $u \cdots = v \cdots$ in $R$. Also, suppose that $P$ embeds into a group $G$ such that $P \subseteq G$ is quasi-lattice ordered in the sense of \cite{Nica}. Then $G_0 = \gekl{e}$ and $\partial C^*_{\lambda}(P)$ is simple.
\etheo
\bproof
In view of Corollary~\ref{Cor:bd-tf->simple}, it suffices to prove $G_0 = \gekl{e}$. Let $g \in G_0$. Assume that $gP \cap P \neq P$. Then $g \in G_0$ implies that this intersection is not empty. Hence, we must have $gP \cap P = pP$ for some $p \in P$ because $P \subseteq G$ is quasi-lattice ordered. If $p \neq e$, then there exists $u \in \Sigma$ with $pP \subseteq uP$. By assumption, there exists $v \in \Sigma$ such that no relation in $R$ is of the form $u \cdots = v \cdots$. Because $(\Sigma,R)$ is complete for $\curvearrowright_R$, we know that $uP \cap vP = \emptyset$ (see \cite[Proposition~3.3]{Deh}), so that $gP \cap vP = \emptyset$. This contradicts $g \in G_0$. Hence, we must have $gP \cap P = P$, and similarly, $g^{-1} P \cap P = P$. These two equalities imply $g \in P^*$. But $P^* = \gekl{e}$ because $P \subseteq G$ is quasi-lattice ordered. Thus $g = e$.
\eproof
Theorem~\ref{THM_sgp-rep->G0=trivial} implies that for every right-angled Artin monoid $A_{\Gamma}^+$ (see \S~\ref{ex:presentations}) with the property that $(A_{\Gamma},A_{\Gamma}^+)$ is graph-irreducible in the sense of \cite{CrLa2}, the boundary quotient $\partial C^*_{\lambda}(A_{\Gamma}^+)$ is simple.

Moreover, assume that we have a cancellative semigroup. By going over to the opposite semigroup, the left regular representation becomes the right regular representation. In this way, our discussion about C*-algebra generated by left regular representations applies to C*-algebras of right regular representations. In particular, we can define boundary quotients for C*-algebras generated by right regular representations of semigroups. For instance, for the Thompson monoid
$$F^+ = \spkl{x_0, x_1, \dotsc \ \vert \ x_n x_k = x_k x_{n+1} \ {\rm for} \ k<n}^+,$$
it is easy to see that Theorem~\ref{THM_sgp-rep->G0=trivial} applies to the opposite monoid, so that the boundary quotient of the C*-algebra generated by the right regular representation of $F^+$ is simple.

We now turn to the property of pure infiniteness. As we mentioned, the boundary quotient $\partial C^*_{\lambda}(\Nz * \Nz)$ is isomorphic to $\cO_2$, a purely infinite C*-algebra. We will now see that this is not a coincidence.

First of all, it is easy to see that for a partial dynamical system $G \curvearrowright X$, the transformation groupoid $G \ltimes X$ is purely infinite in the sense of \cite{Ma2015} if and only if every compact open subset of $X$ is $(G,\cC \cO)$-paradoxical in the sense of \cite[Definition~4.3]{GS}, where $\cC \cO$ is the set of compact open subsets of $X$. We recall that a non-empty subset $V \subseteq X$ is called 
$(G,\cC \cO)$-paradoxical in \cite[Definition~4.3]{GS} if there exist
$$
  V_1, \dotsc, V_{n+m} \in \cO \ {\rm and} \ t_1, \dotsc, t_{n+m} \in G
$$
such that
$$
  \bigcup_{i=1}^n V_i = V = \bigcup_{i=n+1}^m V_i,
$$
and
$$
  V_i \in U_{t_i^{-1}}, \ t_i.V_i \subseteq V, \ {\rm and} \ t_i.V_i \cap t_j.V_j = \emptyset \ {\rm for} \ {\rm all} \ i \neq j.
$$

\btheo
\label{bdP_pi}
The groupoid $G \ltimes \partial \Omega_P$ is purely infinite if and only if there exist $p, q \in P$ with $pP \cap qP = \emptyset$.
\etheo

\bproof
Obviously, if $pP \cap qP \neq \emptyset$ for all $p, q \in P$, then $\partial \Omega_P$ degenerates to a point. 

Let us prove the converse. Every compact open subset of $\widehat{\cJ_P}$ can be written as a disjoint union of basic open sets 
$$U = \menge{\psi \in \partial \Omega_P}{\psi(X) = 1, \psi(X_1) = \dotso = \psi(X_n) = 0},$$
for some $X, X_1, \dotsc, X_n \in \cJ_P$. Hence it suffices to show that $U$ is $(G,\cC \cO)$-paradoxical. Since $(\widehat{\cJ_P})_{\max}$ is dense in $\partial \Omega_P$, there exists $\chi \in (\widehat{\cJ_P})_{\max}$ with $\chi \in U$. As $\chi$ lies in $(\widehat{\cJ_P})_{\max}$, $\chi(X_i) = 0$ implies that there exists $Y_i \in \cJ$ with $X_i \cap Y_i = \emptyset$ and $\chi(Y_i) = 1$ (see Lemma~\ref{chimax=0}). Let
$$Y \defeq X \cap \bigcap_{i=1}^n Y_i.$$
Certainly, $Y \neq \emptyset$ as $\chi(Y) = 1$. Moreover, for every $\psi \in \partial \Omega_P$, $\psi(Y) = 1$ implies $\psi \in U$. Now choose $x \in Y$. By assumption, we can find $p, q \in P$ with $pP \cap qP = \emptyset$. For $\psi \in \partial \Omega_P$, $xp.\psi(xpP) = \psi(P) = 1$. Similarly, for all $\psi \in \partial \Omega_P$, we have $xq.\psi(xqP) = 1$. Thus 
\bglnoz
&& xp.U \subseteq xp.\partial \Omega_P \subseteq U, \ xq.U \subseteq xq.\partial \Omega_P \subseteq U \\
&{\rm and}& (xp.U) \cap (xq.U) \subseteq (xp.\partial \Omega_P) \cap (xq.\partial \Omega_P) = \emptyset
\eglnoz
since $xpP \cap xqP = \emptyset$.
\eproof

\bcor
\label{Cor:bdquot-pis}
If $P$ is not the trivial monoid, $P \neq \gekl{e}$, and if $G_0 \curvearrowright \partial \Omega_P$ is topologically free, then the boundary quotient $\partial C^*_{\lambda}(P)$ is a purely infinite simple C*-algebra.
\ecor
\bproof
First of all, by Corollary~\ref{Cor:bd-tf->simple}, the boundary quotient is simple.

Furthermore, we observe that our assumptions that $P \neq \gekl{e}$ and that $G_0$ acts topologically freely on $\partial \Omega_P$ imply that $P$ is not left reversible: If $P$ were left reversible, then $\partial \Omega_P$ would consist of only one point. Also, if $P$ were left reversible, then we would have $P \subseteq G_0$. Since every element in $P$ obviously leaves $\partial \Omega_P$ fixed, and by our assumption that $P \neq \gekl{e}$, we conclude that $G_0$ cannot act topologically freely on $\partial \Omega_P$ if $P$ were left reversible. Hence Theorem~\ref{bdP_pi} implies that the groupoid $G \ltimes \partial \Omega_P$ is purely infinite.

This, together with \cite[Theorem~4.4]{GS}, implies that the boundary quotient $\partial C^*_{\lambda}(P)$ is purely infinite. This completes our proof.
\eproof

\bcor
If $P$ is not the trivial monoid, $P \neq \gekl{e}$, if $G \ltimes \partial \Omega_P$ is amenable, and if $G_0 \curvearrowright \partial \Omega_P$ is topologically free, then the boundary quotient $\partial C^*_{\lambda}(P)$ is a unital UCT Kirchberg algebra.
\ecor
\bproof
By assumption, our semigroups are countable, so that all the C*-algebras we construct are separable. Clearly, the boundary quotient $\partial C^*_{\lambda}(P)$ is unital. 

Since $G \ltimes \partial \Omega_P$ is amenable, the boundary quotient $\partial C^*_{\lambda}(P)$ is nuclear and satisfies the UCT. 

Now our claim follows from Corollary~\ref{Cor:bdquot-pis}
\eproof
Note that this shows that \cite[Corollary~7.23]{LiNuc} holds without the independence and the Toeplitz condition.

Let us now study simplicity of reduced semigroup C*-algebras. Let $P$ be a semigroup which embeds into a group. If $C^*_{\lambda}(P)$ is simple, then the groupoid $G \ltimes \Omega_P$ must be minimal, as $C^*_{\lambda}(P) \cong C^*_r(G \ltimes \Omega_P)$ (see the isomorphism~\eqref{P-->GxOmega}). In particular, we must have $\Omega_P = \partial \Omega_P$. This equality can be characterized in terms of the semigroup as follows:
\blemma
\label{Omega=bd}
Let $P$ be a monoid. We have $\Omega_P = \partial \Omega_P$ if and only if for every $X_1, \dotsc, X_n \in \cJ_P$ with $X_i \subsetneq P$ for all $1 \leq i \leq n$, there exists $p \in P$ with $pP \cap X_i = \emptyset$ for all $1 \leq i \leq n$.
\elemma
\bproof
Let $\chi_P$ be the character in $\Omega_P$ determined by $\chi_P(X) = 1$ if and only if $X = P$, for all $X \in \cJ_P$. Such a character exists in $\widehat{\cJ_P}$, and our assumption that $P$ has an identity element ensures that $\chi_P$ lies in $\Omega_P$. This is because an equation of the form
$$
  P = \bigcup_{i=1}^n X_i
$$
for some $X_i \in \cJ_P$ implies that $X_i = P$ for some $1 \leq i \leq n$, as one of the $X_i$ must contain the identity element.

First, we claim that $\Omega_P = \partial \Omega_P$ holds if and only if $\chi_P$ lies in $\partial \Omega_P$. This is certainly necessary. It is also sufficient as $\partial \Omega_P$ is $G$-invariant, and
$$
  \menge{p.\chi_P = \chi_{pP}}{p \in P}
$$
is dense in $\Omega_P$ (see Lemma~\ref{chi_P:dense}).

Now basic open subsets containing $\chi_P$ are of the form
$$
  U(P; X_1, \dotsc, X_n) = \menge{\chi \in \Omega_P}{\chi(X_1) = \dotso = \chi(X_n) = 0},
$$
for $X_1, \dotsc, X_n \in \cJ_P$ with $X_i \subsetneq P$ for all $1 \leq i \leq n$.

$\chi_P \in \partial \Omega_P$ if and only if $\chi_P \in \overline{(\widehat{\cJ_P})_{\max}}$ if and only if for all $X_1, \dotsc, X_n \in \cJ_P$ with $X_i \subsetneq P$ for all $1 \leq i \leq n$, there is $\chi \in (\widehat{\cJ_P})_{\max}$ with $\chi \in U(P; X_1, \dotsc, X_n)$. Hence it follows that our proof is complete once we show that there exists $\chi \in (\widehat{\cJ_P})_{\max}$ with $\chi \in U(P; X_1, \dotsc, X_n)$ if and only if there exists $p \in P$ with $pP \cap X_i = \emptyset$ for all $1 \leq i \leq n$.

For \an{$\Rarr$}, assume that $\chi \in (\widehat{\cJ_P})_{\max}$ lies in $\chi \in U(P; X_1, \dotsc, X_n)$. Then $\chi(X_i) = 0$ for all $1 \leq i \leq n$. But this means that there must exist $Y_i \in \cJ_P$, for $1 \leq i \leq n$, such that $\chi(Y_i) = 1$ and $X_i \cap Y_i = \emptyset$ (see Lemma~\ref{chimax=0}). Take the intersection
$$
  Y \defeq \bigcap_{i=1}^n Y_i.
$$
As $\chi(Y) = 1$, $Y$ is not empty. Therefore, we may choose some $p \in Y$. Obviously, $pP \subseteq Y$ as $Y$ is a right ideal. Moreover, for every $1 \leq i \leq n$, we have
$$
  X_i \cap pP \subseteq X_i \cap Y \subseteq X_i \cap Y_i = \emptyset.
$$

For \an{$\Larr$}, suppose that there exists $p \in P$ with $pP \cap X_i = \emptyset$ for all $1 \leq i \leq n$. An easy application of Zorn's Lemma yields a character $\chi \in (\widehat{\cJ_P})_{\max}$ with $\chi(pP) = 1$. Hence $\chi(X_i) = \emptyset$ for all $1 \leq i \leq n$, and it follows that $\chi$ lies in $U(P; X_1, \dotsc, X_n)$.
\eproof

Let us derive some immediate consequences.

\bcor
\label{Cor:char-min_tf}
If $G \curvearrowright \Omega_P$ is topologically free, then $C^*_{\lambda}(P)$ is simple if and only if for every $X_1, \dotsc, X_n \in \cJ_P$ with $X_i \subsetneq P$ for all $1 \leq i \leq n$, there exists $p \in P$ with $pP \cap X_i = \emptyset$ for all $1 \leq i \leq n$.
\ecor
\bproof
This follows immediately from Corollary~\ref{Cor:bdGPD-min} and Lemma~\ref{Omega=bd} (see also Corollary~\ref{Cor:bd-tf->simple}).
\eproof

\bcor
\label{Cor:char-min_P*=e}
Let $P$ be a monoid with identity $e$, and suppose that $P$ embeds into a group. Suppose that $P$ has trivial units $P^* = \gekl{e}$. Then $C^*_{\lambda}(P)$ is simple if and only if for every $X_1, \dotsc, X_n \in \cJ_P$ with $X_i \subsetneq P$ for all $1 \leq i \leq n$, there exists $p \in P$ with $pP \cap X_i = \emptyset$ for all $1 \leq i \leq n$.
\ecor
\bproof
This follows from Corollary~\ref{Cor:char-min_tf} and Lemma~\ref{P*=e->TF}.
\eproof

As an example, the countable free product $P = *_{i=1}^{\infty} \Nz$ satisfies the criterion in Lemma~\ref{Omega=bd}. Moreover, we obviously have $P^* = \gekl{e}$. Hence Corollary~\ref{Cor:char-min_P*=e} applies, and we deduce that $C^*_{\lambda}(*_{i=1}^{\infty} \Nz)$ is simple. Actually, $C^*_{\lambda}(*_{i=1}^{\infty} \Nz)$ is canonically isomorphic to the Cuntz algebra $\cO_{\infty}$.

\section{The Toeplitz condition}
\label{sec:Toeplitz}

So far, we were able to derive all our results about semigroup C*-algebras just using descriptions as partial crossed products. However, it turns out that when we want to compute K-theory or the primitive ideal space, we need descriptions (at least up to Morita equivalence) as ordinary crossed products, attached to globally defined dynamical systems. Let us now introduce a criterion which guarantees such descriptions as ordinary crossed products.

\bdefin
Let $P \subseteq G$ be a semigroup embedded into a group $G$. We say that $P \subseteq G$ satisfies the Toeplitz condition (or simply that $P \subseteq G$ is Toeplitz) if for every $g \in G$ with $g^{-1}P \cap P \neq \emptyset$, the partial bijection
$$
  g^{-1}P \cap P \to P \cap gP, \, x \ma gx
$$
lies in the inverse semigroup $I_l(P)$.
\edefin

We can also think of $I_l(P)$ as partial isometries on $\ell^2 P$. In this picture, we can give an equivalent characterization of the Toeplitz condition. First, using the embedding $P \subseteq G$, we pass to the bigger Hilbert space $\ell^2 G$. Let $1_P$ be the characteristic function of $P$, viewed as an element in $\ell^{\infty}(G)$. Moreover, let $\lambda$ be the left regular representation of $G$ on $\ell^2 G$. Then $P \subseteq G$ is Toeplitz if and only if for every $g \in G$ with $1_P \lambda_g 1_P \neq 0$, we can write $1_P \lambda_g 1_P$ as a finite product of isometries and their adjoints from the set
$$
  \menge{V_p}{p \in P} \cup \menge{V_q^*}{q \in P}.
$$

Let us now explain why the reduced semigroup C*-algebra $C^*_{\lambda}(P)$ is a full corner in an ordinary crossed product if $P \subseteq G$ is Toeplitz. In terms of the partial dynamical system $G \curvearrowright \Omega_P$, this amounts to showing that if $P \subseteq G$ is Toeplitz, then $G \curvearrowright \Omega_P$ has an enveloping action, in the sense of \cite{Aba1}, on a locally compact Hausdorff space. This is because if $P \subseteq G$ is Toeplitz, then $g^{-1}P \cap P$ lies in the semilattice $\cJ_P$. Hence, for every $g \in G$,
$$
  U_{g^{-1}} = \menge{\chi \in \Omega_P}{\chi(g^{-1}P \cap P) = 1},
$$
since among all $s \in I_l(P)\reg$ with $\sigma(s) = g$, $g^{-1}P \cap P$ is the maximal domain. This means that for every $g \in G$, the subspace $U_{g^{-1}}$ is clopen. Whenever this is the case, our partial dynamical system will have an enveloping action on a locally compact Hausdorff space. This follows easily from \cite{Aba1}.

In the following, we give a direct argument describing $C^*_{\lambda}(P)$ as a full corner in an ordinary crossed product in a very explicit way. First, we introduce some notation.

Fix an embedding $P \subseteq G$ of a semigroup $P$ into a group $G$.
\bdefin
We let $\cJ_{P \subseteq G}$ be the smallest $G$-invariant semilattice of subsets of $G$ containing $\cJ_P$.
\edefin

\blemma
\label{Lem:J_PinG}
We have
\bgl
\label{J_PinG}
  \cJ_{P \subseteq G} = \menge{\bigcap_{i=1}^n g_iP}{g_i \in G}.
\egl
If $P \subseteq G$ is Toeplitz, then 
$$
  \cJ_P\reg = \menge{\emptyset \neq Y \cap P}{Y \in \cJ_{P \subseteq G}\reg}.
$$
\elemma
\bproof
Clearly,
$$
  \menge{\bigcap_{i=1}^n g_iP}{g_i \in G}
$$
is a $G$-invariant semilattice of subsets of $G$. It remains to show that it contains $\cJ_P$. It certainly includes the subset $P$ of $G$. Moreover, for every subset $X \in P$ and all $p, q \in P$, we have $p(X) = pX$ and $q^{-1}(X) = q^{-1}X \cap P$. Here, $pX$ and $q^{-1}X$ are products taken in $G$. Therefore, we see that $\cJ_{P \subseteq G}$ is closed under left multiplication and pre-images under left multiplication. But $\cJ_P$ may be characterized as the smallest semilattice of subsets of $P$ containing $P$ and closed under left multiplication and pre-images under left multiplication. Therefore, $\cJ_P$ is contained in $\cJ_{P \subseteq G}$.

Our argument above also shows that we always have 
$$
  \cJ_P\reg \subseteq \menge{\emptyset \neq Y \cap P}{Y \in \cJ_{P \subseteq G}\reg}.
$$
Now let us assume that $P \subseteq G$ is Toeplitz, and let us prove \an{$\supseteq$}. By assumption, the partial bijection
$$
  g^{-1}P \cap P \to P \cap gP, \, x \ma gx
$$
lies in $I_l(P)$ as long as $g^{-1}P \cap P \neq \emptyset$. Therefore, as long as $g^{-1}P \cap P \neq \emptyset$, the image of this partial bijection, $P \cap gP$, lies in $\cJ_P$. Hence it follows, because of \eqref{J_PinG}, that
$$
  \menge{\emptyset \neq Y \cap P}{Y \in \cJ_{P \subseteq G}\reg}
$$
is contained in $\cJ_P\reg$.
\eproof

\bdefin
\label{DEF:D}
We define
$$
  D_{P \subseteq G} \defeq C^*(\menge{1_Y}{Y \in \cJ_{P \subseteq G}}) \subseteq \ell^{\infty}(G).
$$
\edefin

Obviously, $D_{P \subseteq G}$ is $G$-invariant with respect to the canonical action of $G$ on $\ell^{\infty}(G)$ by left multiplication. Therefore, we can form the crossed product $D_{P \subseteq G} \rtimes_r G$. It is easy to see, and explained in \cite[\S~2.5]{CEL1}, that we can identify this crossed product $D_{P \subseteq G} \rtimes_r G$ with the C*-algebra
$$
  C^*(\menge{1_Y \lambda_g}{Y \in \cJ_{P \subseteq G}, \, g \in G}) \subseteq \cL(\ell^2 G)
$$
concretely represented on $\ell^2 G$.

\bprop
\label{Prop:Toeplitz->C*_r=fullcorner}
In the situation above, $1_P$ is a full projection in $D_{P \subseteq G} \rtimes_r G$.

If $P \subseteq G$ is Toeplitz, then
\bgl
\label{C*_r=fullcorner}
  C^*_{\lambda}(P) = 1_P (D_{P \subseteq G} \rtimes_r G) 1_P.
\egl
In particular, $C^*_{\lambda}(P)$ is a full corner in $D_{P \subseteq G} \rtimes_r G$.
\eprop
Equation~\eqref{C*_r=fullcorner} is meant as an identity of sub-C*-algebras of $\cL(\ell^2 G)$.
\bproof
As the linear span of elements of the form
$$
  1_Y \lambda_g, \ Y \in \cJ_{P \subseteq G}, \, g \in G
$$
is dense in $D_{P \subseteq G} \rtimes_r G$, it suffices to show that, for all $Y \in \cJ_{P \subseteq G}$ and $g \in G$,
$$
  1_Y \lambda_g \in \rukl{D_{P \subseteq G} \rtimes_r G} 1_P \rukl{D_{P \subseteq G} \rtimes_r G}
$$
in order to show that $1_P$ is a full projection. Let
$$
  Y = \bigcap_{i=1}^n g_i P.
$$
Then
$$
  1_Y \lambda_g = \rukl{\lambda_{g_1} 1_P} 1_P \rukl{\lambda_{g_1}^* 1_{\bigcap_{i=2}^n g_i P} \lambda_g}
$$
lies in
$$
  \rukl{D_{P \subseteq G} \rtimes_r G} 1_P \rukl{D_{P \subseteq G} \rtimes_r G}.
$$

Let us prove that
$$
  C^*_{\lambda}(P) = 1_P (D_{P \subseteq G} \rtimes_r G) 1_P
$$
if $P \subseteq G$ is Toeplitz. First, observe that \an{$\subseteq$} always holds as for all $p \in P$, we have $V_p = 1_P \lambda_p 1_P$. Conversely, it suffices to show that for every $Y \in \cJ_{P \subseteq G}$ and $g \in G$, $1_P 1_Y \lambda_g 1_P$ lies in $C^*_{\lambda}(P)$. But
$$
  1_P 1_Y \lambda_g 1_P = \rukl{1_P 1_Y 1_P} \rukl{1_P \lambda_g 1_P},
$$
and $1_P 1_Y 1_P$ lies in $C^*_{\lambda}(P)$ as $P \cap Y$ lies in $\cJ_P$ as long as it is not empty by Lemma~\ref{Lem:J_PinG}, and $1_P \lambda_g 1_P$ lies in $C^*_{\lambda}(P)$ because $P \subseteq G$ is Toeplitz.
\eproof

Let us discuss some examples. First, assume that $P$ is cancellative, and right reversible, i.e., $Pp \cap Pq \neq \emptyset$ for all $p, q \in P$. Then $P$ embeds into its group $G$ of left quotients. We have $G = P^{-1}P$. We claim that $P \subseteq G$  is Toeplitz in this case: Take $g \in G$, and write $g = q^{-1}p$ for some $p, q \in P$. Then the partial bijection
$$
  g^{-1}P \cap P \to P \cap gP, \, x \ma gx
$$
is the composition of
$$
  q^{-1}: \: qP \to P, \, qx \ma x \ {\rm and} \ p: \: P \to pP, \, x \ma px.
$$
This is because
$$
  g^{-1}P \cap P = p^{-1}qP \cap P = p^{-1}(qP) \cap P = p^{-1}(\dom(q^{-1})) = \dom(q^{-1}p),
$$
and for $x \in g^{-1}P \cap P = \dom(q^{-1}p)$, we have $gx = q^{-1}px = (q^{-1}p)(x)$.

In particular, if $P$ is the positive cone in a totally ordered group $G$, then $P \subseteq G$ is Toeplitz. Also, the inclusion $B_n^+ \subseteq B_n$ of the Braid monoid into the corresponding Braid group is Toeplitz. Furthermore, if $R$ is an integral domain with quotient field $Q$, then for the $ax+b$-semigroup $R \rtimes R\reg$, we have that $R \rtimes R\reg \subseteq Q \rtimes Q\reg$ is Toeplitz.

Let us discuss a second class of examples. Suppose that we have a monoid $P$ with identity $e$, and that $P \subseteq G$ is an embedding of $P$ into a group $G$. Furthermore, we assume that
$$\cJ_{P \subseteq G}\reg = \menge{gP}{g \in G}.$$
In this situation, we claim that $P \subseteq G$ is Toeplitz.

To see this, take $g \in G$. If $g^{-1}P \cap P \neq \emptyset$, then we can find $p \in P$ such that $g^{-1}P \cap P = pP$. This is because we have $\cJ_{P \subseteq G}\reg = \menge{gP}{g \in G}$ by assumption. Here, we used the hypothesis that $P$ has an identity element. Therefore, we can find $q \in P$ with $g^{-1}q = p$. We now claim that the partial bijection
$$
  g^{-1}P \cap P \to P \cap gP, \, x \ma gx
$$
is the composition of
$$
  p: \: P \to pP, \, x \ma px \ {\rm and} \ q^{-1}: \: qP \to P, \, qx \ma x.
$$
This is because
$$
  g^{-1}P \cap P = qP = \dom(pq^{-1}),
$$
and for $x \in g^{-1}P \cap P = \dom(pq^{-1})$, we have $gx = pq^{-1}x = (pq^{-1})(x)$.

In particular, for every graph $\Gamma$ as in \S~\ref{ex:presentations}, the inclusion $A_{\Gamma}^+ \subseteq A_{\Gamma}$ of the right-angled Artin monoid in the corresponding right-angled Artin group is Toeplitz. For instance, the canonical embedding $\Nz * \Nz \into \Fz_2$ is Toeplitz. Also, the inclusion $B_{k,l}^+ \subseteq B_{k,l}$ of the Baumslag-Solitar monoid into the corresponding Baumslag-Solitar group is Toeplitz, for $k, l \geq 1$. Moreover, the inclusion $F^+ \subseteq F$ of the Thompson monoid into the Thompson group is Toeplitz.

We make the following observation, which is an immediate consequence of our preceding discussion and Lemma~\ref{Lem:J_PinG}:
\bremark
\label{J=pP->J_PinG=gP}
Suppose that $P$ is a monoid which is embedded into a group $G$. If
$$\cJ_P\reg = \menge{pP}{p \in P},$$
then $P \subseteq G$ is Toeplitz if and only if
$$\cJ_{P \subseteq G}\reg = \menge{gP}{g \in G}.$$
\eremark

Let us present two examples of semigroup embeddings into groups which are not Toeplitz. In both of our examples, the semigroup will be given by the non-abelian free monoid $\Nz * \Nz$ on two generators.

First, consider the canonical homomorphism $\Nz * \Nz \to \Fz_2 / \Fz_2''$. Here, $\Fz_2''$ is the second commutator subgroup of the non-abelian free group $\Fz_2$ on two generators. By \cite{Hoch}, this canonical homomorphism $\Nz * \Nz \to \Fz_2 / \Fz_2''$ is injective. We want to see that $\Nz * \Nz \into \Fz_2 / \Fz_2''$ is not Toeplitz.

Let us denote both the canonical generators of $\Nz * \Nz$ and $\Fz_2$ by $a$ and $b$. We use the notation $[g,h] = ghg^{-1}h^{-1}$ for commutators. Obviously,
$$
  [(ab)^{-1},(ba)^{-1}] [ba,bab] [(ab)^{-1},(ba)^{-1}]^{-1} [ba,bab]^{-1} 
$$
lies in $\Fz_2''$. Thus
$$
  (ba) (ab) [(ab)^{-1},(ba)^{-1}] [ba,bab] [(ab)^{-1},(ba)^{-1}]^{-1} [ba,bab]^{-1} (ab)^{-1} (ba)^{-1} 
$$
lies in $\Fz_2''$. Now set
\bglnoz
  p &=& (ab) (ba) (ba) (bab)\\
  q &=& (ab) (ba) (bab) (ba)\\
  x &=& (ba) (ab) (bab) (ba)\\
  y &=& (ba) (ab) (ba) (bab).  
\eglnoz
Then
\bglnoz
  && p q^{-1} y x^{-1}\\
  &=& (ba) (ab) [(ab)^{-1},(ba)^{-1}] [ba,bab] [(ab)^{-1},(ba)^{-1}]^{-1} [ba,bab]^{-1} (ab)^{-1} (ba)^{-1} 
\eglnoz
lies in $\Fz_2''$. Therefore, we have $pq^{-1} = xy^{-1}$ in $\Fz_2 / \Fz_2''$. Now we consider $g = pq^{-1}$. Obviously, $P \cap gP \neq \emptyset$ as $p \in gP$. Moreover, we know that for $P = \Nz * \Nz$, the non-empty constructible right ideals are given by $\cJ_P\reg = \menge{pP}{p \in P}$. Hence by Remark~\ref{J=pP->J_PinG=gP}, if $\Nz * \Nz \into \Fz_2 / \Fz_2''$ were Toeplitz, we would have $P \cap gP = zP$ for some $z \in P$, as $P \cap gP$ must lie in $\cJ_P\reg$.

We already know that $p$ lies in $P \cap gP$. Moreover, $x$ lies in $P \cap gP$ as $x=gy$ in $\Fz_2 / \Fz_2''$. But the only element $z \in P$ with $p \in zP$ and $x \in zP$ is the identity element $z = e$. This is because $p$ starts with $a$ while $x$ starts with $b$.

Hence, if $\Nz * \Nz \into \Fz_2 / \Fz_2''$ were Toeplitz, we would have $P \cap gP = P$, or in other words, $P \subseteq gP$. In particular, the identity element $e \in P$ must be of the form $e = pq^{-1}r$ for some $r \in P$. Hence it would follow that $q = rp$ in $\Fz_2 / \Fz_2''$, and therefore in $\Nz * \Nz$. But this is absurd as $p \neq q$ while $p$ and $q$ have the same word length with respect to the generators $a$ and $b$.

All in all, this shows that $\Nz * \Nz \into \Fz_2 / \Fz_2''$ is not Toeplitz.

Our second example is given as follows: Again, we take $P = \Nz * \Nz$. But this time, we let our group be the Thompson group
$$
  F \defeq \spkl{x_0, x_1, \dotsc \ \vert \ x_n x_k = x_k x_{n+1} \ {\rm for} \ k < n}.
$$
Let $a$ and $b$ be the canonical free generators of $\Nz * \Nz$. Consider the homomorphism
$$
  \Nz * \Nz \to F, \, a \ma x_0, \, b \ma x_1.
$$
This is an embedding. For instance, this follows from uniqueness of the normal form in \cite[(1.3) in \S~1]{BrGe}. We claim that this embedding $\Nz * \Nz \into F$ is not Toeplitz.

To simplify notations, let us identify $\Nz * \Nz$ with the monoid $\spkl{x_0,x_1}^+$ generated by $x_0$ and $x_1$ in $F$. Consider
$$q = x_0^4 x_1 \ {\rm and} \ p = x_0^3.$$
Set $g \defeq pq^{-1}$. Then we have $p \in P \cap gP$. But we also have that $x_0 x_1 x_0^2$ lies in $P \cap gP$ because
$$
  x_0^3 x_1 x_0 x_1 = x_0^4 x_2 x_1 = x_0^4 x_1 x_3 \ {\rm in} \ F,
$$
so that
$$
  pq^{-1} x_0^3 x_1 x_0 x_1 = p q^{-1} x_0^4 x_1 x_3 = p x_3 = x_0^3 x_3 = x_0 x_1 x_0^2 \ {\rm in} \ F.
$$
If $\Nz * \Nz \into F$ were Toeplitz, we would have that $P \cap gP$ is of the form $zP$ for some $z \in P$. The argument is the same as in the previous example. But as we saw that $x_0^3$ and $x_0 x_1 x_0^2$ both lie in $P \cap gP$, our element $z$ can only be either the identity element $e$ or the generator $x_0$.

If $z = e$, then we would have $P \cap gP = P$, hence the identity $e$ must lie in $gP$. This means that there exists $r \in P$ with $e = gr = pq^{-1}r$ and therefore $q = rp$. But this is absurd.

If $z = x_0$, then we would have $P \cap gP = x_0P$, hence $x_0 \in gP$. Thus there must exist an element $r \in P$ with $x_0 = gr = pq^{-1}r$, and thus $qp^{-1}x_0 = r$. We conclude that
$$
  r = qp^{-1}x_0 = x_0^4 x_1 x_0^{-3} x_0 = x_0^4 x_1 x_0^{-2}
$$
so that
$$
  x_0^4 x_1 = r x_0^2.
$$
But this is again absurd.

All in all, this shows that $\Nz * \Nz \into F$ is not Toeplitz.

Looking at the preceding two examples, and comparing with our observation above that the canonical embedding $\Nz * \Nz \into \Fz_2$ is Toeplitz, we get the feeling that it is easier for the universal group embedding of a semigroup to satisfy the Toeplitz condition than for any other group embedding. Indeed, this is true. Let us explain the reason. We need the following equivalent formulation of the Toeplitz condition:
\blemma
\label{equi-from_T}
Let $P$ be a semigroup, and suppose that $P \subseteq G$ is an embedding of $P$ into a group $G$. The inclusion $P \subseteq G$ satisfies the Toeplitz condition if and only if for all $p, q \in P$, there exists a partial bijection $s \in I_l(P)$ with $s(q) = p$ and the intersection $P \cap qp^{-1}P$, taken in $G$, is contained in the domain $\dom(s)$.
\elemma
\bproof
If $g \in G$ satisfies $g^{-1}P \cap P \neq \emptyset$, then there exists $p, q$ in $P$ with $g^{-1}p = q$, i.e., $g = pq^{-1}$. This shows that $P \subseteq G$ is Toeplitz if and only if for all $p, q \in P$, the partial bijection
$$
  qp^{-1}P \cap P \to P \cap pq^{-1}P, \, x \ma pq^{-1}x
$$
lies in $I_l(P)$. But this is precisely what our condition says.
\eproof

\bcor
Suppose that we have a semigroup $P$ with two group embeddings $P \into G$ and $P \into \ti{G}$. Furthermore, assume that there is a group homomorphism $\ti{G} \to G$ such that the diagram
\bgl
\label{PtiGG}
  \xymatrix@C=20mm{
  P \ar@{^{(}->}[r] \ar[dr] & \ti{G} \ar[d] \\
   & G
  }
\egl
commutes. Then if $P \into G$ is Toeplitz, then the inclusion $P \into \ti{G}$ must be Toeplitz as well.
\ecor
\bproof
In our equivalent formulation of the Toeplitz condition (see Lemma~\ref{equi-from_T}), the only part which depends on the group embedding of our semigroup is the intersection $P \cap qp^{-1}P$. In our particular situation, the intersection $P \cap qp^{-1}P$ taken in $\ti{G}$ is given by
$$
  \menge{x \in P}{pq^{-1}x \in P \ {\rm in} \ \ti{G}},
$$
while the intersection $P \cap qp^{-1}P$ taken in $G$ is given by
$$
  \menge{x \in P}{pq^{-1}x \in P \ {\rm in} \ G}.
$$
Because of the commutative diagram \eqref{PtiGG}, the condition $pq^{-1}x \in P$ in $\ti{G}$ implies the condition $pq^{-1}x \in P$ in $G$. Hence the intersection $P \cap qp^{-1}P$, taken in $\ti{G}$, is contained in the intersection $P \cap qp^{-1}P$, taken in $G$, where we view both intersections as subsets of $P$. Our claim follows.
\eproof

As an immediate consequence, we obtain
\bcor
Let $P$ be a semigroup which embeds into a group, and assume that $P \into G_{\rm univ}$ is its universal group embedding. If $P \into G_{\rm univ}$ does not satisfy the Toeplitz condition, then for any other embedding $P \into G$ of our semigroup into a group $G$, we must have that $P \into G$ does not satisfy the Toeplitz condition either.
\ecor

\section{Graph products}
\label{sec:GraphProducts}

We discuss the independence condition and the Toeplitz condition for graph products.

Let $\Gamma = (V,E)$ be a graph with vertices $V$ and edges $E$. Assume that two vertices in $V$ are connected by at most one edge, and no vertex is connected to itself. Hence we view $E$ as a subset of $V \times V$. For every $v \in V$, let $P_v$ be a submonoid of a group $G_v$. We then form the graph products
$$P \defeq \Gamma_{v \in V} P_v$$
and
$$G \defeq \Gamma_{v \in V} G_v,$$
as in \S~\ref{ss:GraphProducts}. As explained in \S~\ref{ss:GraphProducts}, we can think of $P$ as a submonoid of $G$ in a canonical way.

Our goal is to prove that if each of the individual semigroups $P_v$, for all $v \in V$, satisfy the independence condition, then the graph product $P$ also satisfies the independence condition. Similarly, if each of the pairs $P_v \subseteq G_v$, for all $v \in V$, are Toeplitz, then the pair $P \subseteq G$ satisfies the Toeplitz condition as well. Along the way, we give an explicit description for the constructible right ideals of $P$.

We use the same notation as in \S~\ref{ss:GraphProducts}.

\subsection{Constructible right ideals}

Let us start with some easy observations.

\blemma
\label{omit-reduced}
Let $x_1 \dotso x_s$ be a reduced expression for $x \in G$, with $x_i \in G_{v_i}$. Assume that $v_1, \dotsc, v_j \in V^i(x)$. Then for all $1 \leq i \leq j$, $x_1 \dotso x_{i-1} x_{i+1} \dotso x_s$ is a reduced expression (for $x_i^{-1} x$). Similarly, if $v_{s-j}, \dotsc, v_s \in V^f(x)$, then for all $1 \leq i \leq j$, $x_1 \dotso x_{s-i-1} x_{s-i+1} \dotso x_s$ is a reduced expression (for $x x_{s-i}^{-1}$).
\elemma
\bproof
By assumption, the expressions $x_1 \dotso x_s$ and $x_i x_1 \dotso x_{i-1} x_{i+1} \dotso x_s$ are shuffle equivalent. In particular, the latter expression is reduced. Our first claim follows. The second assertion is proven analogously.
\eproof

\blemma
\label{initial-gx}
For $w \in V$, let $g$ be an element in $G_w$. Then for every $x \in G$, we have $g S_w^i(x) = S_w^i(gx)$.
\elemma
\bproof
Let $x_1 \dotso x_s$ be a reduced expression for $x$. If $w \notin V^i(x)$, then Lemma~\ref{XYZ} implies that $g x_1 \dotso x_s$ is a reduced expression for $gx$, and our claim follows. If $w \in V^i(x)$, we may assume that $x_1 = S_w^i(x)$. If $g x_1 \neq e$, then obviously $(g x_1) x_2 \dotso x_s$ is a reduced expression for $gx$, and we are done. If $g x_1 = e$, then $x_2 \dotso x_s$ is a reduced expression for $gx$ by Lemma~\ref{omit-reduced}. Clearly, $w \notin V^i(gx)$, and our claim follows.
\eproof

\bdefin
\label{product-ideals}
Let $W \subseteq V$ be a subset with $W \times W \subseteq E$, i.e., for every $w_1$, $w_2$ in $W$, we have $(w_1,w_2) \in E$. Given constructible right ideals $X_w \in \cJ_{P_w}$ for every $w \in W$, we set
$$
  \rukl{\prod_{w \in W} X_w} \cdot P \defeq \menge{x \in P}{S_w^i(x) \in X_w \fa w \in W}.
$$
If for some $w \in W$, we have $X_w = \emptyset$, then we set $\rukl{\prod_{w \in W} X_w} \cdot P = \emptyset$. If $W = \emptyset$, we set $\rukl{\prod_{w \in W} X_w} \cdot P = P$.
\edefin

By construction, we clearly have
$$
  \rukl{\prod_{w \in W} X_w} \cdot P = \bigcap_{w \in W} (X_w \cdot P).
$$

\blemma
\label{const}
Assume that $X_w = p_1^{-1} q_1 \dotso p_n^{-1} q_n (P_w)$ for some $p_i, q_i \in P_w$. Then we have $X_w \cdot P = p_1^{-1} q_1 \dotso p_n^{-1} q_n (P)$. Here we view $p_i$, $q_i$ as elements of $P$ (via the canonical embedding $P_w \subseteq P$).
\elemma
\bproof
We proceed inductively on $n$. The case $n=0$ is trivial. Let $p_i$, $q_i$ be elements of $P_w$, for $1 \leq i \leq n+1$. Set $Y_w \defeq p_2^{-1} q_2 \dotso p_{n+1}^{-1} q_{n+1} (P_w)$. We compute
\bglnoz
  (q_1 (Y_w)) \cdot P &=& \menge{x \in P}{S_w^i(x) \in q_1 (Y_w)} \\
  &=& \menge{x \in q_1 P}{q_1^{-1} (x) \in Y_w \cdot P} \\
  &=& \menge{x \in q_1 P}{x \in q_1(Y_w \cdot P)} \\
  &=& \menge{x \in P}{x \in q_1 p_2^{-1} q_2 \dotso p_{n+1}^{-1} q_{n+1} (P)}.
\eglnoz
Finally,
\bglnoz
  (p_1^{-1} q_1 (Y_w)) \cdot P &=& \menge{x \in P}{S_w^i(x) \in p_1^{-1} q_1 (Y_w)} \\
  &=& \menge{x \in P}{p_1 S_w^i(x) \in q_1 (Y_w)} \\
  &=& \menge{x \in P}{S_w^i(p_1 x) \in q_1 (Y_w)} \ \ \ \text{ by Lemma~\ref{initial-gx}} \\
  &=& \menge{x \in P}{p_1 x \in (q_1 (Y_w)) \cdot P} \\
  &=& p_1^{-1} (q_1 (Y_w)) \cdot P = p_1^{-1} q_1 \dotso p_{n+1}^{-1} q_{n+1} (P).
\eglnoz
\eproof

\blemma
\label{standard-lem}
Assume that we are given $p \in P$ and $W$, $\menge{X_w}{w \in W}$ as in Definition~\ref{product-ideals}. Assume that $\emptyset \neq p \rukl{\prod_{w \in W} X_w} \cdot P \neq P$. Then there exist $\ti{p}$ in $P$, $\ti{W} \subseteq V$ with $\ti{W} \times \ti{W} \subseteq E$, $\ti{X}_w \in \cJ_{P_w}$ for $w \in \ti{W}$ with
\begin{itemize}
\item $\ti{W} \neq \emptyset$ and $\emptyset \neq \ti{X}_w \neq P_w$ for every $w \in \ti{W}$,
\item either $\ti{p} = e$ or for all $v \in V^f(\ti{p})$, there exists $w \in \ti{W}$ with $(v,w) \notin E$,
\end{itemize}
such that
$$
  p \rukl{\prod_{w \in W} X_w} \cdot P = \ti{p} \rukl{\prod_{w \in \ti{W}} \ti{X}_w} \cdot P.
$$
\elemma
\bproof
We proceed inductively on the length $l(p)$ of $p$. If $l(p) = 0$, i.e., $p = e$, then for all $w \in W$, we must have $X_w \neq \emptyset$, and there must exist $w \in W$ with $X_w \neq P_w$. Thus, we can set
$$\ti{W} \defeq \menge{w \in W}{X_w \neq P_w} \ {\rm and} \ \ti{X}_w \defeq X_w \ {\rm for} \ w \in \ti{W}.$$

Now assume that $l(p)>0$. Without changing $p \rukl{\prod_{w \in W} X_w} \cdot P$, we can replace $W$ by $\menge{w \in W}{X_w \neq P_w}$. So we may just as well assume that for every $w \in W$, we have $\emptyset \neq X_w \neq P_w$. If for every $v \in V^f(p)$, there exists $w \in W$ with $(v,w) \notin E$, then we can just set $\ti{W} = W$ and $\ti{X}_w = X_w$ for all $w \in \ti{W}$. If not, then we choose $v \in V^f(p)$ with $(v,w) \in E$ for every $w \in W$. Let $p_1 \dotso p_r$ be a reduced expression for $p$, with $p_r \in P_v$. Set $X_v \defeq P_v$ if $v \notin W$. Using Lemma~\ref{omit-reduced} and Lemma~\ref{initial-gx}, we deduce
\bglnoz
  && p_r \rukl{\prod_{w \in W} X_w} \cdot P = \menge{y \in P}{y = p_r x \text{ for some } x \in \rukl{\prod_{w \in W} X_w} \cdot P} \\
  &=& \menge{y \in P}{S_w^i(y) \in X_w \fa v \neq w \in W \text{ and } S_v^i(y) \in p_r X_v} \\
  &=& (p_r X_v) \cdot \rukl{\prod_{v \neq w \in W} X_w} \cdot P.
\eglnoz
Thus
$$p \rukl{\prod_{w \in W} X_w} \cdot P = (p_1 \dotsm p_{r-1}) \eckl{(p_r X_v) \cdot \rukl{\prod_{v \neq w \in W} X_w}} \cdot P.$$
Now our claim follows once we apply the induction hypothesis with $p_1 \dotsm p_{r-1}$ in place of $p$.
\eproof

\bdefin
\label{standard-def}
Assume that we are in the situation of Lemma~\ref{standard-lem}, i.e., we are given $p \in P$ and $W$, $\menge{X_w}{w \in W}$ as in Definition~\ref{product-ideals}. Assume that
$$\emptyset \neq p \rukl{\prod_{w \in W} X_w} \cdot P \neq P.$$
Then $$p \rukl{\prod_{w \in W} X_w} \cdot P$$ is said to be in standard form if both conditions from the Lemma~\ref{standard-lem} are satisfied, i.e.,
\begin{itemize}
\item $W \neq \emptyset$ and $\emptyset \neq X_w \neq P_w$ for every $w \in W$,
\item either $p = e$ or for all $v \in V^f(p)$, there exists $w \in W$ with $(v,w) \notin E$.
\end{itemize}
\edefin

\blemma
\label{standard-->reduced}
Assume that $p \rukl{\prod_{w \in W} X_w} \cdot P$ is in standard form. Given reduced expressions $p_1 \dotsm p_r$ for $p$ and $x_1 \dotsm x_s$ for $x \in \rukl{\prod_{w \in W} X_w} \cdot P$, $p_1 \dotsm p_r x_1 \dotsm x_s$ is a reduced expression for $px$. In particular, if in addition $p \neq e$, then for every $v \in V^i(p)$, we have $S_v^i(px) = S_v^i(p)$.
\elemma
\bproof
For our first claim, the case $p=e$ is trivial. So let us assume $p \neq e$. As $X_w \neq P_w$ for all $w \in W$, we know that $W \subseteq V^i(x)$, so that $V^f(p) \cap V^i(x) = \emptyset$ because $p \rukl{\prod_{w \in W} X_w} \cdot P$ is in standard form. Then our assertion that $p_1 \dotsm p_r x_1 \dotsm x_s$ is reduced follows from Lemma~\ref{XYZ}.
\eproof

\bprop
\label{const-ideals}
The non-empty constructible right ideals of $P$ are precisely given by all the non-empty subsets of $P$ of the form $p \rukl{\prod_{w \in W} X_w} \cdot P$, with $p \in P$ and $W$, $\menge{X_w}{w \in W}$ as in Definition~\ref{product-ideals}.
\eprop
\bproof
First, we prove that $p \rukl{\prod_{w \in W} X_w} \cdot P$ is constructible. It certainly suffices to check that $\rukl{\prod_{w \in W} X_w} \cdot P$ is constructible. But Lemma~\ref{const} tells us that $X_w \cdot P$ is constructible for every $w \in W$. Therefore, $\rukl{\prod_{w \in W} X_w} \cdot P = \bigcap_{w \in W} (X_w \cdot P)$ is constructible itself.

Secondly, we show that every non-empty constructible right ideal is of the form $p \rukl{\prod_{w \in W} X_w} \cdot P$. For this purpose, let $\cJ'$ be the set of all non-empty constructible right ideals which are of the form $p \rukl{\prod_{w \in W} X_w} \cdot P$. Clearly, $P$ lies in $\cJ'$. Also, if $\emptyset \neq X \in \cJ'$ and $p \in P$, then obviously $pX$ lies in $\cJ'$. It remains to prove that for $\emptyset \neq X \in \cJ'$ and $q \in P$, we have $q^{-1} (X) \in \cJ'$ if $q^{-1}(X) \neq \emptyset$. Since the set $\cJ_P\reg$ of non-empty constructible right ideals of $P$ is minimal with respect to these properties, this would then show that $\cJ_P\reg \subseteq \cJ'$, as desired. By induction on $l(q)$, we may assume that $q \in P_v$, and it even suffices to consider the case $q \in P_v \setminus P_v^*$. For $X = p \rukl{\prod_{w \in W} X_w} \cdot P$, we want to show that $q^{-1}(X) = \emptyset$ or $q^{-1}(X) \in \cJ'$. We distinguish between the following cases:

1.) $p=e$:

1.a) There exists $w \in W$ with $(v,w) \notin E$. Without loss of generality we may assume that $X_w \neq P_w$ for all $w \in W$. Then for every $x \in P$, $w \notin V^i(qx)$ since $v \in V^i(qx)$. Thus $S_w^i(qx) = e \notin X_w$. Therefore, 
$$q^{-1} \rukl{\prod_{w \in W} X_w} \cdot P = \emptyset.$$

1.b) We have $(v,w) \in E$ for all $w \in W$ and $v \notin W$. Then
\bglnoz
  q^{-1} \rukl{\prod_{w \in W} X_w} \cdot P &=& \menge{x \in P}{S_w^i(qx) \in X_w \fa w \in W} \\
  &=& \menge{x \in P}{S_w^i(x) \in X_w \fa w \in W} \\
  &=& \rukl{\prod_{w \in W} X_w} \cdot P \in \cJ'.
\eglnoz

1.c) We have $(v,w) \in E$ for all $w \in W$ and $v \in W$. Then
\bglnoz
  q^{-1} \rukl{\prod_{w \in W} X_w} \cdot P &=& \menge{x \in P}{S_w^i(qx) \in X_w \fa w \in W} \\
  &=& \eckl{(q^{-1} (X_v)) \cdot \rukl{\prod_{v \neq w \in W} X_w}} \cdot P \in \cJ'.
\eglnoz

2.) $p \neq e$: We can clearly assume that
$$\emptyset \neq p \rukl{\prod_{w \in W} X_w} \cdot P \neq P.$$
By Lemma~\ref{standard-lem}, we may assume that $p \rukl{\prod_{w \in W} X_w} \cdot P$ is in standard form. And because we have already finished the case $p=e$, we can in addition assume that $p \neq e$. Without loss of generality, we may assume $v \in V^i(p)$, as we would otherwise have $q^{-1} \left[ p \left( \prod_{w \in W} X_w \right) \cdot P \right] = \emptyset$ or $q^{-1} \left[ p \left( \prod_{w \in W} X_w \right) \cdot P \right] = p \left( \prod_{w \in W} X_w \right) \cdot P$. Lemma~\ref{standard-->reduced} gives $S_v^i(px) = S_v^i(p)$ for every $x \in \rukl{\prod_{w \in W} X_w} \cdot P$. Now $y$ lies in $q^{-1} \eckl{p \rukl{\prod_{w \in W} X_w} \cdot P}$ if and only if there exists $x \in \rukl{\prod_{w \in W} X_w} \cdot P$ such that $qy = px$. Hence if there exists $y \in q^{-1} \eckl{p \rukl{\prod_{w \in W} X_w} \cdot P}$, we must have
$$q S_v^i(y) = S_v^i(qy) = S_v^i(px) = S_v^i(p)$$
by Lemma~\ref{initial-gx}. Thus $p \in S_v^i(p) P \subseteq qP$. This implies that
$$
  q^{-1} \eckl{p \rukl{\prod_{w \in W} X_w} \cdot P} = (q^{-1}p) \rukl{\prod_{w \in W} X_w} \cdot P \in \cJ'.
$$
\eproof

\subsection{The independence condition}

\blemma
\label{p-in-p}
Assume that
$$\emptyset \neq p \rukl{\prod_{w \in W} X_w} \cdot P \neq P \ {\rm and} \ \emptyset \neq \ti{p} \rukl{\prod_{w \in \ti{W}} \ti{X}_w} \cdot P \neq P$$
are in standard form, with $p \neq e$. If
$$\ti{p} \rukl{\prod_{w \in \ti{W}} \ti{X}_w} \cdot P \subseteq p \rukl{\prod_{w \in W} X_w} \cdot P,$$
then $\ti{p} \in pP$.
\elemma
\bproof
First of all, let us show that $\ti{p} \neq e$. Namely, assume the contrary, i.e., $\ti{p} = e$. Take $\ti{x} \in \rukl{\prod_{w \in \ti{W}} \ti{X}_w} \cdot P$. By assumption, we can find $x \in  \rukl{\prod_{w \in W} X_w} \cdot P$ so that $\ti{x} = px$. Moreover, choose $v \in V^i(p)$. By Lemma~\ref{standard-->reduced}, it follows that $S_v^i(\ti{x}) = S_v^i(px) = S_v^i(p)$. Thus we have proven that every $\ti{x} \in \rukl{\prod_{w \in \ti{W}} \ti{X}_w} \cdot P$ must satisfy $S_v^i(\ti{x}) = S_v^i(p)$. But this is obviously a wrong statement. Thus we must have $\ti{p} \neq e$.

Now we proceed inductively on $l(p)$. We start with the case $l(p) = 1$, i.e., $p \in P_v$. For $\ti{x} \in \rukl{\prod_{w \in \ti{W}} \ti{X}_w} \cdot P$ with $S_v^i(\ti{p} \ti{x}) = S_v^i(\ti{p})$, we can always find $x \in \rukl{\prod_{w \in W} X_w} \cdot P$ so that $\ti{p} \ti{x} = px$. By Lemma~\ref{standard-->reduced}, we deduce that $p = S_v^i(px) = S_v^i(\ti{p} \ti{x}) = S_v^i(\ti{p})$. Therefore, $\ti{p} \in S_v^i(\ti{p}) P = pP$.

For the induction step, take $v \in V^i(p)$. For $\ti{x} \in \rukl{\prod_{w \in \ti{W}} \ti{X}_w} \cdot P$ with $S_v^i(\ti{p} \ti{x}) = S_v^i(\ti{p})$, again choose $x \in \rukl{\prod_{w \in W} X_w} \cdot P$ so that $\ti{p} \ti{x} = px$. Then $S_v^i(p) = S_v^i(px) = S_v^i(\ti{p} \ti{x}) = S_v^i(\ti{p})$. This shows that both $\ti{p}$ and $p$ lie in $S_v^i(p)P$. We deduce that
$$(S_v^i(p)^{-1} \ti{p}) \rukl{\prod_{w \in \ti{W}} \ti{X}_w} \cdot P \subseteq (S_v^i(p)^{-1} p) \rukl{\prod_{w \in W} X_w} \cdot P.$$
Since $l(S_v^i(p)^{-1} p) < l(p)$, we can now apply the induction hypothesis, and we are done.
\eproof

\blemma
\label{p-in-e}
As above, let
$$\emptyset \neq \rukl{\prod_{w \in W} X_w} \cdot P \neq P \ {\rm and} \ \emptyset \neq \ti{p} \rukl{\prod_{w \in \ti{W}} \ti{X}_w} \cdot P \neq P$$
be in standard form, this time with $\ti{p} \neq e$ (and $p=e$). If
$$\ti{p} \rukl{\prod_{w \in \ti{W}} \ti{X}_w} \cdot P \subseteq \rukl{\prod_{w \in W} X_w} \cdot P,$$
then 
$$\ti{p} \in \rukl{\prod_{w \in W} X_w} \cdot P.$$
\elemma
\bproof
For $\ti{x} \in \rukl{\prod_{w \in \ti{W}} \ti{X}_w} \cdot P$ with $S_v^i(\ti{p} \ti{x}) = S_v^i(\ti{p})$, $\ti{p} \ti{x}$ lies in $\rukl{\prod_{w \in W} X_w} \cdot P$ by assumption. Hence, Lemma~\ref{standard-->reduced} tells us that for all $w \in W$, $S_w^i(\ti{p}) = S_w^i(\ti{p} \ti{x})$ lies in $X_w$. Thus $\ti{p}$ lies in $\rukl{\prod_{w \in W} X_w} \cdot P$.
\eproof

\blemma
\label{e-in-e}
Let
$$\emptyset \neq \rukl{\prod_{w \in W} X_w} \cdot P \neq P \ {\rm and} \ \emptyset \neq \rukl{\prod_{w \in \ti{W}} \ti{X}_w} \cdot P \neq P$$
be in standard form. Then
$$\rukl{\prod_{w \in \ti{W}} \ti{X}_w} \cdot P \subseteq \rukl{\prod_{w \in W} X_w} \cdot P$$
if and only if $W \subseteq \ti{W}$ and $\ti{X}_w \subseteq X_w$ for every $w \in W$.
\elemma
\bproof
The direction \an{$\Larr$} is obvious. To prove the reverse direction, first assume that $W \nsubseteq \ti{W}$. Choose for every $\ti{w} \in \ti{W}$ an element $x_{\ti{w}} \in \ti{X}_{\ti{w}}$. Then the product $\prod_{\ti{w} \in \ti{W}} x_{\ti{w}}$ obviously lies in $\rukl{\prod_{w \in \ti{W}} \ti{X}_w} \cdot P$. But for $w \in W \setminus \ti{W}$, we have $S_w^i(\prod_{\ti{w} \in \ti{W}} x_{\ti{w}}) = e \notin X_w$ as $X_w \neq P_w$. This contradicts $\rukl{\prod_{w \in \ti{W}} \ti{X}_w} \cdot P \subseteq \rukl{\prod_{w \in W} X_w} \cdot P$. So we must have $W \subseteq \ti{W}$. If for some $w \in W$, we have $\ti{X}_w \nsubseteq X_w$, then choose $x_w \in \ti{X}_w \setminus X_w$. For all remaining $\ti{w} \in \ti{W} \setminus \gekl{w}$, choose $x_{\ti{w}} \in \ti{X}_{\ti{w}}$. Then the product $\prod_{\ti{w} \in \ti{W}} x_{\ti{w}}$ lies in $\rukl{\prod_{w \in \ti{W}} \ti{X}_w} \cdot P$. But $S_w^i(\prod_{\ti{w} \in \ti{W}} x_{\ti{w}}) = x_w \notin X_w$. This again contradicts $\rukl{\prod_{w \in \ti{W}} \ti{X}_w} \cdot P \subseteq \rukl{\prod_{w \in W} X_w} \cdot P$.
\eproof

\bprop
\label{ind_GraphProducts}
If for every $v \in V$, the semigroup $P_v$ satisfies independence, then the graph product $P$ satisfies independence.
\eprop
\bproof
Let
$$\emptyset \neq p \rukl{\prod_{w \in W} X_w} \cdot P \neq P$$
be in standard form, and let
$$\emptyset \neq p_i \rukl{\prod_{w \in W_i} X^{(i)}_w} \cdot P \neq P$$
be finitely many constructible right ideals of $P$ in standard form. If
$$
  p \rukl{\prod_{w \in W} X_w} \cdot P = \bigcup_i p_i \rukl{\prod_{w \in W_i} X^{(i)}_w} \cdot P,
$$
then either $p=e$ or $p_i \in pP$ for all $i$ by Lemma~\ref{p-in-p}. Hence
$$
  \rukl{\prod_{w \in W} X_w} \cdot P = \bigcup_i (p^{-1} p_i) \rukl{\prod_{w \in W_i} X^{(i)}_w} \cdot P
$$
Therefore, we may without loss of generality assume that $p=e$, i.e. 
\bgl
\label{union}
  \rukl{\prod_{w \in W} X_w} \cdot P = \bigcup_i p_i \rukl{\prod_{w \in W_i} X^{(i)}_w} \cdot P.
\egl
Let $I = \menge{i}{p_i \neq e}$ and $J = \menge{i}{p_i = e}$. By Lemma~\ref{p-in-e}, we have for all $i \in I$ and $w \in W$ that $S_w^i(p_i) \in X_w$. We define for every $i \in I$:
$$p_i' = \prod_{w \in W} S_w^i(p_i).$$
For each $i \in I$, we obviously have
$$p_i \rukl{\prod_{w \in W_i} X^{(i)}_w} \cdot P \subseteq p_i P \subseteq p_i' P \subseteq \rukl{\prod_{w \in W} X_w} \cdot P.$$
Therefore,
$$
  \rukl{\prod_{w \in W} X_w} \cdot P = \bigcup_{i \in I} (p_i' P) \cup \bigcup_{i \in J} \rukl{\prod_{w \in W_i} X^{(i)}_w} \cdot P.
$$
Set $\ti{W}_i \defeq W$ if $i \in I$, $\ti{W}_i \defeq W_i$ for $i \in J$ and
$$
  \ti{X}^{(i)}_w \defeq
  \bfa
    S_w^i(p_i) P_w & \falls i \in I, \ w \in \ti{W}_i, \\
    X^{(i)}_w & \falls i \in J, \ w \in \ti{W}_i.
  \efa
$$
Since $\rukl{\prod_{w \in \ti{W}_i} \ti{X}^{(i)}_w} \cdot P = p_i' P$ for all $i \in I$, we obviously again have
\bgl
\label{union'}
  \rukl{\prod_{w \in W} X_w} \cdot P = \bigcup_i \rukl{\prod_{w \in \ti{W}_i} \ti{X}^{(i)}_w} \cdot P.
\egl
Moreover, $\ti{X}^{(i)}_w \neq P_w$ for all $i$ and $w \in \ti{W}_i$.

By Lemma~\ref{e-in-e}, we must have $\ti{X}^{(i)}_w \subseteq X_w$ for all $i$ and $w \in \ti{W}_i$. Assume that for all $i$ with $\ti{W}_i = W$, there exists $w(i) \in W$ with $\ti{X}^{(i)}_{w(i)} \subsetneq X_{w(i)}$. Choose for every $w \in \gekl{w(i)}_i$ an element
$$x_w \in X_w \setminus \bigcup_{\menge{i}{w(i) = w}} \ti{X}^{(i)}_{w(i)}.$$
This is possible since $\cJ_{P_v}$ is independent for every $v \in V$, so that
$$X_w \setminus \bigcup_{\menge{i}{w(i) = w}} \ti{X}^{(i)}_{w(i)} \neq \emptyset.$$
For all remaining $w \in W$, just choose some $x_w \in X_w$. Then $x \defeq \prod_{w \in W} x_w$ lies in $\rukl{\prod_{w \in W} X_w} \cdot P$, but for all $i$ with $\ti{W}_i = W$, $S_{w(i)}^i(x)$ does not lie in $\ti{X}^{(i)}_{w(i)}$. Therefore, $x$ does not lie in $\rukl{\prod_{w \in \ti{W}_i} \ti{X}^{(i)}_w} \cdot P$ whenever $i$ satisfies $\ti{W}_i = W$. For $i$ with $\ti{W}_i \neq W$, take $\ti{w} \in \ti{W}_i \setminus W$. Then $S_{\ti{w}}^i(x) = e \notin \ti{X}^{(i)}_{\ti{w}}$. Thus also for $i$ with $\ti{W}_i \neq W$, we have $x \notin \rukl{\prod_{w \in \ti{W}_i} \ti{X}^{(i)}_w} \cdot P$. Since this contradicts \eqref{union'}, there must exist an index $i$ with $W_i = W$ and $\ti{X}^{(i)}_w = X_w$ for all $w \in W$. In particular, for that index $i$, we must have $$\rukl{\prod_{w \in W} X_w} \cdot P = \rukl{\prod_{w \in \ti{W}_i} \ti{X}^{(i)}_w} \cdot P.$$

If this index $i$ lies in $I$, then we have shown that $\rukl{\prod_{w \in W} X_w} \cdot P$ is a principal right ideal, and we are done. If this index $i$ lies in $J$, then we have proven that $\rukl{\prod_{w \in W} X_w} \cdot P$ coincides with one of the (constructible right) ideals on the right hand side of \eqref{union} (since $p_i = e$ for $i \in J$), and we are also done.
\eproof

\subsection{The Toeplitz condition}

\bdefin
Let $x \in G$, and assume that $x_1 \dotsm x_s$ is a reduced expression for $x$. We set $S(x) \defeq \gekl{x_1, \dotsc, x_s}$.
\edefin
Note that this is well-defined by Theorem~\ref{graphprod-wordproblem}.

\blemma
\label{syl}
Let $g, x \in G$, $v \in V^f(g)$, and assume $S_v^f(g) S_v^i(x) \neq e$. Then $S_v^f(g) S_v^i(x)$ lies in $S(gx)$.
\elemma
\bproof
Let $g_1 \dotsm g_r$ be a reduced expression for $g$, with $g_r = S_v^f(g)$.

First of all, if $V^f(g) \cap V^i(x) = \emptyset$, then Lemma~\ref{XYZ} tells us that for every reduced expression $x_1 \dotsm x_s$ for $x$, $g_1 \dotsm g_r x_1 \dotsm x_s$ is a reduced expression for $gx$. Hence $g_r = S_v^f(g) S_v^i(x)$ lies in $S(gx)$.

Secondly, assume that $V^f(g) \cap V^i(x) = \gekl{v}$. If $x_1 \dotsm x_s$ is a reduced expression for $x$ with $x_1 = S_v^i(x)$, then since $g_r x_1 \neq e$, Lemma~\ref{XYZ} tells us that $g_1 \dotsm g_{r-1} (g_r x_1) x_2 \dotsm x_s$ is a reduced expression for $gx$. Again, our claim follows.

Finally, it remains to treat the case
$$\emptyset \neq V^f(g) \cap V^i(x) \neq \gekl{v}.$$
We proceed inductively on $l(g)$. The cases $l(g)=0$ and $l(g)=1$ are taken care of by the previous cases. As
$$\emptyset \neq V^f(g) \cap V^i(x) \neq \gekl{v},$$
we can choose $w \in V^f(g) \cap V^i(x)$ with $w \neq v$. If $v$ lies in $V^f(g) \cap V^i(x)$, then choose a reduced expression $g_1 \dotsm g_r$ for $g$ with $g_{r-1} \in G_w$ and $g_r \in G_v$, and let $x_1 \dotsm x_s$ be a reduced expression for $x$ with $x_1 \in G_v$ and $x_2 \in G_w$. Then
$$gx = g_1 \dotsm g_{r-2} (g_r x_1) (g_{r-1} x_2) x_3 \dotsm x_s.$$
Set
$$g' \defeq g_1 \dotsm g_{r-2} g_r \ {\rm and} \ x' \defeq x_1 (g_{r-1} x_2) x_3 \dotsm x_s.$$
By Lemma~\ref{omit-reduced}, we know that $g_1 \dotsm g_{r-2} g_r$ is a reduced expression, so that $g_r = S_v^f(g')$. Also, $x_1 (g_{r-1} x_2) x_3 \dotsm x_s$ is a reduced expression. This is clear if $g_{r-2} x_2 \neq e$, and it follows from Lemma~\ref{omit-reduced} in case $g_{r-2} x_2 = e$. Thus $x_1 = S_v^i(x')$. So we again have $$S_v^f(g') S_v^i(x') = S_v^f(g) S_v^i(x) \neq e.$$
Since $l(g') < l(g)$, induction hypothesis tells us that $S_v^f(g) S_v^i(x) = S_v^f(g') S_v^i(x')$ lies in $S(g'x') = S(gx)$. The case $v \notin V^f(g) \cap V^i(x)$ is treated similarly. Just set $x_1 = e$.
\eproof

For $g \in G$, let us denote the partial bijection
$$
  g^{-1}P \cap P \to P \cap gP, \, x \ma gx
$$
by $g_P$.

\blemma
\label{factor}
Let $g_1 \dotsm g_r$ be a reduced expression for $g \in G$. Then
$$g_P = (g_1)_P \dotsm (g_r)_P.$$
\elemma
\bproof
We proceed inductively on $l(g)$. The case $l(g) = 1$ is trivial. First, we show that for $x \in P$, $gx \in P$ implies $g_r x \in P$. Let $g_r \in G_v$. Then by Lemma~\ref{syl}, $g_r S_v^i(x)$ lies in $S(gx)$ or $g_r S_v^i(x) = e$. Since $gx \in P$, we conclude that in any case, we have $g_r S_v^i(x) \in P_v$. Obviously, $S(g_r x) \subseteq \gekl{g_r S_v^i(x)} \cup S(x)$. So we obtain $g_r x \in P$. Therefore, we compute
\bglnoz
  \dom(g_P) &=& \menge{x \in P}{gx \in P} 
  = \menge{x \in P}{gx \in P \ {\rm and} \ g_r x \in P}\\
  &=& \menge{x \in P}{g_r x \in P \ {\rm and} \ (g_1 \dotso g_{r-1})(g_r x) \in P}\\
  &=& \dom( (g_1 \dotsm g_{r-1})_P \ g_P).
\eglnoz
Hence it follows that $g_P = (g_1 \dotsm g_{r-1})_P \ g_P$.

By induction hypothesis, $(g_1 \dotsm g_{r-1})_P = (g_1)_P \dotso (g_{r-1})_P$, and we are done.
\eproof

\blemma
\label{local-Toeplitz}
For $g \in G_v$, we have $g^{-1}P_v \cap P_v \neq \emptyset$ if and only if $g^{-1}P \cap P \neq \emptyset$.

Assume that this is the case, and that there are $p_i$, $q_i$ in $G_v$ with
$$g_{P_v} = p_1^{-1} q_1 \dotso p_n^{-1} q_n$$
in $I_l(P_v)$. Then
$$g_P = p_1^{-1} q_1 \dotso p_n^{-1} q_n$$
in $I_l(P)$.
\elemma
\bproof
Let us start proving the first claim. Since $P_v \subseteq P$, the implication \an{$\Rarr$} is obvious. For the reverse direction, assume that $g^{-1} P \cap P \neq \emptyset$, i.e., there exists $x \in P$ with $gx \in P$. Then obviously, $S_v^i(x) \in P_v$, and $g S_v^i(x) = S_v^i(gx)$ lies in $P_v$ (here we used Lemma~\ref{initial-gx}), so $g^{-1} P_v \cap P_v \neq \emptyset$.

Secondly, we show $g^{-1} P \cap P = q_n^{-1} p_n \dotso q_1^{-1} p_1 (P)$:
\bglnoz
  g^{-1} P \cap P &=& \menge{x \in P}{gx \in P} = \menge{x \in P}{S_v^i(gx) \in P_v} \\
  &=& \menge{x \in P}{g S_v^i(x) \in P_v} \ \text{ by Lemma~\ref{initial-gx}} \\
  &=& \menge{x \in P}{S_v^i(x) \in g^{-1} P_v \cap P_v} \\
  &=& \menge{x \in P}{S_v^i(x) \in q_n^{-1} p_n \dotso q_1^{-1} p_1 (P_v)}\\
  &=& \rukl{q_n^{-1} p_n \dotso q_1^{-1} p_1 (P_v)} \cdot P \\
  &=& q_n^{-1} p_n \dotso q_1^{-1} p_1 (P) \ \text{ by Lemma~\ref{const}}.
\eglnoz
Therefore, we have
$$
  \dom(g_P) = \dom(p_1^{-1} q_1 \dotsm p_n^{-1} q_n)
$$
as subsets of $P$. Hence it follows that
$$g_P = p_1^{-1} q_1 \dotso p_n^{-1} q_n$$
in $I_l(P)$ because we have $p_1^{-1} q_1 \dotsm p_n^{-1} q_n = g$ in $G_v \subseteq G$. Here we are taking products of $p_i^{-1}$ and $q_i$ as group elements in $G_v$ and $G$.
\eproof

\bprop
\label{Toeplitz_GraphProducts}
If for all $v \in V$, $P_v \subseteq G_v$ is Toeplitz, then $P \subseteq G$ is Toeplitz.
\eprop
\bproof
Let $g_1 \dotsm g_r$ be a reduced expression for $g \in G$, with $g_i \in G_{v_i}$. Assume that $g^{-1} P \cap P \neq \emptyset$. By Lemma~\ref{factor}, we know that
$$g_P = (g_1)_P \dotso (g_r)_P.$$
In particular, $g_i^{-1}P \cap P \neq \emptyset$ for all $1 \leq i \leq r$. By Lemma~\ref{local-Toeplitz}, we conclude that $g_i^{-1} P_{v_i} \cap P_{v_i} \neq \emptyset$ for all $1 \leq i \leq r$. Since for all $1 \leq i \leq r$, the embedding $P_{v_i} \subseteq G_{v_i}$ is Toeplitz, we can find $p_{i,j}$, $q_{i,j}$ in $P_{v_i}$ (for $1 \leq j \leq n_i$) with
$$(g_i)_{P_{v_i}} = p_{i,1}^{-1} q_{i,1} \dotso p_{i,n_i}^{-1} q_{i,n_i} \ {\rm in} \ I_l(P_{v_i}).$$
Lemma~\ref{local-Toeplitz} implies that 
$$
  (g_i)_P = p_{i,1}^{-1} q_{i,1} \dotso p_{i,n_i}^{-1} q_{i,n_i} \ {\rm in} \ I_l(P)
$$
for all $1 \leq i \leq r$. Thus we have, in $I_l(P)$:
$$
  g_P = (g_1)_P \dotso (g_r)_P = \rukl{p_{1,1}^{-1} q_{1,1} \dotso p_{1,n_1}^{-1} q_{1,n_1}} \dotso \rukl{p_{r,1}^{-1} q_{r,1} \dotso p_{r,n_r}^{-1} q_{r,n_r}}.
$$
\eproof

\section{K-theory}
\label{sec:K}

Let us apply the K-theory results from \cite{EchKK} to semigroups and their reduced semigroup C*-algebras.

Let $P$ be a semigroup which embeds into a group. Assume that $P$ satisfies independence, and that we have an embedding $P \subseteq G$ into a group $G$ such that $P \subseteq G$ is Toeplitz. Furthermore, suppose that $G$ satisfies the Baum-Connes conjecture with coefficients.

As $\cJ_{P \subseteq G}\reg = G. \cJ_P\reg$ by Lemma~\ref{Lem:J_PinG}, we can choose a set of representatives $\fX \subseteq \cJ_P$ for the $G$-orbits $G \backslash \cJ_{P \subseteq G}\reg$. For every $X \in \fX$, let
$$
  G_X \defeq \menge{g \in G}{gX = X},
$$
and let
$$
  \iota_X: \: C^*_{\lambda}(G_X) \to C^*_{\lambda}(P), \, \lambda_g \ma \lambda_g 1_X.
$$
Here we identify $C^*_{\lambda}(P)$ with the crossed product $D_{P \subseteq G} \rtimes_r G$ as in Proposition~\ref{Prop:Toeplitz->C*_r=fullcorner}. This is possible because of our assumption that $P \subseteq G$ is Toeplitz.
\btheo
\label{THM:K}
In the situation above, we have that
$$
  \bigoplus_{X \in \fX} (\iota_X)_*: \: \bigoplus_{X \in \fX} K_*(C^*_{\lambda}(G_X)) \overset{\cong}{\lori} K_*(C^*_{\lambda}(P))
$$
is an isomorphism.
\etheo
To see how Theorem~\ref{THM:K} follows from \cite[Corollary~5.19]{EchKK}, we explain how to choose $\Omega$, $I$ and $e_i$, $i \in I$ (in the notation of \cite[Corollary~5.19]{EchKK}). Let $\Omega$ be the spectrum of $D_{P \subseteq G}$ ($D_{P \subseteq G}$ was introduced in Definition~\ref{DEF:D}), so that our semigroup C*-algebra is a full corner in $C_0(\Omega) \rtimes_r G$ by Proposition~\ref{Prop:Toeplitz->C*_r=fullcorner}. Moreover, let $I$ be $\cJ_{P \subseteq G}\reg$, and let $e_X$ be given by $1_X$ for all $X \in \cJ_{P \subseteq G}\reg$. Applying \cite[Corollary~5.19]{EchKK}, with coefficient algebra $A = \Cz$, to this situation yields Theorem~\ref{THM:K}.

If, in addition, $P$ is a monoid and we have $\cJ_P\reg = \menge{pP}{p \in P}$, then we must have $\cJ_{P \subseteq G}\reg = \menge{gP}{g \in P}$, so that we may choose $\fX = \gekl{P}$. Then the stabilizer group $G_P = P^*$ becomes the group of units in $P$. The theorem above then says that the *-homomorphism
$$
  \iota: \: C^*_{\lambda}(P^*) \overset{\cong}{\lori} C^*_{\lambda}(P), \, \lambda_g \ma V_g
$$
induces an isomorphism
$$
  \iota_*: \: K_*(C^*_{\lambda}(P^*)) \overset{\cong}{\lori} K_*(C^*_{\lambda}(P)).
$$

In particular, if we further have that $P$ has trivial unit group, then we obtain that the unique unital *-homomorphism $\Cz \to C^*_{\lambda}(P)$ induces an isomorphism
$$
  K_*(\Cz) \overset{\cong}{\lori} K_*(C^*_{\lambda}(P)).
$$
This applies to positive cones in total ordered groups, as long as the group satisfies the Baum-Connes conjecture with coefficients. It also applies to right-angled Artin monoids, to Braid monoids, to Baumslag-Solitar monoids of the type $B_{k,l}^+$ for $k,l \geq 1$, and to the Thompson monoid.

Let us also discuss the case of $ax+b$-semigroups over rings of algebraic integers in number fields. This case is also discussed in detail in \cite{CunAlgAct}. Let $K$ be a number field with ring of algebraic integers $R$. We apply our K-theory result to the semigroup $P = R \rtimes R\reg$. This semigroup embeds into the $ax+b$-group $K \rtimes K\reg$. All our conditions are satisfied, so that we only need to compute orbits and stabilizers. We have a canonical identification
$$
  G \backslash \cJ_{R \rtimes R\reg \subseteq K \rtimes K\reg} \overset{\cong}{\lori} Cl_K, \, [\mfa \times \mfa\reg] \ma [\mfa].
$$
Moreover, for the stabilizer group $G_{\mfa \times \mfa\reg}$, we obtain
$$
  G_{\mfa \times \mfa\reg} = \mfa \rtimes R^*.
$$
Here, $R^*$ is the group of multiplicative units in $R$.

Hence, our K-theory formula reads in this case
$$
  \bigoplus_{[\mfa] \in Cl_K} K_*(C^*_{\lambda}(\mfa \rtimes R^*)) \overset{\cong}{\lori} K_*(C^*_{\lambda}(R \rtimes R\reg)).
$$

There is a generalization of this formula to $ax+b$-semigroups over Krull rings (see \cite{Li5}). Let us explain this, using the notation from \S~\ref{ss:Krull}.

Let $R$ be a countable Krull ring with group of multiplicative units $R^*$ and divisor class group $C(R)$. Then our K-theory formula gives
$$
  \bigoplus_{[\mfa] \in C(R)} K_*(C^*_{\lambda}(\mfa \rtimes R^*)) \overset{\cong}{\lori} K_*(C^*_{\lambda}(R \rtimes R\reg)).
$$
The reader may also consult \cite{CunAlgAct}.

Building on our discussion of graph products in \S~\ref{ss:GraphProducts} and \S~\ref{sec:GraphProducts}, we can also present a K-theory formula for graph products. 

As in \S~\ref{ss:GraphProducts} and \S~\ref{sec:GraphProducts}, let $\Gamma = (V,E)$ be a graph with vertices $V$ and edges $E$, such that two vertices in $V$ are connected by at most one edge, and no vertex is connected to itself. So we view $E$ as a subset of $V \times V$. For every $v \in V$, let $P_v$ be a submonoid of a group $G_v$. We then form the graph products
$$P \defeq \Gamma_{v \in V} P_v$$
and
$$G \defeq \Gamma_{v \in V} G_v.$$
We have a canonical embedding $P \subseteq G$. 

For every $v \in V$, choose a system $\fX_v$ of representatives for the orbits $G_v \backslash \cJ_{P_v \subseteq G_v}\reg$ which do not contain $P_v$. Moreover, for every non-empty subset $W \subseteq V$, define $\fX_W \defeq \prod_{w \in W} \fX_w$. Combining Proposition~\ref{const-ideals}, Proposition~\ref{ind_GraphProducts}, Proposition~\ref{Toeplitz_GraphProducts} and Theorem~\ref{THM:K}, we obtain
\btheo
Assume that for every vertex $v$ in $V$, our semigroup $P_v$ satisfies independence, and that $P_v \subseteq G_v$ is Toeplitz. Moreover, assume that $G$ satisfies the Baum-Connes conjecture with coefficients. Then the K-theory of the reduced C*-algebra of $P$ is given by
$$
  K_*(C^*_{\lambda}(P^*)) \oplus \bigoplus_{\substack{\emptyset \neq W \subseteq V \\ W \times W \in E}} \bigoplus_{(X_w)_w \in \fX_W} K_*(C^*_{\lambda}(\prod_{w \in W} G_{X_w})) \overset{\cong}{\lori} K_*(C^*_{\lambda}(P)).
$$
\etheo
\bproof
We know that $P$ satisfies independence by Proposition~\ref{ind_GraphProducts}, and we know that $P \subseteq G$ is Toeplitz by Proposition~\ref{Toeplitz_GraphProducts}. Moreover, it is an immediate consequence of Proposition~\ref{const-ideals} that
\bglnoz
  && G \, \backslash \, \cJ_{P \subseteq G}\reg\\ 
  &=& \gekl{P} \sqcup \menge{\eckl{\rukl{\prod_{w \in W} X_w} \cdot P}}{\emptyset \neq W \subseteq V, \, W \times W \subseteq E, \, (X_w)_w \in \fX_W}.
\eglnoz
As we get for the stabilizer groups
$$
  G_{\rukl{\prod_{w \in W} X_w} \cdot P} = \prod_{w \in W} G_{X_w},
$$
our theorem follows from Theorem~\ref{THM:K}.
\eproof
Note that the graph product $G$ satisfies the Baum-Connes conjecture with coefficients if for every vertex $v \in V$, the group $G_v$ has the Haagerup property. This is because, by \cite{AD}, the graph product $G$ has the Haagerup property in this case.

\section[Further developments, outlook, and open questions]{Further developments, outlook, \\ and open questions}

Based on the result we presented, in particular descriptions as partial or ordinary crossed products as well as our K-theory formula, we obtain classification results for semigroup C*-algebras.

For instance, the case of positive cones in countable subgroups of the real line, where these groups are equipped with the canonical total order coming from $\Rz$, have been studied in \cites{Dou,JX,CPPR,Li4}. It turns out that the semigroup C*-algebra of such positive cones remembers the semigroup completely. Actually, we can replace the semigroup C*-algebra by the ideal corresponding to the boundary quotient. It turns out that also these ideals determine the positive cones completely.

For right-angled Artin monoids, a complete classification result was obtained in \cite{ELR}, building on previous work in \cites{CrLa1,CrLa2,Iva,LR1}. The final classification result allows us to decide which right-angled Artin monoids have isomorphic semigroup C*-algebras by looking at the underlying graphs defining
our right-angled Artin monoids. The invariants of the graphs deciding the isomorphism class of the semigroup C*-algebras are explicitly given, and easy to compute in concrete examples.

For Baumslag-Solitar monoids, important structural results about their semigroup C*-algebras were obtained in \cites{Sp1,Sp2}.

In the case of $ax+b$-semigroups over rings of algebraic integers in number fields, partial classification results have been obtained in \cite{Li3}, building on previous work in \cites{CDL,EL}. It turns out that for two number fields with the same number of roots of unity, if the $ax+b$-semigroups over their rings of algebraic integers have isomorphic semigroup C*-algebras, then our number fields must have the same zeta function. In other words, they must be arithmetically equivalent (see \cites{Perlis,PS}).

In addition to these classification results, another observation is that the canonical commutative sub-C*-algebra (denoted by $D_{\lambda}(P)$) of our semigroup C*-algebra often provides interesting extra information. In many situations, the partial dynamical system attached to our semigroup (embedded into a group) is topologically free, and then this canonical commutative sub-C*-algebra is a Cartan subalgebra in the sense of \cite{R08}. For instance, for rings of algebraic integers in number fields, it is shown in \cite{Li6} that Cartan-isomorphism for two semigroup C*-algebras of the $ax+b$-semigroups implies that the number fields are arithmetically equivalent and have isomorphic class groups. This is a strictly stronger statement then just being arithmetically equivalent, as there are examples of number fields which are arithmetically equivalent but have difference class numbers (see \cite{dSP}).

It would be interesting to obtain structural results for semigroup C*-algebras of the remaining examples mentioned in \S~\ref{sec:Ex}.

For instance, for more general totally ordered groups, the semigroup C*-algebras of their positive cones have not been studied and would be interesting to investigate. Their boundary quotients are given by the reduced group C*-algebras of our totally ordered groups. It would be interesting to study the structure of the ideals corresponding to these boundary quotients.

For Artin monoids which are not right-angled, it would be interesting to find out more about their semigroup C*-algebras. For example, the case of Braid monoids would already be interesting. Here the boundary quotients are given by the reduced group C*-algebras of Braid groups. Therefore, the semigroup C*-algebras of Braid monoids cannot be nuclear. But what about the ideals corresponding to the boundary quotients?

It would also be very interesting to study the semigroup C*-algebra of the Thompson monoid. While the boundary quotient of the semigroup C*-algebra attached to the left regular representation is isomorphic to the reduced group C*-algebra of the Thompson group, the boundary quotient of the semigroup C*-algebra generated by the right regular representation is a purely infinite simple C*-algebra (see our discussion after Theorem~\ref{THM_sgp-rep->G0=trivial}, and also Corollary~\ref{Cor:bdquot-pis}). Is it nuclear?

In the case of $ax+b$-semigroups over rings of algebraic integers in number fields, is it possible to find a complete classification result for their semigroup C*-algebras? This means that we want to know when precisely two such $ax+b$-semigroups have isomorphic semigroup C*-algebras. It would be interesting to find a characterization in terms of the underlying number fields and their invariants.

Finally, it seems that not much is known about semigroup C*-algebras of finitely generated abelian cancellative semigroups. However, we remark that it is not difficult to see that all numerical semigroups have isomorphic semigroup C*-algebras. Moreover, subsemigroups of $\Zz^2$ are discussed in \cite{CunToricVar}.

Moreover, apart from the issue of classification, we would like to mention a couple of interesting further questions.

Given a semigroup $P$ which is cancellative, i.e., both left and right cancellative, we can form the semigroup C*-algebra $C^*_{\lambda}(P)$ generated by the left regular representation, and also the semigroup C*-algebra $C^*_{\rho}(P)$ generated by the right regular representation. It was observed in \cites{CEL2,Li5} that these two types of semigroup C*-algebras are completely different. However, strangely enough, they seem to share some properties. For instance, in all the examples we know, our semigroup C*-algebras $C^*_{\lambda}(P)$ and $C^*_{\rho}(P)$ have isomorphic K-theory (see \cites{CEL2,Li5}). There is even an example when this is the case, where our semigroup does not satisfy independence (see \cite{LiNor2}). Is this a general phenomenon? Do $C^*_{\lambda}(P)$ and $C^*_{\rho}(P)$ always have isomorphic K-theory? What other properties do $C^*_{\lambda}(P)$ and $C^*_{\rho}(P)$ have in common? For instance, what about nuclearity?

Looking at Theorem~\ref{THM:nuc-am}, and in particular Corollary~\ref{Cor:PinGam}, the following task seems interesting: Find a semigroup $P$ which embeds into a group, whose semigroup C*-algebra is nuclear, such that $P$ does not embed into an amenable group.

With our discussion of the Toeplitz condition in mind (see \S~\ref{sec:Toeplitz}), it would be interesting to find a semigroup which embeds into a group, for which the universal group embedding is not Toeplitz.

Finally, we remark that it would be an interesting project to try to generalize our K-theory computations to subsemigroups of groups without using the Toeplitz condition.

\bibliography{references}

\end{document}